\documentclass[12pt,a4paper]{amsart}

\usepackage{geometry}

\usepackage{amsmath, amsthm, amssymb, mathabx}
\usepackage[pdftex,bookmarks=true]{hyperref}
\usepackage{color}

\newtheorem{theorem}{Theorem}[section]
\newtheorem{proposition}{Proposition}[section]
\newtheorem{lemma}[theorem]{Lemma}
\newtheorem{corollary}[theorem]{Corollary}
\newtheorem{definition}[theorem]{Definition}

\newtheorem{observation}{Observation}[section]

\def\C{{\mathbb C}}  
\def\S{{\bf S}} 
\def\D{{\textbf{D}}} 
\def\E{{\textbf{E}}}

\def\H{{\textbf{H}}}
 
\def\Z{{\mathbb Z}} 
\def\M{{\bf M}}

\def\R{{\mathbb R}}

\def\P{{\bf P}}

\def\Q{{\bf Q}}
\def\I{{\bf I}}

\def\supp{{\rm supp}}

\def\Z{{\mathbb Z}}

\def\I{{\bf I}}
\def\J{{\bf J}}

\def\<{\left<}
\def\>{\right>}
\def\endproof{$\Box$}
\def\dist{{\rm dist}}

\title[Variational estimates for the bilinear Fourier integral]{Variational
estimates for the \\ bilinear iterated Fourier integral}

\author
[Y. Do, C. Muscalu, C. Thiele]
{Yen Do \ \ \ Camil Muscalu \ \ \ Christoph Thiele }

\address{Yen Do, Department of Mathematics,
The University of Virginia, Charlottesville, VA 22904-4137, USA}
\email{yendo@virginia.edu}

\address{Camil Muscalu, Department of Mathematics,
Cornell University, Ithaca, NY 14852-4201 , USA}
\email{camil@math.cornell.edu}

\address{Christoph Thiele, 
Mathematisches Institut, Universit\"at Bonn,
Endenicher Alle 60, D- 53115 Bonn, Germany,
formerly
Department of Mathematics,
UCLA, Los Angeles, CA 90095-1555, USA
}
\email{thiele@math.uni-bonn.de}


\thanks{Y. D. partially supported by NSF grants DMS-0635607002, DMS-1201456, and DMS-1521293.}
\thanks{C.M. partially supported by NSF grant DMS-0653519.} 
\thanks{C.Th. partially supported by NSF grant DMS 1001535 and by the Hausdorff Center for Mathematics.}
\date{\today}

\begin{document}

\begin{abstract}
We prove pointwise variational $L^p$ bounds for a bilinear Fourier integral operator in a large but not necessarily sharp range of exponents.
This result is a joint strengthening of the corresponding bounds for the 
classical Carleson operator, the bilinear Hilbert transform, the variation norm 
Carleson operator, and the bi-Carleson operator. Terry Lyon's  rough path theory allows for extension of our result to multilinear estimates. 
We consider our result a proof of concept for 
a wider array of similar estimates with possible applications to 
ordinary differential equations.

\end{abstract}

\maketitle

\section{Introduction}
Consider the bilinear iterated Fourier inversion integral
\begin{equation}\label{e.BHTdefn}
B(f_1,f_2)(x) := \int_{\xi_1<\xi_2} \widehat{f}_1(\xi_1) \widehat{f}_2(\xi_2) e^{ix(\xi_1+\xi_2)} d\xi_1d\xi_2 \ . \ \ 
\end{equation}
It is a close relative of the bilinear Hilbert transform, and as such satisfies
$L^p$ bounds as in \cite{lt1997,lt1997b,lt1998,lt1999}. 


Given any $r\in (0,\infty)$, let $T_r$ denote the following stronger operator
\begin{equation}\label{e.Tr-defn} 
\sup_{K, N_0<\dots <N_K} \Big(\sum_{j=1}^K \Big|\int_{N_{j-1}<\xi_1<\xi_2<N_j}\widehat f_1(\xi_1) \widehat f_2(\xi_2)e^{ix(\xi_1+ \xi_2)}d\xi_1d\xi_2\Big|^{r/2}\Big)^{2/r} \ \ .
\end{equation}
Thus $T_r$ is a \emph{variation sum} over truncations of $B$, in particular it dominates  both $B$ and the bi-Carleson operator considered  in 
\cite{mtt2006f}, which essentially is the limit case $r\to \infty$ of $T_r$.

The main result of our paper is the following theorem:
\begin{theorem}\label{t.main} Assume that $r>2$. Then $T_r$ is bounded from $L^{p_1} \times  L^{p_2}$ to $L^{p_3}$ provided that $1/p_3=1/p_1+1/p_2$ and
\begin{equation}\label{e.p1p2-constraint}
 \max(1, \frac {2r}{3r-4})  <p_1,p_2\le\infty \ \ , \ \  \max(\frac 23, \frac{r'}2) < p_3 < \infty \ \ .
\end{equation}
\end{theorem}

 Besides  strengthening  \cite{lt1997,lt1999,mtt2006f}, Theorem~\ref{t.main}  also implies a range of the $L^p$ estimates for the variation norm Carleson theorem  in \cite{osttw2012}. Namely, the variation norm Carleson estimate can be
obtained by a variant of \eqref{e.Tr-defn} without the constraint $\xi_1<\xi_2$,
which in turn can be estimated by the sum of \eqref{e.Tr-defn} and a symmetric version of \eqref{e.Tr-defn}.

The theory of ordinary differential equations  with rough driving signals initiated by T. Lyons \cite{lyons1998} and developed by many, for example \cite{ly2013} discusses similar 
expressions as \eqref{e.Tr-defn} and controls them using additional
information on the Fourier transform of $f_1,f_2$, allows to bootstrap our main theorem in a certain range of exponents to multi(sub)linear estimates: 
\begin{corollary}\label{t.main-multi}Let $k\ge 3$. For any $r>0$ let $T_{k,r}$ denote
\begin{equation}\label{e.Tkr-defn} 
\sup_{K, N_0<\dots <N_K} \Big(\sum_{j=1}^K \Big|\int_{N_{j-1}<\xi_1<\dots<\xi_k<N_j} \prod_{m=1}^k \widehat f_m(\xi_m)e^{ix\xi_m}d\xi_m\Big|^{r/k}\Big)^{k/r} \ \ .
\end{equation}
Then for every $p\in (\frac 32, \infty)$ and $r>\max (2,\frac p{p-1})$ it holds that
\begin{equation}\label{e.Tkr} \|T_{k,r}[f_1,\dots, f_k]\|_{p/k} \lesssim \prod_{m=1}^k \|f_m\|_p \ .
\end{equation}
\end{corollary}
The end-point version $r=\infty$ of this corollary was posed as a problem 
in \cite{mtt2015stein}.

It should be of interest to study variants of  Theorem~\ref{t.main} where 
$e^{ix(\xi_1+\xi_2)}$ is replaced by $e^{ix(\alpha_1 \xi_1 +\alpha_2\xi_2)}$ for general
real parameters $\alpha_1,\alpha_2$ and similarly for Corollary~\ref{t.main-multi}. Such variants of Theorem~\ref{t.main} do not hold if $\alpha_1\alpha_2(\alpha_1+\alpha_2)= 0$, the interesting case $\alpha_1+\alpha_2 = 0$ is discussed in \cite{mtt2003mrl}. It is possible that $\alpha_1\alpha_2(\alpha_1+\alpha_2)\neq 0$ is the only constraint towards such variants of Theorem~\ref{t.main}, this would lead to strong variants of Corollary~\ref{t.main-multi} with interesting consequences for rough ordinary differential equations. For our present purpose, 
we hope that our omission of further parameters $\alpha_1,\alpha_2$ simplifies the readability of our proof and we defer the discussion of potentially rich 
ramifications of the theory of general parameters to the future. 
We refer to \cite{kesler2015} for a discussion of a degenerate trilinear
variant of \eqref{e.BHTdefn}.

While \cite{osttw2012} establishes a sharp range of exponents for the variation norm Carleson operator, we do not prove that the range of exponents of Theorem \ref{t.main} is sharp, though this range is clearly dictated by our method of proof. Part of our ambition was to obtain a sufficiently large range of exponents to allow for instances of Corollary~\ref{t.main-multi}, for example a point $p_1=p_2=2$ and $2<r<3$. 
By the H\"older inequality, the simpler version of $T_{r}$ where the the bilinear Hilbert symbol $1_{\xi_1<\xi_2}$ is replaced by $1$ is controlled by a product of two variation norm Carleson operators, where the new variation-norm exponents $s$ and $t$ satisfy  $2/r=1/s+1/t$. In order to apply the known estimates (from \cite{osttw2012}) to these variation-norm Carleson operators,   we need $p_1>s'$ and $p_2>t'$ and which clearly leads to   the constraint $p_3>r'/2$. We also need $s,t>2$ which together with the relation $2/r=1/t+1/s$ leads to $t',s'>\max(1,2r/(3r-4))$ and from there we then obtain the  lower bound constraints for $p_1,p_2$.

It is not hard to see that for $r\ge 4$ the range \eqref{e.p1p2-constraint} becomes  the classical range for the bilinear Hilbert transform (and in particular independent of $r$). Since variation-norm operators with larger exponents are smaller, we may assume without loss of generality that $2<r<4$ in the rest of the paper.

\section{Outline of the proof}\label{e.outline}
By dualization and monotone convergence, we cound find  measurable functions $K:\R \to \Z_+$, $N_0(x)\le N_1(x)\le \dots$, and $d_0(x),d_1(x),\dots$ such that

(i)  $\sum_{j\ge 0}|d_j|^{r/(r-2)} \equiv 1$;

(ii) for every $x$, if $j>K(x)$ then $d_j(x)=0$ and $N_j(x)=N_{j-1}(x)$; and

(iii) for every $x\in \R$ we have $ |T_r(f_1,f_2)(x)|  \lesssim B_r(f_1,f_2)(x)$, where
\begin{eqnarray*}
B_r(f_1,f_2)(x) &:=& \sum_{j=1}^{K} d_j \iint _{N_{j-1}<\xi_1<\xi_2<N_j} \prod_{j=1,2} \widehat f_j(\xi_j) e^{ix \xi_j} d\xi_1d\xi_2 \ \ .
\end{eqnarray*}
Thus,  it suffices to prove the desired estimate  $B_r$, provided that the implicit constants depend only on $r$ and $p_1,p_2,p_3$. We will fix $K$, $(N_j)$, and $(d_j)$ in the rest of the paper. By  monotone convergence, we may assume that $K$, $N_j$ are bounded.

For any  $M<N$ we will  decompose  $1_{M<\xi_1<\xi_2<N}$ into three components: 
\begin{eqnarray}
\label{e.triangular-symbol-decomp}  
\nonumber && m_{CC}(M,N,\xi_1,\xi_2) +  m_{BC}(M,N,\xi_1,\xi_2)  + m_{LM}(M,N,\xi_1,\xi_2)  
\end{eqnarray}
\begin{itemize}
\item $m_{BC}$ captures the singularity along the line segment from $(M,M)$ to $(N,N)$,
\item $m_{CC}$ captures the singularity along the other two edges of the triangle, and
\item $m_{LM}$ is an error term that has two singularities at $(M,M)$ and $(N,N)$. 
\end{itemize}
Construction of these symbols are detailed in  Section~\ref{s.triangular-symbol-decomp} using a hybrid of the arguments in \cite{osttw2012} and \cite{mtt2006f}.  Applying this decomposition for $(M,N)=(N_{j-1}, N_j)$ and using the triangle inequality,  it follows that $B_r$ is controlled by three corresponding bilinear operators. Via standard arguments (detailed in the Appendices), each of these three operators is in turn controlled by a bounded sum of discrete operators, which will be described in Section~\ref{s.model}. To prove boundedness of the discrete model operators in the desired $L^p$ ranges, we  will use the new $L^p$ theory for outer measures introduced in \cite{dt2014}.  It turns out that    analogues (for outer measure spaces) of classical singular integral operators arise naturally in our proof, and they are effective tools to handle nested levels of time-frequency analysis. In order to study these operators, we adapted an argument in \cite{dt2014} to prove a Marcinkiewicz interpolation theorem (see Lemma~\ref{l.outer-marcinkiewicz}) for (quasi)sublinear maps between outer measure spaces, which generalizes a simpler interpolation result in \cite{dt2014}. 

\subsection{Notational conventions}

By a (standard) dyadic interval we mean $[n2^k, (n+1)2^k)$ for some $n,k\in \mathbb Z$. For any $a\in [0,1)$, by an $a$-shifted dyadic interval we mean an interval of the form $2^k([n,n+1) + (-1)^k a)$. 

For any interval $I$ we denote its midpoint by $c(I)$, its left children by $I_-$ and its right children by $I_+$. For any positive $C$, the $C$-enlargement of $I$ is defined to be the $C$-dilation of $I$ from its center and will be denoted by $CI$. Enlargements of cubes are defined similarly. We will define $\widetilde \chi_I(x)=(1+\frac{|x-c(I)|}{|I|})^{-2}$.

For any set $A\subset \R$, we define $D_y A:=\{2^y \xi: \xi\in A\}$ the dilation of $A$ relative to the origin. For any number $\alpha$ and any interval $I$, let  $I+\alpha := \{x+\alpha: x\in I\}$.

We say that $A\lesssim B$ if there is an absolute  constant $C \in (0,\infty)$ such that $|A|\le C|B|$. If $C$ depends on $t_1,t_2,\dots$ we will say that $A\lesssim_{t_1,t_2,\dots} B$.  (We sometimes suppress some subscripts if the dependence is not important for the relevant discussion.)

If $J$ is an interval we define $M_Jf = \sup_{I \supset J} |I|^{-1} \int_I |f(x)|dx$  for every function $f$.
We will use the following normalizations for  inner product and Fourier transforms:
\begin{eqnarray*}
\<f,g\> := \int_{\R} f(x) \overline{g(x)}dx \ \ , \ \ \widehat f(\xi) :=  (2\pi)^{-1/2}\int_{\R} f(x)e^{-ix\xi}dx \ \ .
\end{eqnarray*}

\section{Decomposition of the triangular symbol}\label{s.triangular-symbol-decomp}
In this section, we will  construct  $m_{CC}$ and $m_{BC}$  ($m_{LM}$ is defined using \eqref{e.triangular-symbol-decomp}).

\subsection{Construction of $m_{CC}$}

\subsubsection{Decomposition of $1_{M<\xi<N}$} \label{s.interval-decomp} Below we adapt a decomposition in \cite{osttw2012}. We partition $(M,N)$ into maximal dyadic intervals $I$ such that
\begin{eqnarray*}
\dist(I, \{M,N\}) &\ge& L_1 |I| \ \ .
\end{eqnarray*}
Let $\mathcal H=\mathcal H(M,N)$ denote the set of these dyadic intervals. It is clear that the length of two neighboring elements of $\mathcal H$ differs by a factor of $2^k$ where $k$ is integer in $[-1,1]$, therefore we may denote the  lengths of the left and right neighbors of $I \in \mathcal H$ by $u(I)|I|$ and $v(I)|I|$ where $u(I), v(I) \in \{ \frac 1 2,1,2\}$ and these two numbers are completely determined from the  following details:
\begin{itemize}
\item the unique  $m\in \Z_+$  such that $M\in I-(m+1)|I|$;
\item the unique   $n\in \Z_+$ such that $N\in I+(n+1)|I|$;
\item whether $I$ is the left or right children of its dyadic parent.
\end{itemize}
Let $A$ denote the set of eligible $(side,m,n)$ for $I\in \mathcal H$, where $side\in\{left,right\}$. For any $\alpha \in A$,  let $\I_{M,N}(\alpha)$ be the set of  corresponding dyadic intervals.

Given any $1/2<c<5/8$ by elementary arguments we could construct  smooth functions $\phi_{u,v}$ indexed by $(u,v) \in \{\frac 1 2, 1, 2\}^2$, all supported in $[-c,c]$, such that
\begin{eqnarray*}
1_{(M,N)}(\xi) &=& \sum_{I\in \mathcal H(M,N)} \phi_{u(I),v(I)}(\frac{\xi-c(I)}{|I|})     \ \ .
\end{eqnarray*}
We will write $\phi_{\alpha,I}(\xi) := \phi_{u,v}(\frac{\xi-c(I)}{|I|})$ if the details of $I$ match with $\alpha$, thus
\begin{eqnarray}\label{e.interval-decomp}
1_{(M,N)} (\xi) &=& \sum_{\alpha \in A} \sum_{I\in \I_{M,N}(\alpha)}  \phi_{\alpha,I}(\xi)  \ \ .
\end{eqnarray}
Note that the sum in the right  converges absolutely pointwise: every term is nonegative, and at every $\xi$ there are only finitely many nonzero terms in the sum.

\subsubsection{Definition of $m_{CC}$} 
Intuitively, $m_{CC}(M,N,\xi_1,\xi_2)$  is a smooth restriction of $1_{M<\xi_1<\xi_2<N}$ to 
\begin{equation}\label{e.R1-defn}
\Big\{(\xi_1,\xi_2) \in [M,N]^2: \min(|\xi_1-M|, |\xi_2-N|) \le    |\xi_1-\xi_2|/{200} \Big\} \ \ ,
\end{equation}
which will be denoted by $R_1(M,N)$. For a more precise statement, see Lemma~\ref{l.CC-properties}.

To motivate, note that using \eqref{e.interval-decomp} we obtain for any $(\xi_1,\xi_2)\in R_1$ 
\begin{eqnarray}
\nonumber 1_{M<\xi_1<\xi_2<N}  &=& 1_{M<\xi_1<N} 1_{M<\xi_2<N}\\
\label{e.product} &=& \sum_{\alpha,\beta\in A} \sum_{I\in \I_{M,N}(\alpha)} \sum_{J\in \I_{M,N}(\beta)} \phi_{\alpha,I}(\xi_1) \phi_{\beta,J}(\xi_2) 
\end{eqnarray}
and we will construct $m_{CC}$ by removing terms supported far from $R_1$ in the sum.

Specifically, let
\begin{eqnarray*}
A_1  = \{m_\alpha \le n_\alpha /4\} \ \ , \ \ 
A_3 = \{m_\alpha \ge 4 n_\alpha\} \ \ , \ \
A_2 = A-A_1-A_3 \ \ ,
\end{eqnarray*}
which are subsets of $A$. As we will see,  $A_2$ is a finite set. 
\begin{definition}\label{d.CCsym} Define
\begin{equation*}
\boxed{m_{CC}(M,N,\xi_1,\xi_2)  := \sum_{k=1}^5 m_{CC,k}(M,N,\xi_1,\xi_2)}
\end{equation*}
where $m_{CC,j}$ are sub-sums of the right hand side of \eqref{e.product} under some extra constraints.  We will always have $I\in \I_{M,N}(\alpha)$ and $J \in \I_{M,N}(\beta)$ and the summands are always $\phi_{\alpha,I}(\xi_1)\phi_{\beta,J}(\xi_2)$, and the following table details the extra contraints:
\begin{table}[h]
\begin{tabular}{|c|c|c|c|}
\hline
Symbols & Conditions on $\alpha$ & Conditions on $\beta$ & Extra conditions on $I$,$J$    \\ \hline\hline
$m_{CC,1}$ & $\alpha\in A_1$ & $\beta\in A_1$  & $|I|\le |J|/{16}$    \\ \hline \hline
$m_{CC,2}$ & $\alpha\in A_1$ & $\beta\in A_2$  & $|I|\le |J|/16$    \\ \hline \hline
 $m_{CC,3}$ & $\alpha\in A_1$ & $\beta\in A_3$  & None    \\ \hline \hline
 $m_{CC,4}$ & $\alpha\in A_2$ & $\beta\in A_3$  & $|I| \ge 16|J|$    \\ \hline \hline
 $m_{CC,5}$ & $\alpha\in A_3$ & $\beta\in A_3$  & $|I| \ge 16|J|$    \\ \hline
\end{tabular}
\end{table}  
\end{definition}

\subsection{Construction of $m_{BC}$}

In this section, we construct the $m_{BC}$ symbol, which may be viewed as a smooth restriction of $\chi_{M<\xi_1<\xi_2<N}$
to a neighborhood of $R_2(M,N)$, which consists of all $M<\xi_1,\xi_2<N$ such that
\begin{equation}\label{e.R2-defn} 
 |\xi_1-\xi_2|< \frac { \min(|\xi_1+  \xi_2-M|, |\xi_1+  \xi_2-N|)}{100}  
\end{equation} 
For a more precise statement, see Lemma~\ref{l.BC-properties}. 

To motivate the construction, we first note that if $(\xi_1,\xi_2) \in R_2$ then
$$\chi_{M<\xi_1<\xi_2<N} = \chi_{\xi_1<\xi_2} \chi_{M< \xi_1+ \xi_2 < N} \ \ .$$

We will construct $m_{BC}$  by writing the product in the right hand side of the last display as a sum of products of wave packets (using a suitable decomposition for each factor), and then removing from this sum essentially those  terms that are supported far from $R_2$. 

\subsubsection{Decomposition of $\chi_{\xi_1<\xi_2}$}\label{s.mbc-def}
We largely follow \cite{lt1999} (see also \cite{thiele2006, ms2013}). 

\underline{Shifted cubes:} For any $b=(b_1,\dots, b_n) \in [0,1)^n$ we say that $S=S_1\times \dots \times S_n$ is a $b$-shifted dyadic cube if $S_j$ is a $b_j$-shifted dyadic interval, and $|S_1|=\dots =|S_n|$. In this paper, unless otherwise specified, the coordinates of  underlying shifts are assumed to be in  $\{0,\frac 1 3, \frac 23\}$.

\underline{Whitney decomposition:} For each $b\in \{0, \frac  13, \frac 2 3\}^2$ consider the collection $\S_b$ of all maximal $b$-shifted squares $S_1\times S_2$ satisfying

(i)  $\dist(S_1\times S_2,\ell)  \ge L_2|S_1|$, here $\ell$  is the line $\{\xi_1=\xi_2\}$.

It is clear that for such square it holds that
\begin{equation}\label{e.whitney}
L_2|S_1| \le  \dist(S_1\times S_2,\ell)  \le (2 L_2+1)|S_1| \ \ .
\end{equation}
 Using a partition of unity argument, we may write
$$\chi_{\xi_1<\xi_2} =  \sum_{b} \sum_{S_1\times S_2 \in \S_b} \phi_{S_1\times S_2}(\xi_1,\xi_2)$$
where $\{\phi_{S_1\times S_2}\}$ is a family of nonnegative $C^\infty$ bump functions, such that $\phi_{S_1\times S_2}$ is supported inside $\frac 4 {5} S_1 \times \frac 45 S_2$.  Note that if $(\xi_1,\xi_2) \in S$ then $\xi_1+ \xi_2  \in   S_1+  S_2$. Note that $ S_1+   S_2$ can be generously covered by $6$ intervals of the form $(3/4)I_j$, $1\le j \le 4$, where $I_j$ are  shifted dyadic intervals having the same length as $|S_1|$.  


 Using partitions of unity, we may find nonnegative smooth functions $\phi_{3,I_j}$ with support  inside $(4/5)I_j$, such that for every  $\xi \in  S_1+  S_2$ it holds that $1=\sum_{j} \phi_{3,I_j}(\xi)$.
We obtain
$$\chi_{\xi_1<\xi_2} =\sum_{S \in \S} \phi_{S_1\times S_2}(\xi_1,\xi_2)\phi_{3,S_3}(\xi_1+ \xi_2)$$
where $\S$ denote the collection of cubes $S_1\times S_2\times S_3$ formed using $S_1\times S_2\in \bigcup \S_b$ and $S_3$ are the covering intervals $I_j$'s discussed above. Note that every element of $\S$ is a shifted dyadic cube, where the underlying shifts are elements of $\{0,\frac 1 3 , \frac 2 3 \}^3$. Expanding   $\phi_{S_1\times S_2}(\xi_1,\xi_2)$ into bilinear Fourier series, we obtain
\begin{equation}\label{e.bht-decomp}
\chi_{\xi_1<\xi_2}(\xi_1,\xi_2) = \sum_{k\in \Z^2} a_k \sum_{S\in \S} \phi_{1,S,k}(\xi_1)\phi_{2,S,k}(\xi_2) \phi_{3,S,k}(\xi_1+ \xi_2) \ \ ,
\end{equation}
here $(a_k)$ is a rapidly decaying sequence and $\phi_{i,S,k}$'s are $C^d$-bump functions uniformly adapted to $S\in \S$, and $\phi_{i,S,k}$ is  supported inside $\frac 5 6 S_i$. (In fact $\phi_{3,S,k} \equiv \phi_{3,S_3}$ is independent of $k$, but we prefer to use $\phi_{3,S,k}$   for later convenience of notation.) Als, $d$ is a fixed finite constant, but could be chosen arbitrarily large.

\subsubsection{Definition of $m_{BC}$}

  For convenience, let $\ell(S)=|S_1|=|S_2|=|S_3|$ denote the side length of $S$. Using \eqref{e.interval-decomp} and \eqref{e.bht-decomp}, it follows that we may write $\chi_{\xi_1<\xi_2}(\xi_1,\xi_2) \chi_{M< \xi_1+ \xi_2< N}$ as
\begin{equation}\label{e.expansionBC}
 =\sum_{k,\alpha, S, I} a_k \phi_{\alpha,I}( \xi_1+\xi_2)  \phi_{3,S,k}(\xi_1+  \xi_2)\prod_{j=1}^2 \phi_{j,S,k}(\xi_j)   \ \ . \end{equation}
Here the summation is over all $k\in \mathbb Z^2$, $\alpha\in A$, $S\in \S$, and $I\in \I_{M,N}(\alpha)$.
\begin{definition} \label{d.BCsym} Define
\begin{eqnarray*}
\boxed{m_{BC}(M,N,\xi_1,\xi_2) := \sum_{\stackrel{k,\alpha, S,I}{\ell(S)\le  |I|}} a_k      \phi_{\alpha,I}(\frac{\xi_1+\xi_2}2)    
\phi_{3,S,k}(\xi_1+\xi_2)
\prod_{j=1}^2 \phi_{j,S,k}(\xi_j)}
\end{eqnarray*}
\end{definition}

\section{The model operators}\label{s.model}

\subsection{Tiles and wave packets} 

A (standard) tile $P=I_P\times \omega_P$ is a  rectangle of area $1$, where the spatial interval $I_P$ is a standard dyadic interval and the frequency interval $\omega_P$ is a standard dyadic interval. We will also use shifted tiles, where $I_P$ is still standard dyadic but $\omega_P$ is a shifted dyadic interval.

Let $\P$ be a tile collection. For $1\le p\le \infty$ we say that the collection of functions $\{\phi_P, P\in \P\}$ is a $L^p$-normalized wave packet collection   if $\widehat \phi_P  \subset (5/4) \omega_P$, and 
\begin{eqnarray*}
\frac{d^n}{dx^n} 
\Big(e^{-ic(\omega_P)x}\phi_P(x)\Big) &\lesssim_{n}&  \frac{1}{|I_P|^{n+1/p}} (1+ \frac{|x-c(I_P)|}{|I_P|})^{-n}
\end{eqnarray*}
uniformly over $P\in \P$, for all $n\ge 0$.  If the   estimate holds only for $0\le n \le n_0$ then we say that the collection is of order $n_0$.

\subsection{Multi-tiles}
We say $Q=(Q_1,\dots, Q_m)$ is an $m$-tile if the tiles $Q_1,\dots, Q_m$ share the same spatial interval,  denoted by $I_Q$.  The cube $\omega_Q:=\omega_{Q_1}\times \dots \times \omega_{Q_m}$ is called the frequency cube of $Q$.

\begin{definition}[Sparse]\label{d.sparse}  A collection $\Q$ of $m$-tiles is sparse relative to a constant $C_0$ if the following hold: for any $Q, R \in \Q$ with   $|I_R|/C_0 \le |I_Q|\le |I_R|$ we must have $|I_Q|=|I_P|$, and furthermore either $\omega_Q=\omega_R$ or $C_0\omega_Q \cap C_0 \omega_{R}  =\emptyset$.
\end{definition}

\begin{definition}[Rank-$1$]\label{d.rank1} A collection $\Q$ of $m$-tiles is of rank $1$ relative to $C_1\ge 1$ if the following holds for any $Q,R \in \Q$:
\begin{itemize}
\item  If there exists $j$ such that $\omega_{Q_j}=\omega_{R_j}$ then $\omega_Q=\omega_R$.
\item  For every $j_0 \in \{1,\dots, m\}$, if  $5\omega_{Q_{j_0}} \subset 5\omega_{R_{j_0}}$ then  $5C_1\omega_{Q}\subset 5C_1 \omega_{R}$. If furthermore $|I_R|<|I_Q|$ then $5\omega_{Q_j}  \cap 5\omega_{R_j} = \emptyset$ for every $j\ne j_0$.
\end{itemize}
\end{definition}

 \subsection{Rigid triples of intervals}\label{s.rigid}
We say that $\{(I_{lower}, I, I_{upper}), I \in \I\}$ is a rigid collection of interval triples if there are integers $L_1 < m,n \lesssim L_1$ such that one of the following situations happens:
\begin{itemize}
\item[(i)] For every $I\in \I$ we have $I_{lower} = I - m|I|$ and $I_{upper} = I+n|I|$.\\
\item[(ii)] For every $I\in \I$ we have $I_{lower} = I - m|I|$ and $I_{upper} = [c(I)+(n-1/2)|I|,\infty)$.\\
\item[(iii)] For every $I\in \I$ we have $I_{lower} = (-\infty, c(I)- (m-1/2)|I|)$ and $I_{upper} = I+n|I|$.
\end{itemize}
We say that two collections have the same structure if the same situation (i.e. (i) or (ii) or (iii)) holds for both, with possibly different pairs $(m,n)$.

\subsection{Description of the model operators}
To bound $T_r$, we will show that it suffices to bound the following four types of operators: $T_{C\times C}$ (product of Carleson operators), $T_{CC}$ (paraproduct of Carleson operators), $T_{BC}$ (composition of bilinear Hilbert transform and Carleson operators), $T_{LM}$ (variational bi-linear Carleson operators). We will define these operators shortly. For convenience of notation, in the following  we denote 
\begin{eqnarray}\label{e.modpacket}
\widetilde\phi_{1,P,j}(x) = \phi_{1,P}(x)1_{N_{j-1}(x) \in \omega_{1,P,lower}} 1_{N_j (x)\in \omega_{1,P,upper}}
\end{eqnarray}
and we define $\widetilde \phi_{2,P,j}$, $\widetilde \phi_{3,P,j}$, $\widetilde \phi_{1,Q,j}$, \dots similarly. As a convention, $\P$s will denote tile collections and $\Q$s will denote shifted tri-tile collections, and all underlying collection of wave packets are $L^1$-normalized. All interval triples will be rigid, and in type $CC$ we demand that the two underlying rigidity types are the same. Without loss of generality we assume that all tile and tritiles collections are finite and sufficiently sparse (all estimates are uniform over these collections), and for $\Q$ we assume that the tri-tiles share the same shift.

\begin{eqnarray*}
T_{C\times C}(f_1,f_2) &=&\sum_{j=1}^{K}  d_j \Big(\sum_{P \in \P_1} |I_P|\<f_1,\phi_{1,P}\> \widetilde \phi_{1,P,j} \Big)  \Big(\sum_{P\in \P_2} |I_{P}|\<f_2,\phi_{2, P}\> \widetilde \phi_{2, P,j}  \Big) \\
T_{CC}(f_1,f_2) &=&\sum_{j=1}^K d_j  \sum_{P \in \P_1} |I_P| \<f_1,\phi_{1,P}\> \widetilde \phi_{1,P,j}   \sum_{\stackrel{P' \in  \P_2}{|I_{P'}|\le |I_P|/16}} |I_{P'}|\<f_2,\phi_{2,P'}\>\widetilde \phi_{2,P',j}   \\
T_{BC}(f_1,f_2)&=&\sum_{j=1}^K d_j \sum_{Q\in \Q} |I_Q|\<f_1,\phi_{1,Q}\>\<f_2,\phi_{2,Q}\>   \sum_{P\in \P: |I_P|\le |I_Q|} |I_P|\<\phi_{3,Q},\phi_P\>\widetilde \phi_{P,j}  \\
T_{LM}(f_1,f_2)&=&\sum_{j=1}^K  d_j  \sum_{Q \in \Q} |I_Q| \<f_1,\phi_{1,Q}\> \<f_2,\phi_{2,Q}\>  \phi_{3,Q}(x)    1_{N_{j-1}, N_j  \ constraints} \ \ .
\end{eqnarray*}
The contraints in $T_{LM}$  read as follow: 
\begin{itemize}
\item $2N_{j-1}\in \omega_{3,Q,lower}$, $2N_j \in \omega_{3,Q,upper}$, 
\item $N_{j-1}\in \omega_{k,Q,lower}$ and  $N_j \in \omega_{k,Q,upper}$ for each $k=1,2$. 
\end{itemize}
Furthermore, in $T_{LM}$ $\Q$ is not shifted and it satisfies the following  rigidity constraint:  for $m_2,m_3$ fixed bounded  integers it holds for every $Q\in \Q$ that
\begin{itemize}
\item $\omega_{Q_2}$ is the translation of $\omega_{Q_1}$ by $m_2|\omega_{Q_1}|$,  and 
\item  $\{x/2: x\in \omega_{Q_3}\}$ is the translation of the left children of $\omega_{Q_1}$ by $m_3|\omega_{Q_1}|/2$,
\end{itemize}

\subsection{Reduction to model operators}

In the Appendix, we will show that $B_r(f_1,f_2)$ is bounded by  a finite average of discrete operators of the above types.  By the H\"older inequality (see also the discussion at the end of the introduction), the desired bounds for $T_{C\times C}$ follow from known $L^p$ estimates for discrete variation-norm Carleson operators \cite{osttw2012}: in Section~\ref{s.carl-embed} we will also deduce these estimates (see Theorem~\ref{t.varCarl}) as a byproduct of several generalized Carleson embedding estimates and the $L^p$ theory for outer measure introduced in \cite{dt2014}.  The proof of Theorem~\ref{t.varCarl} will also serve as a model for the unfortunately more technical treatments for $T_{CC}$, $T_{BC}$, and $T_{LM}$.

\section{Some background on outer measure spaces}

We recall several notions from \cite{dt2014}, with some simplifications for the setting of the current paper. An outer measure space $(X,S,\mu)$ consists of:

 (i) A countable set $X$. Often we will assume $X$ is finite,  in that case the underlying estimates are independent of the size of $X$. 
 
 (ii) An outer measure $\mu$ generated using countable coverings from a pre-measure  on a fixed collection $\E$ of non-empty subsets of $X$, which in particular covers $X$. 
 
 (iii) A \emph{size} $S$ which assigns a number in $[0,\infty]$ to each pair $(f,E)$ where $f:X \to \mathbb C$ Borel measurable and $E\in  \E$, such that
\begin{eqnarray*}
S(f+g)(E) \quad \lesssim\quad S(f)(E)+S(g)(E)  \quad, \quad S(\lambda f)(E) \quad=\quad |\lambda|S(f)(E)  \ \ ,
\end{eqnarray*}
and if $|f|\le |g|$ then $S(f)(E) \le S(g)(E)$ for all $E\in \E$.

Given an outer measure space,  we may define $\mathcal L^\infty \equiv \mathcal L^{\infty,\infty}$ by
$$\|f\|_{\mathcal L^\infty(X,S,\mu)} = outsup_X S(f) \equiv \sup_{E\in \E} S(f)(E)  \ \ ,$$
and  $\mathcal L^p$  may be defined as follows: 
\begin{itemize}
\item for Borel measurable $F\subset X$ we define $outsup_F S(f) :=\|f1_F\|_{\mathcal L^\infty(X,S,\mu)}$.
\item for any $\lambda \in \R$ let $\mu(S(f)>\lambda):=\inf\{\mu(F):  \ outsup_{X\setminus F} S(f) \le \lambda\}$. Then
\begin{eqnarray*}
\|f\|_{\mathcal L^{p,\infty}(X,S,\mu)} &=& \sup_{\lambda>0} \lambda \mu(S(f)>\lambda)^{1/p} \ \ , \\
\|f\|_{\mathcal L^{p}(X,S,\mu)} &=& \Big(p \int_0^\infty \lambda^{p-1} \mu(S(f)>\lambda)d\lambda\Big)^{1/p} \ \ .
\end{eqnarray*}
\end{itemize}
Many standard properties of classical $L^p$ spaces can be proved for outer $L^p$ spaces, see \cite{dt2014} for details. We summarize several estimates from \cite{dt2014}.

\begin{proposition}[Outer Radon--Nikodym]
Assume  that $\|f\|_{\mathcal L^\infty(X,S,\mu)}<\infty$, and for some Borel measure $\nu$ on $X$ it holds for every $E\in \E$ that
$$\int_E |f|d\nu \le C_1 \mu(E) S(f)(E) \ \ .$$
Then it holds that (the implicit constant  does not depend on  $\|f\|_{\mathcal L^\infty(X,S,\mu)}$)
$$\int_X |f|d\nu \lesssim_{C_1} \|f\|_{\mathcal L^1(X,S,\mu)}\ \ .$$

\end{proposition}

\begin{proposition}[Outer H\"older]
Suppose that  
$S(f_1f_2)(E) \le \prod_j S_j(f_j)(E)$ for all $E \in \E$.
Let $p_1,p_2,p_3\in (0,\infty]$ such that $1/p_1+1/p_2=1/p_3$. Then 
$$\|f_1f_2\|_{\mathcal L^{p_3}(X,S,\mu)} \le 2 \prod_j \|f_j\|_{\mathcal L^{p_j}(X,S_j,\mu)} \ \ .$$
\end{proposition}

\begin{proposition}[Convexity]
If $p_1<p<p_2$ and $1/p = \alpha_1/p_1 + \alpha_2/p_2$ with $\alpha_1,\alpha_2\in (0,1)$ and $\alpha_1+\alpha_2=1$, then
$$\|f\|_{\mathcal L^p(X,S\mu)} \le C  \prod_j \|f\|_{L^{p_j,\infty}(X,S,\mu)}^{\alpha_j}  \ \ .$$
\end{proposition}

The following Lemma  generalizes \cite[Proposition 3.5]{dt2014}. Below $(X,S,\mu)$ and $(Y,S',\nu)$ are given outer measure spaces.

\begin{lemma}\label{l.outer-marcinkiewicz} 
Let $K$ be an operator mapping Borel measurable functions on $X$ to Borel measurable functions on $Y$,  with the following properties:
\begin{itemize}
\item[(i)] Scaling invariance: for any $\lambda\ge 0$, $|K(\lambda f)|=|\lambda K(f)|$;
\item[(ii)] Quasi sublinear: $|K(f+g)| \le C |K(f)| + C|K(g)|$;
\item[(iii)] Bounded from $L^{p_j}(X,S,\mu)$ to $L^{q_j,\infty}(Y,S',\nu)$, i.e.  $\|K f\|_{\mathcal L^{q_j,\infty}}  \le M_j \|f\|_{\mathcal L^{p_j}}$.
\end{itemize}
Let $p_0<p_1$ and $q_0<q_1$ such that $0<p_j\le q_j\le\infty$. For $\theta \in (0,1)$ assume that
$$\frac 1 {p} = \frac{\theta}{p_1} + \frac{1-\theta}{p_0} \ \ , \ \  \frac 1 {q} = \frac{\theta}{q_1} + \frac{1-\theta}{q_0} \ \ . $$
Then  $K$ maps $L^{p}(X,S,\mu)$ into $L^{q}(Y, S', \nu)$  with norm controlled by $M_0^{1-\theta} M_1^{\theta}$.
\end{lemma}
\proof 
By scaling invariance, we may assume $M_0=M_1=M=1$. 
Let $$\alpha :=  (\frac{q}{p})  \frac{1/q_0 - 1/q_1}{1/p_0-1/p_1}$$

\underline{\bf Case I: $q_1<\infty$.}   It follows that $p_1 <\infty$ since $p_j\le q_j$.    Without loss of generality, assume that $\|f\|_{\mathcal L^{p}(X,S,\mu)} =  1$.

For each $\lambda>0$, by definition there exists a set $U$ such that $S(f1_{U^c}) \le \lambda^\alpha$ and $\mu(U) \le 2 \mu(S(f)>\lambda^\alpha)$. Let  $f_{\lambda,l} = f1_{U}$ and$f_{\lambda,s} = f 1_{U^c}$. Here $s$ stands for small and $l$ stands for large.

Since $f_{\lambda,l}$ is supported on $U$ and $|f_{\lambda,l}|\le |f|$ and monotonicity of size, we have
\begin{eqnarray*}
\|f_{\lambda,l}\|_{\mathcal L^{p_0}(X,S,\mu)}^{p_0}  &\lesssim&  \Big(\int_0^{\lambda^\alpha}  + \int_{\lambda^\alpha}^\infty\Big) t^{p_0-1}\mu(S(f_{\lambda,l})>t) dt \\
&\lesssim&    \lambda^{p_0 \alpha} \mu(S(f)>\lambda^\alpha) + \int_{\lambda^\alpha}^\infty t^{p_0-1} \mu(S(f)>t) dt \ \ .
\end{eqnarray*}
Using $\|f\|_{\mathcal L^{p}}=1$, it follows in particular that  $\|f_{\lambda,l}\|_{\mathcal L^{p_0}} \lesssim  \lambda^{1-q/q_0}$, and
\begin{eqnarray*}
\int_0^\infty \lambda^{(p-p_0)\alpha-1} \|f_{\lambda,l}\|_{\mathcal L^{p_0}(X,S,\mu)}^{p_0} d\lambda &\lesssim& \int_0^\infty \beta^{p-1} \mu(S(f)>\beta)d\beta \quad \lesssim \quad 1 \ \ .
\end{eqnarray*}
For $\|f_{\lambda,s}\|_{p_1}$, using $S(f_{\lambda,s}) \le \lambda^\alpha$ and  monotonicity of size we similarly obtain
\begin{eqnarray*}
\|f_{\lambda,s}\|_{\mathcal L^{p_1}(X,S,\mu)}   &\lesssim&   (\int_0^{\lambda^\alpha} t^{p_1-1}\mu(S(f)>t) dt)^{1/p_1}  \quad \lesssim \quad \lambda^{1-q/q_1}\ \ ,
\end{eqnarray*}
\begin{eqnarray*}
\int_0^\infty \lambda^{(p-p_1)\alpha-1}\|f_{\lambda,s}\|_{\mathcal L^{p_1}(X,S,\mu)}^{p_1} d\lambda &\lesssim& \int_0^\infty \beta^{p-1} \mu(S(f)>\beta)d\beta \quad \lesssim \quad 1 \ \ .
\end{eqnarray*}

Since $q_0 \ge p_0$ and $q_1\ge p_1$, using quasi linearity and scaling invariance and the given assumption on bounds for $K$ at the endpoints, it follows that
\begin{eqnarray*}
\nu(S'(K(f))> \lambda)  &\lesssim& \lambda^{-q_0} \|f_{\lambda,l}\|_{\mathcal L^{p_0}(X,S,\mu)}^{q_0} +  \lambda^{-q_1} \|f_{\lambda,s}\|_{\mathcal L^{p_1}(X,S,\mu)}^{q_1} \\
&\lesssim&  \lambda^{(p - p_0)\alpha - q}  \|f_{\lambda,l}\|_{\mathcal L^{p_0} (X,S,\mu)}^{p_0} + \lambda^{(p - p_1)\alpha - q}  \|f_{\lambda,s}\|_{\mathcal L^{p_1} (X,S,\mu)}^{p_1}  \ \ ,
\end{eqnarray*}
here we have used the definition of $\alpha$.  It follows that
\begin{eqnarray*}
&&\|Kf\|_{\mathcal L^{q}(Y,S',\nu)}^{q}  \quad \lesssim \quad \int_0^\infty \lambda^{q-1} \nu(S'(K(f))> \lambda)  d\lambda \\
&\lesssim& \int_0^\infty \lambda^{(p-p_0)\alpha-1} \|f_{\lambda,l}\|_{\mathcal L^{p_0}(X,S,\mu)}^{p_0} d\lambda + \int_0^\infty \lambda^{(p-p_1)\alpha-1}\|f_{\lambda,s}\|_{\mathcal L^{p_1}(X,S,\mu)}^{p_1} d\lambda
\quad \lesssim \quad  1\ \ .
\end{eqnarray*}

\underline{\bf Case I: $q_1=\infty$.}  As before, we decompose $f=f_{\lambda,l}+f_{\lambda,s}$ where $S(f_{\lambda,s}) \le c\lambda^\alpha$ for $c>0$ small (chosen later),  and $f_{\lambda,l}$ is supported in a set $U$  such that $\mu(U)\le 2\mu(S(f)>c\lambda^\alpha)$. Since $M_1=1$, by choosing $c>0$ sufficiently small we obtain $\nu(S'(Kf)>\lambda) 
\le \nu(S'(Kf_{\lambda,l}) \gtrsim \lambda^\alpha)$.  The rest of the proof is similar. 
\endproof

The following Lemma is a multilinear extension of Lemma~\ref{l.outer-marcinkiewicz}. Below $(X,S,\mu)$ and $(X_j,S_j,\mu_j)$ are outer measure spaces, $j=1,\dots, n$. Let $0<s_j<t_j\le \infty$, $j=1,\dots,n$. Let $A$ be the $n$-dimensional rectangle $\{(x_1,\dots,x_n): 1/t_j \le x_j \le 1/s_j\}$.

\begin{lemma}\label{l.outer-multi-interpolation} Let  $K$ maps measurable function $(f_1,\dots,f_n)$ on $X_1\times \dots \times X_n$ to  measurable functions on $X$. Assume that $K$ has the following properties:
\begin{itemize}
\item[(i)] Scaling invariance: for any $\lambda\ge 0$ and  $1\le j \le n$, 
$$|K(\dots \lambda f_j  \dots)|=|\lambda K(\dots f_j \dots)|$$
\item[(ii)] Quasi sublinear:  for any $1\le j\le n$ 
$$|K(f_1,\dots f_j+g_j \dots, f_n)| \le C |K(f_1,\dots f_j \dots, f_n)| + C|K(f_1,\dots g_j  \dots f_n)|$$
\item[(iii)] For every $(p,p_1,\dots, p_n)$ such that $1/p=1/p_ 1+\dots + 1/p_n$ and $(1/p_1,\dots,1/p_n)$ is one of the vertices of $A$ it holds that
$$\|K f\|_{\mathcal L^{p,\infty}}  \lesssim_{p,p_1,\dots, p_n} \prod \|f_j\|_{\mathcal L^{p_j}}$$ 
\end{itemize}
Then for every $(p,p_1,\dots, p_n)$ such that $1/p=\sum1/p_j$ and $(1/p_1,\dots,1/p_n)$ is in the interior of $A$ it holds that
$$\|K f\|_{\mathcal L^{p}}  \lesssim_{p,p_1,\dots, p_n} \prod \|f_j\|_{\mathcal L^{p_j}} \ \ .$$
\end{lemma}
\proof For simplicity we will show the proof for $n=2$, the general case is similar.

Let $(1/p_1,1/p_2)$ be in the interior of $A$ and $f_j\in \mathcal L^{p_j}$ for $j=1,2$, and assume $1/p=1/p_1+1/p_2$. We  normalize $\|f_j\|_{p_j}=1$ and we will  show that $\|K(f_1,f_2)\|_p \lesssim 1$.

\underline{Case 1: $t_1,t_2<\infty$.}

For every $\lambda>0$ consider the decomposition $f_j=f_{j,\lambda,s}+f_{j,\lambda,l}$ where $f_{j,\lambda,l}$ is the restriction of $f_j$ to a some $U_j\subset X_j$ chosen such that $S_j(f_j1_{U_j^c}) \le \lambda$ and $\mu_j(U_j) \le 2 \mu_j(S_j(f_j)>\lambda)$. From the proof of Lemma~\ref{l.outer-marcinkiewicz}, using $p_j<t_j$ we have
\begin{eqnarray*}
\int_0^\infty  \lambda^{p_j-t_j-1} \|f_{j,\lambda,s}\|_{t_j}^{t_j} d\lambda  &\lesssim& \|f_j\|_{p_j}^{p_j}=1 \ \ , \\
\int_0^\infty  \lambda^{p_j-t_j-1} \|f_{j,\lambda,l}\|_{t_j}^{t_j} d\lambda  
&\lesssim&   \|f_j\|_{p_j}^{p_j} =  1
\end{eqnarray*}
Now, given $\lambda>0$ we will decompose $K(f_1,f_2)$ by decomposing $f_1=f_{1,\alpha_1,s}+f_{1,\alpha_1,l}$ and similarly $f_2=f_{2,\alpha_2,s}+f_{2,\alpha_2,l}$ with $\alpha_1=\lambda^{p/p_1}$ and $\alpha_2=\lambda^{p/p_2}$ (clearly $\alpha_1\alpha_2=\lambda$). This leads to a decomposition of $K(f_1,f_2)$ into four terms, and we will estimate each of them  using the known weaktype estimate at one suitable vertex of $A$. We  show below the treatment for $K(f_{1,\alpha_1,s},f_{2,\alpha_2,s})$, for which we will use the vertex $(t_1,t_2)$.  Letting $t$ denote $t_1t_2/(t_1+t_2)$, it follows that
\begin{eqnarray*}
\int_0^\infty \lambda^{p-1} \mu(S(K(f_{1,\alpha_1,s},f_{2,\alpha_2,s}))>\lambda) d\lambda 
&\lesssim& \int_0^\infty \lambda^{p-t-1} \prod_{j=1,2}\|f_{j,\alpha_j,s}\|_{t_j}^{t_1t_2/(t_1+t_2)} d\lambda\\
(\text{classical H\"older})\quad &\lesssim&  \prod_{j=1,2}\Big(\int_0^\infty \lambda^{p-t_j-1} \|f_{j,\alpha_j,s}\|_{t_j}^{t_j} d\lambda\Big)^{t_{3-j}/(t_1+t_2)}\\
(\text{using $\lambda^{p-1}d\lambda = C \alpha_j^{p_j-1}d\alpha_j$}) \quad &=& C \prod_{j=1,2}\Big(\int_0^\infty \alpha_j^{p_j-t_j-1} \|f_{j,\alpha_j,s}\|_{t_j}^{t_j} d\alpha_j\Big)^{t_{3-j}/(t_1+t_2)}\\
&\lesssim& \prod_{j=1,2} \|f_j\|_{p_j}^{p_j t_{3-j}/(t_1+t_2)} \quad \lesssim \quad 1 \ \  .
\end{eqnarray*}
The other terms (in the decomposition for $K(f_1,f_2)$ could be treated similarly, thus by quasilinearity of size it follows immediately that
\begin{eqnarray*}
\|K(f_1,f_2)\|_{\mathcal L^p(X,S,\mu)} \quad = \quad (p\int_0^\infty \lambda^{p-1} \mu(S(K(f_1, f_2))>\lambda) d\lambda)^{1/p} &\lesssim& 1 \ \ .
\end{eqnarray*}

\underline{Case 2: Exactly one of $t_1,t_2$ is $\infty$.}

Without loss of generality assume that $t_1<t_2=\infty$. In this case we still carry out the same decompositions as before. The two terms that does not involve $f_{2,\alpha_2,s}$ could betreated as before. For $K(f_{1,\alpha_1,s},f_{2,\alpha_2,s})$ using the assumed weak-type estimate at $(t_1, t_2=\infty)$ we have
\begin{eqnarray*}
\mu(S(K(f_{1,\alpha_1,s},f_{2,\alpha_2,s}))>c\lambda) &\lesssim& \lambda^{-t_1} (\|f_{1,\alpha_1,s}\|_{t_1}\|f_{2,\alpha_2,s}\|_{\infty})^{t_1}\\
&\lesssim& \alpha_2^{t_1} \lambda^{-t_1}  \|f_{1,\alpha_1,s}\|_{t_1}^{t_1} \quad = \quad  \alpha_1^{-t_1}  \|f_{1,\alpha_1,s}\|_{t_1}^{t_1} \ \ ,
\end{eqnarray*}
\begin{eqnarray*}
\int_0^\infty \lambda^{p-1} \mu(S(K(f_{1,\alpha_1,s},f_{2,\alpha_2,s}))>c\lambda) d\lambda &\lesssim& \int_0^\infty \alpha_1 ^{p_1-t_1-1}\|f_{1,\alpha_1,s}\|_{t_1}^{t_1}  d\alpha_1 \quad \lesssim \quad 1 \ \ .
\end{eqnarray*}
The term $K(f_{1,\alpha_1,l},f_{2,\alpha_2,s})$ could be treated similarly.

\underline{Case 2: $t_1=t_2=\infty$.}

In this case we modify the decompositions slightly so that $S_j(f_{j,\alpha_j,s}) \le c\alpha_j$ for both $j=1,2$, where $c>0$ is sufficiently small. It follows from the assumed weak-type estimate at $(t_1,t_2)=(\infty,\infty)$ that
\begin{eqnarray*}
\|K(f_{1,\alpha_1,s},f_{2,\alpha_2,s})\|_{\infty} &\lesssim& O(c^2\alpha_1 \alpha_2) \quad = \quad O(c^2\lambda)
\end{eqnarray*}
therefore by choosing $c$ sufficiently small we obtain, for some $C>0$ large,
\begin{eqnarray*}
\nu(S'(K(f_1,f_2))>\lambda) &\le& \nu(S'(K(f_{1,\alpha_1,s},f_{2,\alpha_2,s}))>\lambda/C)  +\\
&+& \nu(S'(K(f_{1,\alpha_1,l},f_{2,\alpha_2,s}))>\lambda/C)  +\\
&+& \nu(S'(K(f_{1,\alpha_1,s},f_{2,\alpha_2,l}))>\lambda/C)  \ \ .
\end{eqnarray*}
The three   terms on the right hand side could be treated as before.
\endproof

\section{Generalized Carleson embeddings and outer $L^p$ estimates for discrete variation-norm Carlerson operators} \label{s.carl-embed}
Let  $\{\phi_{P}, P\in X\}$ be $L^1$-normalized Fourier wave packets where $\P$ is finite sparse, such that $supp\widehat \phi_{P} \subset (5/4)\omega_{P}$, and  $\{(\omega_{P,lower} , \omega_P, \omega_{P,upper}), P \in \P\}$ is   rigid. For simplicity  we assume that $\omega_{P,upper}$ is finite and $\omega_{P,lower}$ is a half line (the  other settings could be handled similarly).  For technical convenience, assume that  $1000\omega_P$ is strictly between $\omega_{P,lower}$ and $\omega_{P,upper}$ for every $P\in \P$.

Let $f$ be a Schwarz function on $\R$, and $\alpha_j:\R\to [0,\infty]$, and define
$$T_1 f(P) := \<f,\phi_P\> \ \ , \ \  T_2f(P) := \<f,  \sum_{1\le j \le K} d_j \widetilde \phi_{P,j} 1_{|I_P|\le \alpha_j(x)} \> \ \ ,$$
where $\widetilde \phi_{P,j}$ is defined by \eqref{e.modpacket}.

In this section, we consider embedding estimates for $T_1$ and $T_2$ from $L^p(\R)$ to outer measure spaces on $\P$.  Following \cite{dt2014} we will refer to these estimates   as generalized Carleson embeddings.

\subsection{Outer measure spaces}\label{s.outerP}
For every $P\in \P$,  let $\widetilde{\omega_P}$ be the convex hull of $50\omega_P$ and $50\omega_{P,upper}$.

\underline{Generating subsets:}  A nonempty $E\subset \P$ is a generating set (i.e. $E\in \E$) if  there exists a dyadic interval $I_E$ and  $\xi_E \in \R$ such that  for every $P\in E$ we have
\begin{eqnarray*}
I_P\subset I_E&,&
(\xi_E -\frac 1{2|I_E|}, \xi_E + \frac 1{2|I_E|}) \subset \widetilde {\omega_P}
\end{eqnarray*}
We say that $E$ is lacunary if furthermore $\xi_E  \in 50 \omega_{P,upper}$ for every $P \in E$, and $E$ is overlapping if $\xi_E  \in \widetilde \omega_P \setminus 50 \omega_{P,upper}$ for every $P\in E$.

\underline{Outer measure:} The outer measure $\mu$  will be generated from $\mu(E)=\inf \sum_j |I_{E_j}|$, infimum taken over all countable  coverings of $E$ by  generating sets.  

\underline{Size:}  For any $0<t<\infty$ and any generating set $E \subset \P$, let $S_t$ be the size 
$$S_{t}(f)(E) = \sup_{S\subset E} \Big(\frac 1 {|I_S|} \sum_{P\in S} |f(P)|^t |I_P| \Big)^{1/t} \ \ .$$  
If the supremum is taken over only lacunary subsets $S$, we call $S_{t,lac}$ the corresponding size, and define $S_{t,overlap}$ similarly.   When $t=\infty$ the three sizes agree
$$S_{\infty,lac}(f)(E) = S_{\infty, overlap}(f)(E) = S_\infty(f)(E) =\sup_{P \in E} |f(P)|  \ \  .$$
It is clear that $S_t$ and $S_{t,lac}$ and $S_{t,overlap}$ are decreasing functions of $t$.

\underline{Strongly disjointness:} Let $(E_m)$ be lacunary. We say they are strongly disjoint if for any  $m\ne n$:
\begin{itemize}
\item If   $P_1\in E_m$ and $P_2\in E_{n}$  such that $10\omega_{P_1} \cap 10 \omega_{P_2} \ne \emptyset$ and $|I_{P_1}|\le |I_{P_2}|$ then $I_{P_1} \cap I_{E_n} = \emptyset$. (In particular $I_{P_1} \cap I_{P_2}=\emptyset$.) 
\end{itemize}
The following estimate follows from a standard argument,  see e.g. \cite{osttw2012} or  \cite{lt1999}.
\begin{proposition}\label{l.wavelet} Let $(E_m)$ be strongly disjoint lacunary and $\P_1=\bigcup E_m$. Then
\begin{eqnarray*}\|\sum_{P\in \P_1} |I_P| f(P) \phi_P\|_{L^2(\R)}  &\lesssim&  (\sum_{P\in \P_1} |I_P| |f(P)|^2 )^{1/2} \ \  + \ \  \sup_{P\in \P_1} |f(P)| (\sum_m |I_{E_m}|)^{1/2} \ \ .
\end{eqnarray*}
\end{proposition}

\subsection{Embeddings for $T_1$}
The following is a discrete version of \cite[Theorem 5.1]{dt2014} and its $L^{2,\infty}$ endpoint is also reformulation of the well-known size lemma in \cite{lt1999}.
\begin{theorem}\label{t.std-embed} For any $2<p\le \infty$ it holds that
$\|T_1 f\|_{\mathcal L^p(\P, S_{2,lac},\mu)} \lesssim \|f\|_{L^p(\R)}$,  
and the weak-type estimate holds at $p=2$.
\end{theorem}
\proof Without loss of generality assume $\|f\|_{L^p(\R)}=1$. The endpoint $p=\infty$ is a consequent of Lemma~\ref{l.s2lac}, so by interpolation it suffices to consider the weak-type estimate  at $p=2$.

 Fix any $\lambda>0$. We will show that there exists $\P'\subset \P$ such that $\mu(\P')  \lesssim \lambda^{2}$ and $\|T_1 f\|_{\mathcal L^\infty(\P\setminus \P', S_{2,lac}, \mu)}  \le \lambda$. Without loss of generality  assume that $\|T_1f\|_{\mathcal L^\infty}>2\lambda$. We define  $\P'=\bigcup_i F_i$ where $F_i$ are selected as follows:

(i) If there is $E\subset \P$ lacunary with $\big(\frac 1{|I_E|} \sum_{P\in E} |T_1(P)|^2 |I_P|\big)^{1/2}>\lambda$,
we select one such $E_1$ with smallest possible  $\xi_E$. 

(ii) Let $F_1= \{P \in \P: (\xi_{E_1}-\frac 1{2|I_{E_1}|}, \xi_{E_1}+ \frac 1 {2|I_{E_1}|}) \in \widetilde \omega_P, \ I_P\subset I_{E_1}\}$.

We remove $F_1$ from $\P$ and repeat the above argument and select $E_m$, $F_m$, $m=2,3,\dots$. Since $\P$ is finite the  process will stop, and by geometry  and sparseness of $\P$ it is clear that $(E_m)_{i\ge 1}$ is strongly disjoint.  Let $\P_1 =\bigcup E_m$. It follows that
\begin{eqnarray*}
\lambda^2 \sum_{m} |I_{E_m}| &\lesssim&   \sum_{P\in \P_1} |I_P| |\<f,\phi_P\>|^2\\
&\lesssim& \|\sum_{P\in \P_1} |I_P| \<f,\phi_P\> \phi_P\|_2 \qquad \text{(Cauchy-Schwarz, $\|f\|_2=1$)}\\
&\lesssim& (\sum_m |I_{E_m}| )^{1/2} \|T_1 f\|_{\mathcal L^\infty(\P,S_{2,lac},\mu)} \qquad \text{(using Lemma~\ref{l.wavelet})}\\
&\lesssim& \lambda  (\sum_m |I_{E_m}| )^{1/2}  \qquad \text{(since   $\|T_1f\|_{\mathcal L^\infty}<2\lambda$).} \
\end{eqnarray*}
It follows that $\lambda^2 \mu(\P') \le \lambda^2 \sum_m |I_{E_m}| \lesssim 1$, as desired. \endproof

For any $g$ consider the following tile maximal average
\begin{eqnarray}\label{e.tileHLmax}
M_N(\P, f) := \sup_{P\in \P}  \frac 1{|I_P|} \int \widetilde \chi_{I_P}^N f  \ \ .
\end{eqnarray}
The following Lemmas are reformulation of  standard estimates, see e.g. \cite{mtt2004f}.
\begin{lemma}\label{l.s2lac-char} It holds that 
\begin{eqnarray*}
\|T_1 f\|_{\mathcal L^\infty(\P, S_{2,lac},\mu)}  &\lesssim&   \sup_{E\subset \P \  lacunary} \frac{1}{|I_E|}\|(\sum_{P\in E} |T_1 f(P)|^2 1_{I_P})^{1/2}\|_{L^{1,\infty}} \ \ .
\end{eqnarray*}
\end{lemma}
By Calderon--Zygmund theory,  for every lacunary generating set $E$ it holds that
\begin{eqnarray*}
\|(\sum_{P\in E} |T_1 f(P)|^2 1_{I_P})^{1/2}\|_{L^{1,\infty}} &\lesssim_N& \|\widetilde \chi_{I_E}^N f\|_1 
\end{eqnarray*}
Consequently, Lemma~\ref{l.s2lac} implies the following standard corollary:
\begin{lemma}\label{l.s2lac} It holds that
\begin{eqnarray*}\|T_1 f\|_{\mathcal L^\infty(\P, S_{2,lac},\mu)}  &\lesssim&  M_{N+4}(\P,f) \quad \lesssim \quad \sup_{P\in \P} \inf_{x\in I_P} M_{I_P}(\widetilde \chi_{I_P}^{N}f)(x) \ \ . 
\end{eqnarray*}
\end{lemma}

When $\P$ is a generating set,  we have the following standard estimates, part (i) is a reformulation of \cite[Proposition 3.4]{mptt2004}. For convenience, we will sketch a proof.
\begin{lemma}\label{l.std-embed-tree}
(i) Suppose that $E$ is a generating set. Then for every $1<p\le \infty$
\begin{eqnarray*}
 \|T_1 f\|_{\mathcal L^{p}(E,S_{2,lac},\mu)} &\lesssim_{p,N}& \|f\widetilde \chi_{I_E}^N\|_p \ \ ,
\end{eqnarray*}
and the weak-type endpoint $p=1$ holds.

(ii) If $E$ is lacunary then for every $1<p<\infty$
\begin{eqnarray*}
\|\sum_{P\in E} |I_P| T_1f(P) \phi_P\|_{L^p(\R)} &\lesssim_p&  \|T_1 f\|_{\mathcal L^{p}(E,S_{2,lac},\mu)} \ \ .
\end{eqnarray*}
\end{lemma}
\proof   (i) For any  $F\subset E$ by Lemma~\ref{l.s2lac} and Lemma~\ref{l.s2lac-char} we have:
\begin{eqnarray*}
\|T_1 f\|_{\mathcal L^\infty(F,S_{2,lac},\mu)} &\le&   C_N \sup_{S\subset F \ lacunary}  \inf_{x\in I_S} M  (f\widetilde \chi_{I_E}^N)(x) \ \ , \end{eqnarray*}
thus the endpoint $p=\infty$ follows. By interpolation it suffices to consider the weak-type endpoint at $p=1$.  Fix $\lambda>0$.
We select a sequence of generating subsets $S_1, S_2, \dots$ of $E$ as follows. If there exists  $J\subset I_E$ dyadic such that $ \inf_{x\in J} M  (f\widetilde \chi_{I_E}^N)(x)  >\lambda/C_N$ 
we choose a maximal $J$, and let $S_1=\{P\in E: I_P\subset J\}$ and remove $S_1$ from $E$, then repeat the above selection algorithm. Since $E$ is finite the algorithm will stop and we obtain our sequence $S_1,S_2,\dots, S_k$. Clearly on $E - S_1-\dots -S_k$ we have $S_{2,lac}(T_1 f)  \le \lambda$. Now, $I_{S_j}$'s are pairwise disjoint and $S_j$ are  generating sets, therefore
\begin{eqnarray*}
\lambda \mu(\bigcup_j  S_j) &\lesssim& \lambda |\{M(f\widetilde \chi_{I_E}^N)>C_N \lambda\}| \quad \lesssim \quad \|f\widetilde \chi_{I_E}^N\|_1  \ \ . 
\end{eqnarray*}

(ii) Let $g\in L^{p'}$, by outer Radon-Nikodym/H\"older we have
\begin{eqnarray*}
\<\sum_{P\in E} |I_P| T_1f(P) \phi_P,g\> &\lesssim& \|T_1f\|_{\mathcal L^p(E,S_{2,lac},\mu)} \|T_1 g\|_{\mathcal L^{p'}(E,S_{2,lac},\mu)}\\
&\lesssim_p&   \|T_1f\|_{\mathcal L^p(E,S_{2,lac},\mu)}  \|g\|_{L^{p'}(\R)}
\end{eqnarray*}
thanks to part (i). The desired estimate follows from duality.
\endproof

\subsection{Embeddings for $T_2$}

The following Theorem generalizes the variation-norm density lemma   in \cite{osttw2012}, which in turn generalizes the density lemma in \cite{lt1999}.

\begin{theorem}\label{t.var-embed}  Assume that $s\ne 2$, and let $q=\min(2,s')$, $s'=s/(s-1)$.
 
For any $p \in (s,\infty]$ it holds that
\begin{eqnarray*}
\|T_2f\|_{\mathcal L^{p}(\P, S_{1,overlap} + S_{q',lac}, \mu)} &\lesssim& M_N(\P, 1_{supp(f)})^{1/p'} \|f  (\sum_j |d_j|^s)^{1/s}\|_p \ \ ,
\end{eqnarray*}
and the weak-type estimate holds at  $p=s$.
\end{theorem}

By interpolation, it suffices to consider the $L^\infty$ estimate and the weaktype estimate at $p=s$. (Note that we could fix $G=supp(f)$ and invoke  interpolation theorems for the linear map $T_2$ from functions on $G$ to outer measure spaces, the factor $M_N(\P,1_G)^{1/p'}$ should be thought of as an estimate for the norm of $T_2$.)

The proof strategy will involve three Lemmas: Lemma~\ref{l.var-infty-lac} and Lemma~\ref{l.var-infty-ovl} will be used to handle the $L^\infty$ endpoint, and Lemma~\ref{l.ms-select} will be used to handle the weak type estimate at $p=2$.

For convenience of notation, for each $\P$ let  $\overline{\P}$ denote the following completion
\begin{eqnarray}\label{e.completion-def}
\overline{\P} &=& \{R:\ \  \exists P_1,P_2\in \P \text{ with } I_{P_1} \subset I_R \subset 30 I_{P_2}, \ \widetilde \omega_R \cap \widetilde \omega_{P_j} \ne \emptyset\}
\end{eqnarray} 

For convenience of notation,  for any measurable sequence $(g_j)$  we denote
\begin{eqnarray}
\label{e.ms-def}
m_{s,N}(\P, (g_j)) &:=&  \sup_{R \in \overline{\P}} \frac 1{|I_R|}\int \widetilde \chi_{I_R}^N   (\sum_{j: N_j \in \widetilde \omega_R}  |g_j|^s)^{1/s} \ \ ,
\end{eqnarray}
Clearly, $m_{\infty,N} \lesssim m_{s,N}$.  

Now, the $L^\infty$ endpoint of Theorem~\ref{t.var-embed} follows from
Using Lemma~\ref{l.var-infty-lac} and Lemma~\ref{l.var-infty-ovl}, which gives the estimate
\begin{eqnarray*}\|T_2f\|_{\mathcal L^\infty(\P, S_{1,overlap}+S_{q',lac},\mu)} &\lesssim& m_{s,N}(\P, (fd_j))\ \ .
\end{eqnarray*}
The weaktype estimate at $p=s$  of Theorem~\ref{t.var-embed} follows from the following result.

\begin{lemma}\label{l.ms-select} Let $p \in [s,\infty]$ and assume $\supp (\sum_j |g_j|^s)^{1/s}) \subset G$. 
 For any $\lambda>0$  there exists $\P_1 \subset \P$ such that
\begin{eqnarray*}
&& \lambda \mu(\P \setminus \P_1)^{1/p} \quad \lesssim_{N,p,s} \quad  M_N(\P, 1_G)^{1/p'}  \|(\sum_{j}|g_j|^{s})^{1/s}\|_{L^{p}(\R)} \\
\text{and} && m_{s,N}(\P, (g_j)) \quad \le \quad \lambda \ \ .
\end{eqnarray*}
\end{lemma}
\proof The endpoint $p=\infty$ is trivial, while the endpint $p=s$ follows from 
$$m_{s,N}(\P, (g_j)) \lesssim M_N(\P, 1_G)^{1/s'} \sup_{R\in \overline{\P}} (\frac 1 {|I_R|} \int \widetilde \chi_{I_R}^N \sum_{j: N_j\in \widetilde \omega_R} |g_j|^s)^{1/s}$$
and  a variation-norm version of the standard density lemma (see e.g. \cite{osttw2012}). The general case could be obtained by a simple interpolation argument: let $A=\{x: (\sum_j |g_j(x)|^s)^{1/s} > \lambda/C\}$ for $C>0$ large,  we use the $p=s$ and $p=\infty$ endpoints to respectively treat $g_j1_A$ and  $g_j1_{A^c}$.  We omit the details.
\endproof

\begin{lemma}\label{l.var-infty-lac}Assume that $s\ne 2$. Let $q=\min(2,s')$. Then
\begin{eqnarray*}
\|T_2f\|_{\mathcal L^\infty(\P, S_{q',lac},\mu)} &\lesssim& m_{s,N}(\P, (fd_j)) \ \ .
\end{eqnarray*}
\end{lemma}
  
\proof 
It suffices to show that for every  lacunary $E$ and $g:E\to \C$ it holds that
\begin{eqnarray}\label{e.lactree}
 \sum_{P \in E} |I_P|  T_2f(P) g(P) &\lesssim& |I_E|  \|g\|_{\mathcal L^\infty(E,S_{q,lac},\mu)}  m_{s,N}(E,(fd_j)) \ \ .
\end{eqnarray}
Indeed, taking $g(P)=\overline{T_2f(P)} |T_2f(P)|^{q'-2}$  and applying the above estimate for all lacunary subsets of $\P$, we obtain an equivalent form of the desired estimate:
\begin{eqnarray*}
\Big(S_{q',lac}(T_2f)(\P)\Big)^{q'} \lesssim \Big(S_{q',lac}(T_2f)(\P)\Big)^{q'-1} m_{s,N}(\P,(fd_j)) \ \ .
\end{eqnarray*}
Below, for brevity we write $\|g\|_\infty$ for $\|g\|_{\mathcal L^\infty(E,S_{q,lac},\mu)}$, and $m_s$ for $m_{s,N}$.

For every $P$ and $x$, it is clear that at most one $j$ will satisfy $N_{j-1}\in \omega_{P,lower}$ and $N_j\in \omega_{P,upper}$. Let $d_P(x)=d_j(x)1_{|I_P|\le \alpha_j}$ if such $j$ exists, otherwise $d_P(x)=0$. Let $\J$ be the collection of all maximal dyadic  $J$ such that $3J$ does not contain any $I_P$, $P\in E$. Clearly,  $\J$ is a partition of $\R$, so the left hand side of (\ref{e.lactree}) is bounded by
\begin{eqnarray}\label{e.Jdecomp-varcar}
&\le& \sum_{J \in \J} \|\sum_{P\in E} |I_P| g(P) \phi_{P} \, f d_P\|_{L^1(J)} \quad \le  \quad A+B\ \ ,
\end{eqnarray}
where $A$ denotes the contribution of   $(J,P)$ such that $|I_P|\le 2|J|$ and  the rest is in $B$. Since $\|\widetilde \chi_{I_P}^N fd_P\|_{1} \lesssim |I_P| m_{\infty,N}(E,(fd_j))$, we obtain
\begin{eqnarray*}
A &\lesssim& \sup_P |g(P)| m_{\infty}(E,(fd_j)) \sum_{J\in \J} \sum_{\stackrel{P: |I_P|\le 2|J|}{I_P\not\subset 3J}} |I_P| \sup_{x\in J} \widetilde \chi_{I_E}(x)^2 \widetilde \chi_{I_P}(x)^{2}   \\
&\lesssim&  \sup_P |g(P)|  m_{\infty}(E,(fd_j))  \sum_{J\in J} \sup_{x\in J} \widetilde \chi_{I_E}(x)^2\\
&\lesssim& |I_E|  \sup_P |g(P)| m_{\infty}(E,(fd_j))   \ \ .
\end{eqnarray*}
Note that the above estimate remains true when $E$ is overlapping.

We now estimate $B$ i.e. the contribution of $(J,P)$ with $|I_P|\ge 4|J|$. 

First, by geometry, such $J$ must be a subset of $3I_E$. Furthermore, by maximality there is $P_1\in E$ such that   $I_{P_1}\subset 3\pi(J)$ where $\pi(J)$ is the dyadic parent of $J$.  Let $R$ be a tile such that $(I_{P_1} \cup \pi(J)) \subset I_R$ and $|I_R|=4|J|$, and  $\xi_E\in \omega_R$. We will show that $R\in \overline{\P}$. To see this,  note that  for \eqref{e.lactree} we may assume  that $I_{P_2}=I_E$ for some $P_2\in E$. It follows that  $I_{P_1}\subset I_R \subset 10J \subset 30 I_{P_2}$; it is also clear that $\widetilde \omega_R \cap \widetilde \omega_{P_j} \ne \emptyset$ for each $j=1,2$ since they all contain $\xi_E$, thus $R\in \overline{\P}$ as claimed.

Now, if $P\in E$ such that $|I_P|\ge |I_R|=4|J|$, by sparseness it follows that $\omega_{P,upper}\subset \widetilde \omega_R$.  Now, for  $x\in J$ by H\"older's inequality we have
\begin{eqnarray*}
&& \sum_{P\in E: \ |I_P|>2|J|} |I_P| g(P) \phi_{P}(x) f(x) d_P(x)  \\
&\lesssim& \big(\sum_{j: N_j \in \widetilde \omega_R}  |fd_j|^{s}\big)^{1/s} \big(\sum_{j} |\sum_{\alpha_j(x)\ge |I_P|>2|J|} |I_P| g(P) \widetilde \phi_{P,j}|^{s'} \big)^{1/s'} \ \ .
\end{eqnarray*}
 
Since $T$ is lacunary and $\P$ is sparse, we can find a sequence of integers $O(1)+\log_2(|J|) \le m_1(x) \le n_1(x) \le   \dots \le m_K(x) \le n_K(x)$ such that
$$\{\alpha_j \ge |I_P|>2|J|: N_{j-1} \in \omega_{P,lower}, N_j \in \omega_{P,upper}\} = \{P \in E: 2^{m_j} < |I_P| \le 2^{n_j}\} \ \ .$$
(Note that when $m_j=n_j$  the set is understood to be empty, one example when this may happen is when $\alpha_j$ is too small or $2|J|$ is too large relative to the range of $|I_P|$ imposed by the $N_j, N_{j-1}$ constraints.)
Now, let $g_E(x) = \sum_{P\in E} |I_P| g(P) \phi_{P}(x)$. 
For every $x \in J$ we have
$$\big(\sum_{1\le j \le K} |\sum_{\alpha_j \ge |I_P|>2|J|} |I_P| g(P) \widetilde \phi_{P,j}|^{s'} \big)^{1/s'} = \big(\sum_{1\le j \le K} |(\Pi_{n_j}  - \Pi_{m_j}) g_E|^{s'}\big)^{1/s'}$$
where  $\Pi_n$s   are Fourier projections onto the relevant frequency scales of $E$ (essentially projecting onto $|\xi-\xi_E|\lesssim 2^{-n}$, thus larger values of $n$ means narrower bands). Using $m_1>\log_2|J| +O(1)$ and Minkowski's inequality, the last display is bounded   by $M_J( \|\Pi_{k}g_E\|_{V^{q}_k(\mathbb Z)})$.
It follows that
\begin{eqnarray*}
B &\lesssim& \sum_{J \in \J: J\subset 3I_E} M_J(\|\Pi_k g_E\|_{V^{s'}_k(\Z)})  \|\big(\sum_{j: N_j \in \widetilde \omega_R}  |fd_j|^{s}\big)^{1/s} \|_{L^1(J)}\\
&\lesssim& \|M(\|\Pi_k g_E\|_{V^{s'}_k(\Z)})\|_{L^1(3I_E)} m_{s}(E,(fd_j)) \ \ .
\end{eqnarray*}
Since $s\ne 2$ we either have $s'>2$ or $s'<2$. If $s'>2$ then by the (continuous) L\'epingle inequality \cite{lepingle1976,bourgain1989,jsw2008} we obtain
\begin{eqnarray*}
B &\lesssim& |I_E|^{1/2} \|g_E\|_2 m_{s}(E,(fd_j)) \\
&\lesssim&  |I_E| \|g\|_{L^\infty(S_{2,lac})}  m_{s}(E,(fd_j)) \ \  .
\end{eqnarray*}
If $s'<2$ then by the continuous Pisier--Xu inequality (see  \cite{dmt2012}) we obtain
\begin{eqnarray*}
B &\lesssim& |I_E|^{1/s} \|(\sum_k |\sum_{P\in E: |I_P|=2^k} |I_P| g(P)\phi_P|^{s'})^{1/s'}\|_{s'} m_s(E,(fd_j))\\
&\lesssim& |I_E| \|g\|_{\mathcal L^\infty(E,  S_{s',lac},\mu)}  m_{s}(E,(fd_j)) \ \ . \qquad \text{\endproof}
\end{eqnarray*}

\begin{lemma} \label{l.var-infty-ovl} If $1<q \le \infty$ then
\begin{eqnarray*}
\|T_2f\|_{\mathcal L^\infty(\P, S_{1,overlap}, \mu)} &\lesssim& m_{\infty,N}(\P, (fd_j)) \ \ .
\end{eqnarray*}
\end{lemma}

\proof The proof is entirely similar to the proof of Lemma~\ref{l.var-infty-lac}, it suffices to show for any overlapping $E$ the following estimates
\begin{eqnarray}\label{e.ovltree}
 \sum_{P \in E} |I_P|  T_2f(P) g(P) &\lesssim& |I_E|  \sup_P |g(P)|  m_{\infty,N}(E,(fd_j)) \ \ .
\end{eqnarray}
We define $\J$ as before and invoke \eqref{e.Jdecomp-varcar} again, and $A$ could be estimated as before. To estimate $B$ we use the following observation: since $E$ is overlapping, for every $x$ there is at most one $j$ and at most one scale of $E$  that  contributes to  $ \sum_{|I_P|>2|J|} |I_P| g(P)\phi_{P}(x) d_P(x)$. Let $R\in \overline{E}$ be as before. It follows that
\begin{eqnarray*}
B &\lesssim& \sum_{J\subset 3I_E} \|\sup_{j: N_j \in \widetilde \omega_R} \sup_{\alpha>2|J|} \sum_{P\in E: \ |I_P|=\alpha} |g(P)|\widetilde \chi_{I_P}^{N+4} d_j f\|_{L^1(J)}\\
 &\lesssim& \sum_{J\subset 3I_E} \|\sup_{j: N_j\in \widetilde \omega_R} |d_{j} f|\|_{L^1(J)} \sup_{P} |g(P)| \\
&\lesssim_N& |I_E| \sup_P |g(P)|  m_{\infty,N}(E,(fd_j)) \ \ . \ \  \text{\endproof}
\end{eqnarray*}

\subsection{Outer $L^p$ estimates for discrete variation-norm Carleson operators}

Let $r\ne 2$ and $q=\min(2,r)$, and $\widetilde \phi_{P,j}$'s are defined relative to  $(N_j(x))$ and $K(x)$.  We consider the variation-norm operator $$V_{s'}(g) = \big(\sum_j  \sup_\alpha|\sum_{P\in \P: |I_P|\le \alpha} |I_P| g(P) \widetilde \phi_{P,j}|^{r}\big)^{1/r} $$

\begin{theorem}\label{t.varCarl} Let $F$ be a subset of $\R$.

(i) For any $1< p <r$  it holds that
\begin{eqnarray*}
\|V_{r}(g)\|_{L^p(F)} &\lesssim_{p,N}&  M_N(\P, 1_F)^{1/p} \|g\|_{\mathcal L^p(\P, S_{q,lac},\mu)}  \ \ .
\end{eqnarray*}
(ii) If $r>2$   then for all $p>r'$ it holds that
\begin{eqnarray*}
 \|  V_{r}(T_1 f)\|_p  &\lesssim&   \|f\|_{L^p(\R)}  \ \ .
\end{eqnarray*}
\end{theorem}
\proof (i) Let $s=r'$. We may find  $(\alpha_j)$  and    $(g_j)$  measurable functions with $\alpha_j\ge 0$ and $\sum_{1\le j\le K} |g_j|^{s} = O(1)$  such that
$$V_{r}(g) \lesssim  \sum_{P} |I_P| g(P) \sum_j \widetilde \phi_{P,j} 1_{|I_P|\le \alpha_j(x)} g_j$$
Let $h\in L^{p'}(\R)$ and $T_2 h(P)=\<\sum_j \widetilde \phi_{P,j} g_j 1_{|I_P|\le \alpha_j(x)}, h1_F \>$. It suffices to show that
\begin{eqnarray*}
 \sum_{P} |I_P| g(P) T_2h(P) &\lesssim_{p,N}& M_N(\P, 1_F)^{1/p} \|g\|_{\mathcal L^p(\P, S_{q,lac},\mu)} \|h\|_{L^{p'}(\R)}
\end{eqnarray*}
Via applications of the classical H\"older inequality it follows that 
$$S_1(gT_2h)\lesssim S_{2,lac}(g) (S_{2,lac}+S_{1,overlap})(T_2 h)\ \ . $$
Thus, using outer Radon-Nykodym and outer H\"older inequalities, we obtain
\begin{eqnarray}\label{e.bilinear-VC}
 \sum_{P} |I_P|g(P)T_2h(P)  &\lesssim& \|g\|_{\mathcal L^p(\P, S_{q,lac},\mu)} \|T_2h\|_{\mathcal L^{p'}(\P, S_{q',lac}+S_{1,overlap},\mu)}\ \ .
\end{eqnarray}
Now, using  Theorem~\ref{t.var-embed} and noticing $s'=r$ it follows that
\begin{eqnarray*}
\|T_2h\|_{\mathcal L^{p'}(\P, S_{q',lac}+S_{1,overlap},\mu)} &\lesssim& M_N(\P, 1_F)^{1/p} \|h(\sum_j |u_j|^{s})^{1/s}\|_{L^{p'}(\R)}\\
&\lesssim& M_N(\P, 1_F)^{1/p} \|h\|_{L^{p'}(\R)}  \ \ ,
\end{eqnarray*} 
provided that $p'> s$.  This completes the proof of part (i). 

(ii)  Since $r>2$, we obtain $q=2$.  We say that a subset  is   major if it has at least half of the total measure. By restricted weak type interpolation \cite{mtt2002} (see also Section~\ref{s.extendrange-BC}), it suffices to show that we could find $\beta$ arbitrarily close to $0$ and also arbitrarily close to $1/s$ such that $\<V_r(T_1 f),h\>$ is of restricted weak type $(\beta,1-\beta)$. That is,  there exists $j_0\in \{1,2\}$ depending only on $\beta$ such that given any $F_1, F_2\subset\R$ with finite positive Lebesgue measures we could find $S_1 \subset F_1$ and $S_2\subset F_2$ both major subsets  and furthemore $S_{j_0}=F_{j_0}$ and
$$\<V_{r}(T_1 f),h\> \lesssim |F_1|^\beta |F_2|^{1-\beta} \quad \text{ if $|f|\le 1_{S_1}$ and $|h|\le 1_{S_2}$.}$$

To get $\beta$ near $1/s$, we let $S_1=F_1$ and $S_2=F_2\setminus E$ where $E:=\{M 1_{F_1}>C|F_1|/|F_2|\}$ and $C$ is large enough, and $S_2=F_2$. Without loss of generality  assume that $1+ \frac{\dist(I_P,E^c)}{|I_P|} \sim 2^k$ provided that we have enough decay in the estimate.  By convexity, for $2<p<s'$  it follows from part (i) that 
\begin{eqnarray*}
\<V_{s'}(g),h\> &\lesssim& \|g\|_{\mathcal L^p(\P,S_{2,lac},\mu)}  M_N(\P, 1_{F_2})^{1/p}   \|h\|_{p'}\\
&\lesssim& \|g\|_{\mathcal L^\infty(\P, S_{2,lac},\mu)}^{1-\frac 2p} \|g\|_{\mathcal L^{2,\infty}(\P, S_{2,lac},\mu)}^{\frac 2 p} M_N(\P, 1_{F_2})^{1/p}   |F_2|^{1/p'}\\
&\lesssim& 2^{-Nk} (\sup_{P\in \P} \sup_{x\in I_P} M1_{F_1}(x))^{1-\frac 2p} |F_1|^{1/p} |F_2|^{1/p'}\\
&\lesssim& 2^{-k} (|F_1|/|F_2|)^{1- 2/p} |F_1|^{1/p} |F_2|^{1/p'} \quad = \quad 2^{-k} |F_1|^{\beta} |F_2|^{1-\beta} \ \ ,
\end{eqnarray*}
where $\beta:=1/p'$, which is arbitrarily close to $1/s$ if  $p$ is sufficiently close to $s'$.

To get $\beta$ near $0$, we let $S_1=F_1\setminus E$ and $S_2=F_2$ where $E:=\{M1_{F_2}> C|F_2|/|F_1|\}$ with $C$ sufficiently large. Similarly we obtain
\begin{eqnarray*}
\<V_{s'}(g),h\> &\lesssim& \|g\|_{\mathcal L^\infty(\P, S_{2,lac},\mu)}^{1-\frac 2p} \|g\|_{\mathcal L^{2,\infty}(\P, S_{2,lac},\mu)}^{\frac 2 p} M_N(\P, 1_{F_2})^{1/p}   |F_2|^{1/p'}\\
&\lesssim& 2^{-Nk}   |F_1|^{1/p} (2^k |F_2|/|F_1|)^{1/p} |F_2|^{1/p'}\\
&\lesssim&  2^{-k} |F_2|\ \ . \qquad \text{\endproof}
\end{eqnarray*}

\section{Estimates for $BC$ model operators}\label{s.BC-estimate}

\begin{theorem}\label{t.BC-estimate} Let $T$ be a  BC model operator. Then $\|T\|_{L^{q_1}\times L^{q_2} \to L^q} < \infty$
for every  $1/q  = 1/ {q_1} + 1/{q_2}$ such that  $\frac{2r}{3r-4} < q_1,q_2\le \infty$ and $q> \frac {r'}2$.
\end{theorem}

For simplicity, we assume that $\omega_{P,upper}$ is finite and $\omega_{P,lower}$ is a halfline for each $P\in \P$; other situations are either symmetric or could be reduced to this setting.

The outer measure spaces on $\P$ are defined as in Section~\ref{s.outerP}. Below we discuss the outer measure spaces on  $\Q$, which are similar to the settings on $\P$, thus we only discuss  the needed changes. By further decomposition if necessary, we may assume that $\P$ and $\Q$ are very sparse.

For any two tiles $R, R'$ we  say that $R < R'$ if  $I_{R} \subsetneq I_{R'}$ and $5\omega_{R'} \subset 5\omega_{R}$.

We say that $R\le R'$ if   $R<R'$ or $R=R'$. Clearly, $\le$  and $<$ are transitive.

\underline{Generating subsets of $\Q$:}  A nonempty $E\subset \Q$ is a generating set if for some tritile $Q_E$ (with the same rigidity) the following holds: for every $Q\in E$ there is $j=j(Q)\in \{1,2,3\}$ such that $Q_j \le Q_{E,j}$,
where $Q_j$ and $Q_{E,j}$ are the $j$-tiles of $Q$ and $Q_E$.  We denote $I_E \equiv I_{Q_E}$ and $\xi_E = c(\omega_{Q_E})$.

For any fixed $j \in \{1,2,3\}$, we say that $E$ is $j$-overlapping if $j(Q)=j$  for every $Q\in E$. We say that $E$ is $j$-lacunary if it is $k$ overlapping for some $k\in \{1,2,3\} \setminus \{j\}$.

\underline{Outer measure on $\Q$:} Let $\sigma$ be generated  from $\sigma(E)=\inf \sum_j |I_{E_j}|$, infimum taken over all countable  coverings of $E$ by  generating sets.  

\underline{Size:}  For every $0<t \le\infty$, we define $S_{t}$ just as in the setting for $\P$. If the defining supremum is taken over all $j$-lacunary $S$, we obtain $S_t^{[j]}$. Similarly,  for $S_t^{[j_1,j_2]}$  the supremum is taken over all $S$ that is both $j_1$ and $j_2$ lacunary, or equivalently $k$ overlapping where $k\ne j_1,j_2$.

\underline{Strongly disjointness:} For any $j\in \{1,2,3\}$,  a collection $(E_m)$ of $j$-lacunary sets of tritiles is strongly disjoint if   the following holds for any $E_m$, $E_n$, $m\ne n$:
\begin{itemize}
\item If $Q\in E_m$ and $R\in E_{n}$  such that $3\omega_{Q_j} \cap  3 \omega_{R_j} \ne \emptyset$ and $|I_{Q}| \le |I_R|$ then $I_{Q} \cap I_{E_n} = \emptyset$. (In particular $I_Q \cap I_{R}=\emptyset$.) 
\end{itemize}
An analogue of Lemma~\ref{l.wavelet} also holds in the current setting.

In the rest of this section, for every $E\subset \Q$ let $\P_E$ denote the set of all $P\in \P$ such that there exists at least one $Q\in E$ with
$|I_P| \le  |I_Q| \quad \text{  and } \quad \frac 54\omega_P  \cap \frac 56 \omega_{3,Q} \ne \emptyset$.
The following observation from \cite{mtt2004f} will be useful in the proof.

\begin{observation} \label{o.decoupleBC} Let $E$ be $3$-lacunary, then the following holds for every $P\in \P_E$ and  $Q\in E$: if $\frac 54\omega_P  \cap \frac 56 \omega_{3,Q} \ne \emptyset$ then $|I_P| \le |I_Q|$.
\end{observation}
\proof Assume the contrary, that is for some $P\in \P_E$ and $Q, R \in E$ it holds that  $\frac 54\omega_P  \cap \frac 56 \omega_{3,Q} \ne \emptyset$, $\frac 54 \omega_P \cap \frac 56 \omega_{R_3} \ne \emptyset$, and $|I_R| \ge |I_P|\ge 2|I_Q|$.  It follows that $|I_R| \ge 2|I_Q|$, therefore using sparseness of $\Q$ it is clear that  $5\omega_{R_3} \subset 5\omega_{Q_3}$, which contradicts the fact that $E$ is $3$-lacunary and  $\Q$ is sparse. \endproof

Now, fix a Schwarz function $f_3$ on $\R$. Let $a_j(Q):=\<f_j,\phi_{j,Q}\>$ for $j=1,2$, and $d(P)= \<f_3,\widetilde{\phi}_{P}\>$ where  $\widetilde \phi_P:= \sum_{j} \widetilde \phi_{P,j}d_j$. Let $K$ and its adjoint be defined by
\begin{eqnarray*}
(K f)(Q) &:=& \sum_{P\in \P: \, |I_P| \le |I_Q|} |I_P| f(P) \<\phi_P, \phi_{3,Q}\> \ \ , \ \ Q\in \Q \ \ , \\
K^\ast f(P) &:=& \sum_{Q \in \Q: \ \ |I_P| \le |I_Q|}  |I_Q| f(Q) \<\phi_{3,Q},\phi_P\> \ \ , \ \ P \in \P \ \ .
\end{eqnarray*}
Now $\<T_{BC}(f_1,f_2),f_3\> = \sum_{Q \in \Q} |I_Q| a_1(Q)a_2(Q) \overline{(K d)(Q)}$, which will be estimated using outer measure techniques.

\subsection{Outer $L^p$ estimates for $K$ and $K^\ast$}

Reall that $\Q$ has rank $1$. The following two Lemmas are the main estimates of this section.

\begin{lemma}\label{l.K-BC} Let $1< q\le 2$. For every $p\in (q',\infty)$ we have
\begin{eqnarray*}
\|K  g\|_{\mathcal L^{p}(\Q, S_2^{[3]},\sigma)} &\lesssim& \|g\|_{\mathcal L^{p}(\P, S_{q',lac} + S_{1,overlap},\mu)} \ \ .
\end{eqnarray*}
\end{lemma}

\begin{lemma}\label{l.Kadj-BC} Let $q\in (1,2]$. Then for any $1<p<q$ we have
\begin{eqnarray*}
\|K^\ast f\|_{\mathcal L^p(\P, S_{q,lac}, \mu)} &\lesssim& \|f\|_{\mathcal L^p(\Q, S_2^{[3]}+ S_1^{[1,2]}, \sigma)} \ \ .
\end{eqnarray*}
\end{lemma}

We will deduce Lemma~\ref{l.K-BC}  from Lemma~\ref{l.Kadj-BC} using a simple duality argument. 
By interpolation it suffices to consider weak-type estimates. Fix $\lambda>0$ and $q'<p<\infty$. Without loss of generality assume $\|K g\|_{\mathcal L^\infty(\Q, S_2^{[3]}, \sigma)} \le 2\lambda$. Similar to the proof of Theorem~\ref{t.std-embed},  we may select\footnote{For more details  see the proof of the $L^2$ case of Lemma~\ref{l.K-disjoint-BC}.}   a strongly disjoint collection of $3$-lacunary sets $(E_m)$ and $\Q'$ containing all $E_m$ such that $\|K g\|_{\mathcal L^\infty(\Q \setminus \Q', S_2^{[3]}, \sigma)} \le \lambda$
and 
$$\lambda^2 \sigma(\Q')  \lesssim  M:= \sum_m  \sum_{Q\in E_m} |I_Q| |K g(Q)|^2 \ \  \ \ . $$
Now, let $f(Q)=K g(Q)$ for $Q\in \bigcup_m E_m$ and zero elsewhere, we obtain via applications of the outer Radon--Nikodym/H\"older inequalities
\begin{eqnarray*}
M &=& \sum_{Q\in \Q} |I_Q|  f(Q)\overline{(K  g)(Q)}  = \sum_{P\in \P} |I_P|  (K^\ast f)(P) \overline{g(P)}\\
&\lesssim& \|K^\ast f\|_{\mathcal L^{p'}(\P, S_{q,lac},\mu)} \|g\|_{\mathcal L^{p}(\P, S_{q',lac}+S_{1,overlap},\mu)} \ \ .
\end{eqnarray*}
Therefore, by Lemma~\ref{l.Kadj-BC} we obtain
\begin{eqnarray*}
M &\lesssim&  \|f\|_{\mathcal L^{p'}(\Q, S_2^{[3]}+S_1^{[1,2]},\sigma)} \|g\|_{\mathcal L^{p}(\P, S_{q',lac}+S_{1,overlap},\mu)}  \ \ .
\end{eqnarray*}
Note that $f$ is supported on $\Q'':=\bigcup E_m$. It is clear that   any $3$-overlapping subset of $\Q''$ is essentially an union of spatially disjoint tritiles. Therefore  $S_1^{[1,2]}(f) \lesssim S_\infty(f) \lesssim S_2^{[3]}(f) \lesssim \lambda$. We obtain
\begin{eqnarray*}
M &\lesssim& \lambda \sigma(\Q')^{1/p'} \|g\|_{\mathcal L^{p}(\P, S_{q',lac}+S_{1,overlap},\mu)}  
\end{eqnarray*}
 Collecting estimates we obtain the desired weak-type estimate
$$\lambda \sigma(\Q')^{1/p} \lesssim \|g\|_{\mathcal L^{p}(\P, S_{q',lac}+S_{1,overlap},\mu)}  \ \ .  $$

In the rest of the section, we prove Lemma~\ref{l.Kadj-BC}.

\proof[Proof of Lemma~\ref{l.Kadj-BC}] By interpolation, it suffices to prove weak-type estimates for any fixed $p\in [1, q)$.   Without loss of generality, assume  $\|f\|_{\mathcal L^{p}(\Q, S_2^{[3]} + S_1^{[1,2]},\sigma)}=1$.

Fix any $\lambda>0$. We will show that there exists $\P'\subset \P$ such that $\mu(\P')  \lesssim \lambda^{p}$ and $\|K^\ast f\|_{\mathcal L^\infty(\P\setminus \P', S_{q,lac}, \mu)}  \le \lambda$. Without loss of generality we may assume that
\begin{equation}\label{e.apriori-size}
\|K^\ast f\|_{\mathcal L^\infty(\P, S_{q,lac}, \mu)}  < 2 \lambda \ \ .
\end{equation}
The construction of  $\P'$ is similar to the proof of Theorem~\ref{t.std-embed}, and we also obtain a strongly disjoint collection of lacunary sets   $(E_m)_{m\ge 1}$ contained inside $\P'$ with the following property: 
\begin{eqnarray*} \lambda^q \mu(\P') \le  \lambda^q \sum_{m} |I_{E_m}| \le   \sum_{P\in \bigcup E_m} |I_P| |(K^\ast f)(P)|^{q} \ \ .
\end{eqnarray*}

Let  $h(P) :=   K^\ast f(P)  |K^\ast f(P)|^{q-2}$ for $P\in \P'$ and zero elsewhere, and let $M$ be the last right hand side, which could be rewritten as
$$M  = \sum_{P\in \P} |I_P| (K^\ast f)(P) \overline{h(P)}  = \sum_{Q\in \Q} |I_Q| f(Q)  \overline{K(h)(Q)} \ \ .$$  
 By the classical H\"older inequality we have
\begin{eqnarray*}
S_1(fg) \quad \lesssim \quad S_2^{[3]}(f) S_2^{[3]}(g) + S_1^{[1,2]}(f) S_{\infty}(g) &\lesssim& \big(S_2^{[3]}  + S_1^{[1,2]}\big)(f)  S_2^{[3]}(g)  \ \ .
\end{eqnarray*}
Using outer Radon--Nikodym, outer H\"older, the normalization $\|f\|_{\mathcal L^p}=1$, we have
\begin{eqnarray*}
M &\lesssim& \| f\ K h \|_{L^1(\Q, S_1, \sigma)} \quad \lesssim\quad   \|K  h\|_{L^{p'}(\Q, S_2^{[3]},\sigma)} \ \ .
\end{eqnarray*}
Using Lemma~\ref{l.K-disjoint-BC} and \eqref{e.apriori-size} and the definition for $h$, it follows that
\begin{eqnarray*}
\|K h\|_{L^{p'}(\Q,S_2^{[3]},\sigma)} 
&\lesssim&  \lambda^{q/q'} (\sum_{k} |I_{E_k}|)^{1/p'} \quad \lesssim \quad  \lambda^{q/q'-q/p'}  M^{1/p'}  \  \ .
\end{eqnarray*}
Collecting estimates we obtain $M^{1/p} \lesssim \lambda^{q/q'-q/p'}$, therefore
\begin{eqnarray*}\mu(\P')^{1/p} 
&\lesssim& \lambda^{-q/p} M^{1/p}   \quad \lesssim \quad \lambda^{-1}  \ \ .   \ \  \text{\endproof}
\end{eqnarray*}

\begin{lemma}\label{l.K-disjoint-BC} Let $q\in (1,2]$ and $(E_k)$ be strongly disjoint lacunary in $\P$.
Then for $q_1\in (q',\infty]$ and every $h:\P\to \C$ that vanishes outside $\P':=\bigcup E_k$ it holds that
\begin{eqnarray*}
\|K h\|_{\mathcal L^{q_1}(\Q,S_2^{[3]},\sigma)} 
&\lesssim& \|h\|_{\mathcal L^\infty(\P_1, S_{q',lac},\mu_1)} (\sum_k |I_{E_k}|)^{1/q_1}
\end{eqnarray*}
and weak type estimates hold at $q_1=q'$.

\end{lemma}
Remark: While  Lemma~\ref{l.K-disjoint-BC} is weaker than Lemma~\ref{l.K-BC}, it will be directly proved as part of the proof of Lemma~\ref{l.Kadj-BC}. 

\proof  

It suffices to prove  the following stronger estimate, which holds for  $q_1>2$:
\begin{eqnarray}\label{e.K-stronger}
\qquad \|K h\|_{\mathcal L^{q_1}(\Q,S_2^{[3]},\sigma)}  &\lesssim&  (\sum_{P\in \P} |I_P| |h(P)|^{q_1})^{1/q_1} \ \ + \ \  \sup_{P\in \P} |h(P)| (\sum_{k} |I_{E_k}|)^{1/q_1} \ \ .
\end{eqnarray}
Clearly, if $q_1\ge q'$ then the desired conclusion follows from  \eqref{e.K-stronger}.

By interpolation, it suffices to prove weak-type estimates at $q_1=2$ and $q_1=\infty$. Since the right hand side of \eqref{e.K-stronger} is not technically an $L^p$ norm, we will detail the interpolation argument. Assume that the weak type estimates hold at $q_1=2$ and $q_1=\infty$. Let $2<q_1<\infty$,  for each $\lambda>0$ we decompose $h(P) = h_1(P)+h_2(P)$ where $h_1(P)=h(P)1_{|h(P)|\le \lambda/C}$ and $C$ is  sufficiently large but about the size of the  norm for the assumed case  $q_1=\infty$ of \eqref{e.K-stronger}. Using sublinearity of $K$, it follows that 
$$\sigma(S^{[3]}_2(K h)>\lambda) \le \sigma(S^{[3]}_2(Th_2) \gtrsim \lambda) \ \ . $$
By the assumed $L^2$ case of \eqref{e.K-stronger}, the last display is bounded above by
\begin{eqnarray*}
&\lesssim& \lambda^{-2}\Big(\sum_{P\in \P} |I_P| |h(P)|^2 1_{\lambda < C|h(P)|} \quad  + \quad  \sup_{P} (|h(P)|^2 1_{\lambda <C|h(P)|})  (\sum_{k} |I_{E_k}|)\Big)\\
&=& \lambda^{-2}\Big(\sum_{P\in \P} |I_P| |h(P)|^2 1_{\lambda <C |h(P)|} \quad + \quad \sup_{P} |h(P)|^2 1_{\lambda <C\sup_P |h(P)|} (\sum_{k} |I_{E_k}|)\Big)
\end{eqnarray*}
here in the second term we are able to move the sup inside because for any $a\ge 0$ $g(x):=x1_{a<x}$ is increasing for $x\in [0,\infty)$. Now, multiplying both side with $\lambda^{q_1-1}$ and integrate over $\lambda \in (0,\infty)$ we will obtain the desired estimate.
This completes the interpolation argument.

\underline{\bf Case 1: $q_1=\infty$.} 
By a standard characterization for $S_2$ (see e.g. \cite[Lemma 6.4]{mtt2004f}) it suffices to show that if $E \subset \Q$ is $3$-lacunary then
\begin{equation}\label{e.K-infinity}
\frac{1}{|I_E|} \|(\sum_{Q\in E} |Kh(Q)|^2 1_{I_Q})^{1/2}\|_{L^{1,\infty}(\R)} \lesssim \sup_{P\in \P} |h(P)|  \ \ .
\end{equation}
Since $\Q$ is very sparse, it is clear that  for any  interval $I$ there is at most one  $Q \in E$ such that $I_Q=I$. This remark will be used implicitly below.

Let $\P_E$ be defined as  in Observation~\ref{o.decoupleBC} and let $\P'_E=\P'\cap \P_E$. It follows that $K  h(Q) =  \<\phi_{3,Q}, h_E\>$ where $h_E(x) := \sum_{P\in \P'_E} |I_P| \overline{h(P)}\phi_{P}(x)$, for every $Q\in E$. Also, since $(E_k)$  is strongly disjoint lacunary,  it is clear that the intervals  $\{I_P, P\in \P'_E\}$ are pairwise disjoint.  

Now, by a standard argument (see e.g.  \cite[Lemma 6.8]{mtt2004f}),  it follows that
$$\frac{1}{|I_E|} \|(\sum_{Q\in E} |Kh(Q)|^2 1_{I_Q})^{1/2}\|_{L^{1,\infty}(\R)}  \lesssim_N  \frac{1}{|I_E|}\int \widetilde \chi_{I_E}(x)^N |h_E(x)|  dx \ \ .$$
Decompose $\P'_E =  \P'_{E,1} \cup \P'_{E,2}$ where $\P'_{E,1}=\{P\in \P'_E: I_P \cap 3I \ne \emptyset\}$.  It suffices to show that for every interval $I$ of the same length as $I_E$ it holds for $j=1,2$ that
\begin{eqnarray*}
A_j &:=& \frac{1}{|I|} \int_I |\sum_{P \in \P'_{E,j}} |I_P| h(P) \phi_{P}(x)| dx \quad \lesssim \quad \sup_P |h(P)| \ \  .
\end{eqnarray*} 
For $A_1$, notice first that for every  $P \in \P'_{E,1}$ we have $I_P \subset 5I$. 
Since $\phi_P$ is $L^1$-normalized and $\{I_P, P\in \P'_{E,1}\}$ are disjoint, we obtain 
$$A_1 \lesssim \frac 1{|I|}\sum_{P \in \P'_{E,1}} |I_P||h(P)|  \lesssim \sup_{P\in \P}  |h(P)| \ \ .$$
To estimate $A_2$,    we decompose the summation over $\P'_{E,2}$ according to the length of $I_P$. 
 Recall that  $\{I_P: P\in \P'_{E,2}\}$  are disjoint and disjoint from $3I$. The  desired estimate for $A_2$ follows from the following pointwise estimate on $I$:
$$\sum_{P\in \P'_{E,2}} |I_P||\phi_{P_1}(x)| \lesssim \sum_n \sum_{P: |I_P|=2^n} (\frac{|I_P|}{|I|})^4  \lesssim 1 \ \ .$$

\underline{\bf Case 2: $q_1=2$.}  
Fix $\lambda>0$. We need to show existence of $\Q'\subset \Q$ such that 
\begin{eqnarray}\label{e.K-L2}
\lambda^2 \sigma(\Q') &\lesssim& N:= \sum_{P\in \P'} |I_P| |h(P)|^2 + \sup_{P\in \P'} |h(P)|^2 \sum  |I_{E_k}|   
\end{eqnarray}
and $\|K h\|_{\mathcal L^\infty(\Q\setminus \Q')} \le \lambda$. Without loss of generality, assume that $\|K  h\|_{\mathcal L^\infty(\Q)} \le 2\lambda$. 

Now, we will construct $\Q'=\Q'_1\cup \Q'_2$ where $\Q'_j$ is union of a collection $\S_j$ of strongly disjoint $j$-lacunary sets, to be selected below.  We first collect  $\S_1$ using the following algorithm:

 (i) If there exists $G \subset \Q$ is $1$-lacunary and $(\frac 1 {|I_G|}\sum_{Q\in G} |I_Q| |K  h(Q)|^{2})^{1/2} > \lambda$,
we select one such $G_1$ with smallest possible $c(\omega_{Q_G})$.

(ii) Remove  all $Q\in \Q$ such that for some $j=1,2,3$ we have $Q_j \le Q_{G_1,j}$.

We repeat the above argument and continue selecting $G_2, G_3, \dots$. Since $\Q$ is finite the selection argument will stop. By geometry, $(G_j)$ is strongly disjoint. Let $\Q'_1$ be the set of all tritiles removed from $\Q$.

We similarly collect $\S_2$ a collection of strongly disjoint $2$-lacunary  sets, the difference is we maximize $c(\omega_{Q_E})$ in step (i). Without loss of generality, assume that $\mu(\Q')\lesssim \mu(\Q'_1)$, which we will estimate  below.

Let $k(Q)=K  h(Q)$ for $Q\in \Q'_1$ and zero elsewhere, clearly
\begin{eqnarray*}
\lambda^2 \mu(\Q') &\lesssim& \lambda^2 \mu(\Q_1')  \lesssim \lambda^2 \sum_{j} |I_{G_j}| \lesssim M \ \ ,
\end{eqnarray*}
\begin{eqnarray*}
M  &:=& \sum_{Q\in \Q'_1} |I_Q| |k(Q)|^2 = \sum_{P\in \P'} \sum_{Q\in \Q'_1}1_{|I_P|\le |I_Q|} |I_Q| |I_P| \overline{k(Q)} h(P) \<\phi_P,\phi_{3,Q}\>  \ \ .
\end{eqnarray*}
Now, to show \eqref{e.K-L2} for $\Q'_1$ it suffices to show that
\begin{eqnarray}\label{e.desiredM}
M &\lesssim&  N^{1/2} M^{1/2}  \ \ .
\end{eqnarray}

We first estimate the analogous double sum $M_{free}$ where  we don't include the coupling condition $|I_P|\le |I_Q|$. By Cauchy--Schwarz and Lemma~\ref{l.wavelet} and the assumptions $\|K h\|_{\mathcal L^\infty(\Q, S_2^{[3]}, \mu)} \lesssim  \lambda$, it follows that
\begin{eqnarray*}
M_{free} &\lesssim& \|\sum_{P\in \P'} |I_P|h(P)\phi_P\|_2  \|\sum_{Q\in \Q'_1}|I_Q| k(Q) \phi_{3,Q}\|_2  \quad \lesssim \quad N^{1/2} M^{1/2} \ \ .
\end{eqnarray*}
 
 We now consider the diagonal sum $M_{diag}$ where $|I_P|=|I_Q|$. We may further assume that $|I_P|+\dist(I_Q,I_P) \sim 2^n |I_Q|$ provided that there is an extra decaying factor  in the estimate. Note that in order for $\<\phi_P,\phi_{3,Q}\>$ to be nonzero the frequency support of $\phi_{P}$ and $\phi_{3,Q}$ must overlap. Thus essentially $Q$ is determined from $P$ and vice versa. Using $|\<\phi_P,\phi_{3,Q}\>| \lesssim 2^{-n} |I_P|^{-1/2} |I_Q|^{-1/2}$ and Cauchy--Schwarz,
\begin{eqnarray*}
M_{diag} &\lesssim& 2^{-n} (\sum_{P\in \P'} |I_P| |h(P)|^{2})^{1/2} (\sum_{Q\in \Q'_1} |I_Q| |k(Q)|^2)^{1/2} \quad \lesssim \quad 2^{-n}N^{1/2}  M^{1/2} \ \ .
\end{eqnarray*}
Thus, to prove  \eqref{e.desiredM}  the roles of $\P'$ and $\Q'_1$ are fairly symmetric.  We will assume below that  $\sum_{m} |I_{G_m}| \le \sum_{k} |I_{E_k}|$, the proof for the other case is entirely similar. Using Cauchy Schwarz, it suffices to show that for every $G\in \{G_1,G_2,\dots\}$   we have
\begin{eqnarray}
\nonumber M_G &:=&\sum_{Q\in G}\sum_{P\in \P'}1_{|I_P|\le |I_Q|} |I_P| |I_Q| h(P) \overline{k(Q)} \<\phi_P,\phi_{3,Q}\> \\
\label{e.MG} &\lesssim&   \sup_{P} |h(P)| \Big( (\sum_{Q\in G} |k(Q)|^2 |I_Q|)^{1/2} |I_G|^{1/2} + \sup_{Q} |k(Q)| |I_G|\Big)  \ \ .
\end{eqnarray}
Now, similar to Observation~\ref{o.decoupleBC},  we may write
$$M_G=\sum_{Q\in G}|I_Q| \overline{k(Q)} \<h_G, \phi_{3,Q}\> \quad, \quad h_G :=  \sum_{P\in \P'_G} h(P)|I_P| \phi_P \ \ ,$$ 
and  $\P'_G$ contains all $P\in \P'$ such that for some $Q\in G$ we have $\frac 54\omega_P \cap \frac 56\omega_{3,Q} \ne \emptyset$ and $|I_P|\le |I_Q|$.  Using strong disjointness, it is clear that the intervals of the elements of $\P'_G$ are essentially pairwise disjoint.

We now estimate the contribution of those $P\in  \P'_G$  such that $I_P\cap4 I_G \ne  \emptyset$. Clearly we will have $I_P \subset 6I_G$. Let $h_{E,0}$ be the corresponding subsum of $h_G$. By Cauchy Schwarz and standard Calderon--Zygmund theory, we have
\begin{eqnarray*}
\sum_{Q\in G}|I_Q| \overline{k(Q)} \<h_{G,0}, \phi_{3,Q}\> 
&\lesssim& (\sum_{Q\in G} |I_Q| |k(Q)|^2)^{1/2} (\sum_{P\in \P'_G: I_P \subset 6I_G} |I_P| |h(P)|^2)^{1/2}\\
&\lesssim& (\sum_{Q\in G} |I_Q| |k(Q)|^2)^{1/2} |I_G|^{1/2} \sup_{P} |h(P)|  \ \ ,
\end{eqnarray*}
in the last estimate we used disjointness of the intervals of elements of $\P'_G$.

Consider the contribution of other $P$'s.  Let $\P'_{G,k}$ contains all $P\in \P'_G$ such that $I_P\cap 2^{k+2}I_{G}\ne \emptyset$ but $I_P \cap 2^{k+1} I_{G} =\emptyset$. Let $\Lambda:=\sup_Q |k(Q)| \sup_P |h(P)|$. It suffices to show that for any $Q\in G$ we have
\begin{eqnarray*}
  \sum_{P\in  \P'_{G,k}: |I_P| \le |I_Q|/C} |I_Q| |I_P| \overline{k(Q)} h(P) \<\phi_P, \phi_{3,Q}\> 
&\lesssim& 2^{-k} \Lambda (\frac{|I_Q|}{|I_G|})^2 |I_G|
\end{eqnarray*} 
We decompose $\phi_{3,Q}=\phi_{3,Q}1_{2^k I_G} + \phi_{3,Q} (1-1_{2^k I_G})$. Since $k(Q)h(Q)\lesssim \Lambda$ and since $I_P$'s are essentially pairwise disjoint and contained in $2^{k+3}I_G\setminus 2^k I_G $,  the left hand side of the above display is bounded above by
\begin{eqnarray*}
&\lesssim&  \Lambda \sum_{P} |I_P| |I_Q| \Big(\|\phi_{P}1_{2^k I_G} \|_1  \|\phi_{3,Q}\|_\infty + \|\phi_{P}\|_1  \|\phi_{3,Q}1_{\R\setminus 2^k I_G}\|_\infty  \Big)\\
&\lesssim&  2^{-2k}\Lambda   \sum_{P: \ |I_P| \le |I_Q|} |I_P| |I_Q| \Big((\frac {|I_P|} {2^k |I_G|})^2 \frac 1{|I_Q|}  + \frac 1 {|I_Q|} (\frac {|I_Q|} {2^k |I_G|})^2  \Big)\\
&\lesssim&  2^{-2k} \Lambda  (\frac {|I_Q|}{|I_G|})^2  \sum_{P: \ I_P\subset 2^{k+3}I_G} |I_P| \quad \lesssim \quad 2^{-k} \Lambda   (\frac{|I_Q|}{|I_G|})^2 |I_G| \ \ . \quad  \text{\endproof}
\end{eqnarray*}

\subsection{Outer $L^\infty$ estimate for $K d$}
When $d(P)=\<f_3,\widetilde \phi_P\>$, we have
\begin{lemma}\label{l.K-infty}  For any $\epsilon>0$ we gave
\begin{eqnarray*}
\|K  d\|_{\mathcal L^\infty(\Q, S_2^{[3]},\sigma)} &\lesssim_{\epsilon,N}& \sup_{Q\in \Q} (\frac 1{|I_Q|}\int \widetilde \chi_{I_Q}^N |f_3|^{1+\epsilon})^{1/(1+\epsilon)} \ \ .
\end{eqnarray*}
\end{lemma}

\proof
Let $E\subset \Q$ be $3$-lacunary and define $\P_E$ as in Observation~\ref{o.decoupleBC}. It follows that  $K d(Q) = \<T_{\P_E}^\ast f_3,\phi_{3,Q}\>$ where $T_{\P_E}^\ast$ is the adjoint of
\begin{eqnarray*}
T_{\P_E} f &:=& \sum_j \sum_{P\in \P_E} |I_P| \<f,\phi_P\> \widetilde \phi_{P,j}  d_j  \ \ .
\end{eqnarray*}
Therefore, similar to the $L^{2,\infty}$ case of Lemma~\ref{l.Kadj-BC},  it suffices to show that
\begin{eqnarray*}
\frac 1 {|I|} \int_I  |(T_{\P_E}^\ast  f_3)(x)| dx &\lesssim& (\frac 1 {|I|}\int \widetilde \chi_I^N |f_3|^{1+\epsilon})^{1/(1+\epsilon)} \ \ ,
\end{eqnarray*}
for every interval $I$ of the same length as $I_E$.

Now, it is clear that $\P_E$ is an overlapping generating subset of $\P$. It follows that for any $x$ only one $j$ and one scale of $\P_E$ would contribute to the defining summation of $T_{\P_E}f$. Thus, $T_{\P_E}f$ is controlled by the maximal function of $\sum_{P\in \P_E} |I_P|\<f,\phi_P\>\phi_P$ and so by Calderon--Zygmund theory, $T_{\P_E}$ is bounded on $L^p(\R)$ for any $1<p<\infty$. By duality $T_{\P_E}^\ast$ is also bounded on $L^p(\R)$. Furthermore,  we also have
\begin{eqnarray}\label{e.T-pointwise}
\<T^\ast _{\P_E}g, h\> &\lesssim& \sum_{P\in \P_E}  \frac 1{|I_P|}\<|g|,\widetilde \chi_P^{N}\> \<\widetilde \chi_{I_P}(x)^{N}, |h|\> \ \ , \ \  N>0  \ \  , 
\end{eqnarray}
from there we obtain a pointwise estimate for $T^\ast_{\P_E}g(x)$ which holds for a.e. $x$.

Decompose $f_3=f_31_{7I} + f_31_{(7I)^c}$. By boundedness of $T^\ast_{\P_E}$ and H\"older inequality, the contribution of $f_31_{7I}$ could be easily controlled. For the contribution of $f_31_{(7I)^c}$, let $\P_{E,1} = \{P\in \P_E: I_P \subset 5I\}$, and $\P_{E,2}=\P_E\setminus \P_{E,1}$. In $\P_{E,1}$ clearly $\widetilde \chi_{I_P}\lesssim \widetilde \chi_I$, thus using the resulting pointwise estimate resulting from \eqref{e.T-pointwise} we obtain
\begin{eqnarray*}
T^\ast_{\P_{E,1}}(f_31_{(7I)^c})(x) &\lesssim&   \sum_{P\in \P_{E,1}}|I_P|^{-1} \<|f_3|,  (\frac{|I_P|}{|I|})^{2} \widetilde \chi_I^N\> \widetilde \chi_{I_P}(x)^2 \quad \lesssim \quad \frac 1{|I|} \<|f_3|, \widetilde \chi_I^N\>    \ \ .
\end{eqnarray*}
For $\P_{E,2}$ it is clear that  $\widetilde \chi_{I_P}(x)\widetilde \chi_{I_P}(y) \le (|x-y|/|I_P|)^{-2} \lesssim \widetilde\chi_I(y)$ for every $x\in I$ and $y\in (7I)^c$. It follows that
\begin{eqnarray*}
\|T_{\P_{E,2}}^\ast f_3\|_{L^1(I)} &\lesssim&  
\int |f_3(y)|  \widetilde \chi_I(y) ^N \sum_{P}  \frac 1{|I_P|} \widetilde \chi_{I_P}(y)^2 (\frac{|I_P|}{|I|})^{-2}   \|\widetilde \chi_{I_P}^2\|_{L^1(I)} dy\\
&\lesssim& \int |f_3(y)|  \widetilde \chi_I(y) ^N \sum_{P}   \widetilde \chi_{I_P}(y)^2 (\frac{|I_P|}{|I|})^{2}   dy\\
&\lesssim& \int |f_3(y)|  \widetilde \chi_I(y) ^N dy\ \ . \quad \text{\endproof}
\end{eqnarray*}

\subsection{Proof of Theorem~\ref{t.BC-estimate}}\label{bcmain-proof-section}

\subsubsection{The basic range} 

Let  $q=\min(2,r/2) = r/2 \in (1,2)$. Let $a_3=\overline {K d}$.  Assume $q_1,q_2,q_3\ge 1$ such that $\sum 1/q_j=1$. By classical H\"older, we obtain
\begin{eqnarray*}
S_1(a_1a_2 a_3)  &\lesssim&   S_2^{[1]}(a_1) S_2^{[2]}(a_2) S_2^{[3]}(a_3) \ \ .
\end{eqnarray*}
Using outer Radon-Nikodym/H\"older inequalities, it follows that 
\begin{eqnarray}
\label{e.trilinearBC1}  \Lambda_{\P,\Q}    &\lesssim&      \prod_{j=1,2,3}  \|a_j\|_{\mathcal L^{q_j}(\Q, S_{2}^{[j]},\sigma)}  \  \ ,
\end{eqnarray}
provided that $\sum 1/q_j = 1$ and $q_3 >q'$. 
If   $q_1,q_2>2$ then via Lemma~\ref{l.K-BC} and Carleson embeddings  (Theorem~\ref{t.std-embed}, Theorem~\ref{t.var-embed}) we obtain
$\Lambda_{\P,\Q} \lesssim \prod_j \|f_j\|_{q_j}$.

\subsubsection{Extending the range}\label{s.extendrange-BC}

To extend the range, we will use restricted weak-type interpolation, following \cite{mtt2002}. Let $\alpha=(\alpha_1,\alpha_2,\alpha_3)$ be such that $\sum_j \alpha_j=1$ and $\alpha_j\ge 0$. Then we say that a trilinear form $\Lambda(f_1,f_2,f_3)$ satisfies restricted weak-type estimates with exponents $\alpha$  if the following holds: there exists $j_0 \in \{1,2,3\}$ such that for every $F_1, F_2, F_3\subset \R$ finite Lebesgue measures we could find $B\subset F_{j_0}$ with less than half of the measure, so that $|\Lambda_{\P,\Q}(f_1,f_2,f_3)| \lesssim |F_1|^{\alpha_1}|F_2|^{\alpha_2}|F_3|^{\alpha_3}$ whenever $|f_j| \le 1_{F_j}$ for every $j$ and furthermore $|f_{j_0}|\le 1_{F_{j_0}-B}$. When exactly one of the index is negative, say $\alpha_{k}<0$, we  say that $\Lambda(f_1,f_2,f_3)$ satisfies restricted weak-type estimates with exponents $\alpha$ if the previous claim holds with $j_0=k$.

From \cite{mtt2002}, if $T(f_1,f_2)$ is a bilinear operator such that $\<T(f_1,f_2),f_3\>$ satisfies restricted weak-type estimates for a finite collection $V$ of triples then 
$$\|T(f_1,f_2)\|_{p_3} \lesssim \|f_1\|_{L^{p_1}(\R)} \|f_2\|_{L^{p_2}(\R)}  $$
for any triples of exponents $0<p_1,p_2,p_3\le \infty$ such that $(1/p_1,1/p_2,1/p_3')$ is in the interior of the convex hull of $V$.

Thus, to show Theorem~\ref{t.BC-estimate}, it suffices to show that $\Lambda_{\P,\Q}(f_1,f_2,f_3)$ satisfies Let $\H\subset\{\alpha_1+\alpha_2+\alpha_3=1\}$  with vertices:
\begin{eqnarray*}
H_1 (\frac 1 2, -\frac 1 2,1) & , & H_2(-\frac 1 2, \frac 1 2, 1)  \ \ ,\\
H_3(\frac 1 2 - \frac 2 {q'}, \frac 12 + \frac 1 {q'}, \frac 1{q'})  & , &  H_4 (\frac 1 2, \frac 1 2 + \frac 1 {q'}, -\frac 1{q'})  \ \ , \\
H_5(\frac 1 2 +\frac 1{q'}, \frac 12, -\frac 1{q'}) & , &  H_6 (\frac1 2 + \frac 1{q'}, \frac 12 - \frac 2 {q'}, \frac 1 {q'})    \ \ .
\end{eqnarray*}

Our proof below will be fairly symmetric accross the vertices, so we will show the claim only for $H_2$ and $H_3$.
For any $\epsilon>0$ and $q_1,q_2>2$, $q_3\ge q'+\epsilon$, it follows from \eqref{e.trilinearBC1}, Carleson embeddings, and convexity that
\begin{eqnarray}
\label{e.trilinearBC2}  \Lambda_{\P,\Q}  &\lesssim& \|a_1\|_\infty^{1-\frac 2{q_1}}  \|a_2\|_\infty^{1-\frac 2 {q_2}} \|a_3\|_\infty^{1-\frac {q'+\epsilon}{q_3}}  \prod_{j=1,2,3}|F_j|^{1/q_j} \ \ .
\end{eqnarray}
By Lemma~\ref{l.s2lac},  for $j=1,2$   we have
\begin{eqnarray}\label{e.aj-infty}
  \|a_j\|_{\mathcal L^\infty(\Q, S_2^{[j]}, \sigma)} 
&\lesssim_N&   \sup_{Q\in \Q} \frac 1{|I_Q|}\int \widetilde \chi_{I_Q}(x)^N |f_j(x)|dx \ \ .
\end{eqnarray}
For $\|a_3\|_{\infty}$ we will use Lemma~\ref{l.K-infty} and obtain
\begin{eqnarray}\label{e.a3-infty}
  \|a_3\|_{\mathcal L^\infty(\Q, S_2^{[3]}, \sigma)} 
&\lesssim_N&   \sup_{Q\in \Q} (\frac 1{|I_Q|}\int \widetilde \chi_{I_Q}(x)^N |f_3(x)|^{1+\epsilon}dx)^{\frac 1 {1+\epsilon}} \ \ .
\end{eqnarray}

We now consider neightborhood of   $H_1$ and $H_2$.  Near these vertices we have $\alpha_1<0$, so we may choose $G_1=F_1\setminus B$ where  $B = \bigcup_{j=2}^3 \{M(1_{F_j}) > C|F_j|/|F_1|\}$ and $C$ is large enough to ensure $|B|<|F_1|/2$.  

Without loss of generality assume that $2^m\le 1+\dist(I_Q, B^c)/|I_Q|<2^{m+1}$  for some $m\ge 0$ integer,
provided that we could obtain extra decaying factors. 

From \eqref{e.aj-infty} we obtain $\|a_2\|_{\mathcal L^\infty}  \lesssim 2^m \sup_{x\in B^c} M1_{F_2}(x) \lesssim 2^m |F_2|/|F_1|$
and similarly using \eqref{e.a3-infty}  we have $\|a_3\|_{\mathcal L^\infty}  \lesssim_\epsilon  (2^m |F_2|/|F_1|)^{1/(1+\epsilon)}$.
It follows from \eqref{e.trilinearBC2} that
\begin{eqnarray*} 
\Lambda_{\P,\Q} &\lesssim&  2^{-m}  |F_1|^{\alpha_1}|F_2|^{\alpha_2}|F_3|^{\alpha_3}
\end{eqnarray*}
$$\alpha_1 = \frac 1{q_1} + \frac 2{q_2} +\frac {q'+\epsilon}{q_3(1+\epsilon)}-1-\frac 1{1+\epsilon} \ \ , \ \ \alpha_2 = 1 -\frac 1 {q_2} \ \ , \ \ \alpha_3=\frac 1{1+\epsilon} - \frac {q'+\epsilon}{q_3(1+\epsilon)} + \frac 1{q_3} \ \ .$$
Thus, letting $(\epsilon,q_1,q_3)$ close to $(0,2,q')$, we can make $\alpha$ arbitrarily close to $H_3$. Similarly, letting $(\epsilon, q_1,q_3)$ close to $(0,2,\infty)$ we can make $\alpha$ arbitrarily close to $H_2$. This completes the proof of Theorem~\ref{t.BC-estimate}.

\section{Estimates for $CC$ model operators}\label{s.CC-estimate}

\begin{theorem}\label{t.CC-estimate} Let $T$ be a $CC$ model operator. Then $\|T\|_{L^{q_1}\times L^{q_2} \to L^{q_3}} < \infty$
for every $1/q=1/q_1+1/q_2$ such that $\frac{2r}{3r-4}<q_1,q_2 < \infty$ and $q_3> \frac {r'}2$.
\end{theorem}
For simplicity, we will assume that $\omega_{j,P,upper}$ are finite intervals and $\omega_{j,P,lower}$ are halflines for $j=1,2$;  the other cases are either symmetric or could be reduced to this situation. We will use the same set up for outer measures space in Section~\ref{s.outerP} for both $\P_1$ and $\P_2$. For convenience of notation, let $\mu_1$ and $\mu_2$ be the corresponding outer measures.

For every $x$, $P\in \P_1$, and $P'\in \P_2$,   there is at most one $1\le j  \le K(x)$ such that  
$$N_{j-1} \in \omega_{1,P,lower} \cap \omega_{2,P',lower} \quad , \quad N_j \in \omega_{1,P,upper} \cap \omega_{2,P',upper}  \ \  . $$
Let $d_{P,P'}(x) = d_j(x)$ if such $j$ exists,  and zero otherwise; define $d_P$ and $d_{P'}$ similarly.

Fix $f_3$ Schwartz on $\R$ and let $a_1(P) = \<f_1,\phi_{1,P}\>$ and $a_2(P') = \<f_2, \phi_{2,P'}\>$, and
\begin{eqnarray*} 
Kf(P) &:=& \sum_{P'  \in  \P_2: |I_P| \ge 16 |I_{P'}|} |I_{P'}| f(P') \<\phi_{1,P}   \phi_{2,P'} d_{P,P'}, f_3\> \\
K^\ast g(P') &:=& \sum_{P   \in  \P_1: |I_P| \ge 16|I_{P'}|} |I_{P}| g(P) \<\phi_{1,P}   \phi_{2,P'} d_{P,P'}, f_3\> \ \ .
\end{eqnarray*}
Clearly, $\<T(f_1,f_2),f_3\>=   \sum_{P\in \P_1} |I_P| a_1(P)  (Ka_2)(P)$. To prove Theorem~\ref{t.CC-estimate} we first establish outer $L^p$ estimates for $K$. 

For convenience of notation, assume that $f_3$ is supported on a fixed set $F_3$. Let $b_P(x) :=  f_3 d_j T_{j,|I_P|/16}f$ if there exists (a unique) $j$ such that $N_{j-1}\in \omega_{1,P,lower}$ and $N_j\in \omega_{1,P,upper}$, and let $b_P(x)=0$ otherwise. Also, define
\begin{eqnarray*}
M_N(F_3) &:=& M_N(\P_1,1_{F_3})^{1/p_1} M_N(\P_2,1_{F_3})^{1/p_2} \ \ .
\end{eqnarray*}
(Recall the definition of $M_N$ from \eqref{e.tileHLmax}.)

\subsection{Outer $L^p$ estimates for $K$.}
The following Lemma is the main estimate of the current section, and we will always assume that $1/p_1+1/p_2+1/p_3=1$ and $2<p_1,p_2<2r/(4-r)$,  and $\infty> p_3> r/(r-2)$.  All implicit constants may depend on these exponents.

\begin{lemma}\label{l.Kf-CC}  It holds that
\begin{eqnarray*}
\|K f\|_{\mathcal L^{p_1'}(\P_1, S_{1,overlap} + S_{2,lac}, \mu_1)} 
&\lesssim_{N}&  M_N(F_3)  \|f\|_{\mathcal L^{p_2}(\P_2,S_{2,lac},\mu_2)} \|f_3\|_{p_3} \ \ .
\end{eqnarray*}
\end{lemma}

The proof of  Lemma~\ref{l.Kf-CC} consists of two parts. The factor $M_N(F_3)$ should be thought of an estimate for the norm of $K: \mathcal L^{p_2}(\P_2) \times L^{p_3}(F_3) \to \mathcal L^{p_1}(\P_1)$, capturing the interaction of  $supp (f)$ (i.e. $\P_2$), $supp(Kf)$ (i.e. $\P_1$) and $supp(f_3)$. Thus, by interpolation (Lemma~\ref{l.outer-multi-interpolation}) it suffices to show the weak-type estimates (note that the constant $M_N(F_3)$ could be naturally absorbed inside the (outer)measures on the right hand side), and we will treat the contribution of $S_{1,overlap}$ in  Section~\ref{s.Kf-overlap-CC} and the contribution of $S_{2,lac}$ in Section~\ref{s.Kf-lac-CC}.  

We first fix some notations. For each $j$ and any $f:\P_2 \to \C$ and $\alpha>0$ let
\begin{eqnarray*}
T_{j,\alpha}(f) &=& \sum_{P'\in \P_2: |I_{P'}| \le \alpha} |I_{P'}| f(P')  \widetilde \phi_{2,P',j}  \ \ , 
\ \ T^\ast_j f \quad = \quad \sup_{\alpha} |T_{j,\alpha}f |\ \ .
\end{eqnarray*}
For every $E\subset \P_1$ let $\P_2(E)$ contains all $P'\in \P_2$ such that  for some $P \in E$ we have $|I_{P'}|\le |I_{P}|/16$ and $\frac 5 4\omega_{1,P_1,upper}\cap \frac 5 4 \omega_{2,P',upper}\ne\emptyset$.

\subsubsection{The  overlapping setting}\label{s.Kf-overlap-CC} In this section we prove that
\begin{eqnarray*}
\|Kf\|_{\mathcal L^{p_1',\infty}(\P_1,S_{1,overlap},\mu_1)} &\lesssim_N& M_N(\P_1, 1_{F_3})^{1/p_1} \|f\|_{\mathcal L^{p_2}(\P_2,S_{2,lac}, \mu_2)}  \|f_3\|_{L^{p_3}(\R)} \ \ .
\end{eqnarray*}
Using H\"older inequality, it follows from Lemma~\ref{l.Kf-overlap-CC-infty} that
\begin{eqnarray*}
&& \|Kf\|_{\mathcal L^{\infty}(\P_1,S_{1,overlap}, \mu_1)} \\
&\lesssim_N& M_N(\P_1, 1_{F_3})^{1/p_1} \sup_{R\in \overline{\P_1}} (\frac 1 {|I_R|} \int \widetilde \chi_{I_R}^N \sum_{j: N_j \in \widetilde \omega_R} |f_3 d_jT^\ast_j(f)|^{p_1'} dx)^{1/p_1'} 
\end{eqnarray*}
Therefore, it follows from Lemma~\ref{l.ms-select} that
\begin{eqnarray*}
 \|Kf\|_{\mathcal L^{p_1',\infty}(\P_1,S_{1,overlap}, \mu_1)} 
&\lesssim_N& M_N(\P_1, 1_{F_3})^{1/p_1} \| f_3 (\sum_{1\le j\le K} |d_jT^\ast_j(f)|^{p_1'})^{1/p_1'}\|_{p_1}
\end{eqnarray*}
so  using $\sum_j |d_j|^{r/(r-2)} = O(1)$ it follows that
\begin{eqnarray*}
 \|Kf\|_{\mathcal L^{p_1',\infty}(\P_1,S_{1,overlap}, \mu_1)} 
&\lesssim& M_N(\P_1, 1_{F_3})^{1/p_1}  \|f_3\|_{p_3}  \| (\sum_{1\le j \le K}  |T^\ast_j(f)|^s)^{1/s}\|_{L^{p_2}(F_3)}  \ \ ,
\end{eqnarray*}
where $s>0$ such that $1/p_1' \le 1/s + (r-2)/r$. Since $1/p'_1=1/p_2+1/p_3 < 1/p_2 + (r-2)/r$, we could choose $s$ such that $s>p_2$ (which is larger than $2$). The desired estimate now follows from Theorem~\ref{t.varCarl}.   

\begin{lemma}\label{l.Kf-overlap-CC-infty}    Uniform over $\P\subset \P_1$ it holds that
\begin{eqnarray*}
\|K(f)\|_{\mathcal L^{\infty}(\P, S_{1,overlap}, \mu_1)} &\lesssim_N&   m_{\infty,N}(\P, (f_3 d_j T^\ast_j f)) \ \ .\end{eqnarray*}
\end{lemma}

\proof   This follows from a simple adaptation of the proof of Lemma~\ref{l.var-infty-ovl} and the fact that for every $j$ and every $P$ we have $|T_{j,|I_P|/16}f| \le T^\ast_j f$. \endproof

\subsubsection{The lacunary setting} \label{s.Kf-lac-CC}  In this section we prove that
\begin{eqnarray*}
\|Kf\|_{\mathcal L^{p'_1,\infty}(\P_1,S_{2,lac},\mu_1)} &\lesssim_N& M_N(F_3)  \|f\|_{\mathcal L^{p_2}(\P_2,S_{2,lac}, \mu_2)}  \|f_3\|_{L^{p_3}(\R)} \ \ .
\end{eqnarray*}

Fix $\lambda>0$. Without loss of generality assume that $\|Kf\|_{\mathcal L^\infty(S_{2,lac})} \le 2\lambda$. 

Apply a variant of the selection argument in the proof of Theorem~\ref{t.std-embed} we may find $A\subset \P_1$ such that
$\|Kf\|_{\mathcal L^\infty(\P_1\setminus A, S_{2,lac},\mu_1)} \le \lambda$
and a strongly disjoint collection of generating subsets $(E_m)$ covering $A$  such that 
\begin{eqnarray*}
\lambda^2 \mu_1(A) \quad \lesssim \quad \lambda^2 \sum_m |I_{E_m}|  
&\lesssim&  M:= \sum_{P\in E_m} |I_P| |Kf(P)|^2 \ \ .
\end{eqnarray*}
Let  $g(P) = \overline{Kf(P)}$ for $P\in \Q_1 :=\bigcup E_m$ and $g(P)=0$ otherwise.  We obtain
\begin{eqnarray*}
M &=& \sum_{P\in \Q_1} |I_{P}| Kf(P)  g(P)  \quad = \quad  \sum_{P'\in \P_2} |I_{P'}| f(P') (K^\ast g)(P')  \ \ , 
\end{eqnarray*}
Using Lemma~\ref{l.Kdual-CC} and the assumption that $\|Kf\|_{\mathcal L^\infty(\Q_1, S_{2,lac}, \mu_1)} \le 2\lambda$ we obtain
\begin{eqnarray*}
\|K^\ast g\|_{\mathcal L^{p'_2}(\P_2,S_{2,lac} + S_{1,overlap},\mu_1)}  &\lesssim_N& M_N(F_3) \lambda (\sum_{m} |I_{E_m}|)^{1/p_1} \|f_3\|_{p_3}  \ \ ,
\end{eqnarray*}
Using outer Radon-Nikodym/H\"older, it follows that
\begin{eqnarray*}
 M &\lesssim& \|f\|_{\mathcal L^{p_2}(\P_2,S_{2,lac},\mu_2)}  \|K^\ast g\|_{\mathcal L^{p'_2}(\P_2,S_{2,lac} + S_{1,overlap},\mu_1)} \\
&\lesssim&  M_N(F_3) \lambda \|f\|_{\mathcal L^{p_2}(\P_2,S_{2,lac},\mu_2)} \|f_3\|_{p_3} (\lambda^{-2}M)^{1/p_1} \ \ ,
\end{eqnarray*}
from this the desired estimates for $M$ (and hence for $\mu_1(A)$) easily follow.

\begin{lemma}\label{l.Kdual-CC} It holds that (with $\|g\|_\infty =   \|g\|_{\mathcal L^{\infty}(\Q_1, S_{2,lac},\mu_1)}$)
\begin{eqnarray*}
\|K^\ast g\|_{\mathcal L^{p_2'}(\P_2,  S_{1,overlap}+S_{2,lac},\mu_2)} 
 &\lesssim_N&  M_N(F_3) \|f_3\|_{p_3}  \|g\|_{\infty}(\sum_m |I_{E_m}|)^{1/p_1}   \ \ .
\end{eqnarray*}
\end{lemma}

\proof  By simple modifications of the argument in Section~\ref{s.Kf-overlap-CC} together with Theorem~\ref{t.varCarl} part (i), we obtain
\begin{eqnarray*}
\|K^\ast g\|_{\mathcal L^{p_2'}(\P_2,  S_{1,overlap},\mu_2)} &\lesssim&  M_N(\P_2,1_{F_3})^{1/p_2} \|f_3\|_{p_3}  \|V_{p_1}(g)1_{F_3}\|_{p_1}\\
&\lesssim&  M_N(\P_1,1_{F_3})^{1/p_1} M_N(\P_2,1_{F_3})^{1/p_2} \|f_3\|_{p_3} \|g\|_{\mathcal L^{p_1}(S_{2,lac})}  \ \ ,
\end{eqnarray*}
so it remains to  consider   the contribution of  $S_{2,lac}$, and by interpolation it suffices to consider weak-type estimates.

Fix $\alpha>0$. Without loss of generality we may assume that  $\|K^\ast g\|_{\mathcal L^\infty(S_{2,lac})} \le 2\alpha$.
By a standard argument, we may find   $B\subset \P_2$ with $\|K^\ast g\|_{\mathcal L^\infty(B, S_{2,lac},\mu_2)} \le \lambda$  and a collection of strongly disjoint lacunary generating sets $(G_n)$ covering $B$ such that
\begin{eqnarray*}
\alpha^2 \mu_2(B) \quad \lesssim \quad  \alpha^2 \sum_n |I_{G_n}| &\lesssim& N:=\sum_{n} \sum_{P'\in G_n} |I_{P'}| |K^\ast g(P')|^2 \ \ .
\end{eqnarray*}
Let $h(P') = \overline{K^\ast g(P')}$ for  $P'\in \Q_2 = \bigcup G_n$ and let $h(P')=0$ otherwise. We obtain
\begin{eqnarray*}
N &=&  \sum_{\stackrel{P \in \Q_1 \ , \ P'\in \Q_2}{|I_P| \ge |I_{P'}|/16}} |I_P| |I_{P'}| h(P') g(P) \<\phi_{1,P}\phi_{2,P'} d_{P,P'}, f_3\> \ \ ,
\end{eqnarray*}
thus it suffices to show that
\begin{eqnarray*}N  &\lesssim&  M_N(F_3) \|f_3\|_{p_3} \|h\|_{\mathcal L^\infty(S_{2,lac})} \|g\|_{\mathcal L^\infty(S_{2,lac})}   (\sum_m |I_{E_m}|)^{\frac 1{p_1}} (\sum_{n} |I_{G_n}|)^{\frac 1{p_2}} \ \ .
\end{eqnarray*}
Note that since $\Q_1$ and $\Q_2$ are unions of strongly disjoint lacunary sets, all overlapping sets are essentially a collection of spatially disjoint tiles, therefore $S_{1,overlap}\lesssim S_{2,lac}$ in each of them. This observation will be used implicitly below.

We first show that  the unconstrained double sum $N_{free}$ over $P,P'$ (where there is no constraint between $P$ and $P'$) satisfies the desired estimate.   Indeed, since $2<p_1, p_2<2r/(3r-4)$ and $1/p_1+1/p_2=1-1/p_3 > 2/r$ we may find $s,t>2$ such that $2/r=1/s+1/t$ and $t>p_1$ and $s>p_2$. Using Theorem~\ref{t.varCarl}   we have
\begin{eqnarray*}
N_{free}  &\lesssim_N& M_N(\P_1, 1_{F_3})^{1/p_1} \|g\|_{\mathcal L^{p_1}(S_{2,lac})} \|(\sum_{j}|\sum_{P \in \Q_1}  |I_P| g(P) \widetilde \phi_{2,P',j}d_j f_3|^{t'})^{1/t'}\|_{p'_1} \\
&\lesssim& M_N(\P_1, 1_{F_3})^{1/p_1}\|g\|_{\mathcal L^{p_1}(S_{2,lac})} \|h\|_{\mathcal L^{p_2}(S_{2,lac})} M_N(\P_2, 1_{F_3})^{1/p_2} \|f_3\|_{p_3} \ \ .
\end{eqnarray*}

We now consider  diagonal   sums when $|I_{P'}| = C|I_P|$ for some fixed $C\in [2^{-4}, 2^4]$.  The proof below is symmetric in $P$, $P'$, so we will assume $C\le 1$. We say that $P$ is linked to $P'$ if the corresponding summand  is nonzero, clearly that all linked pairs satisfy $\omega_{1,P,upper}\cap \omega_{2,P',upper}\ne \emptyset$. Since these are dyadic intervals and $|\omega_{1,P,upper}|\le |\omega_{2,P',uppper}|$, we obtain $\omega_{1,P,upper} \subset \omega_{2,P',upper}$. Without loss of generality, assume that for some $k$  it holds that $2^k |I_P| \le |I_P|+\dist(I_P, I_{P'}) < 2^{k+1} |I_P|$ for all linked pairs, provided that  we have sufficient decay in the estimates. It follows that for each $P'$ there is at most $O(1)$ linked $P$ and vice versa, and by further dividing if necessary we may assume that exactly one $P$ is linked to exactly one $P'$, and let $P'=F(P)$ and $P=F^{-1}(P')$.  It follows that $h(P')$ and $c(P,P'):=|I_{P'}| \<\phi_{1,P}\phi_{2,P'} d_{P,P'}, f_3\>$ are functions on $\P$. Using outer Radon--Nikodym/H\"older and the triangle inequalities, it follows that the corresponding sum is bounded by
\begin{eqnarray*}
N_{diag,k} &\lesssim&    \|g\|_{\mathcal L^{p_1}(\Q_1,S_{2},\mu_1)} \|h\circ F\|_{\mathcal L^{p_2}(\Q_1, S_{2}, \mu_1)} \|c\|_{\mathcal L^{p_3}(\Q_1, S_\infty, \mu_1)}  \  \  .
\end{eqnarray*}
Now, we note that since $\Q_1$ is an union of strongly disjoint  lacunary sets, all the overlapping generating subsets of $\Q_1$ has $O(1)$ elements, therefore $S_2 \lesssim S_{2,lac}$ on $\Q_1$.  For any generating set $E_2 \subset \Q_2$ it is clear that $F^{-1}(E_2)$ could be covered by $O(1)$  generating sets of $\Q_1$, whose top intervals are contained in some bounded enlargement of $2^k I_{E_2}$.  Therefore $\mu_1(F^{-1}(E_2))\lesssim 2^k \mu_2(E_2)$. Conversely, if $E_1$ is a generating set in $\Q_1$ then $F(E_1)$ could be generously covered by $O(2^k)$ generating sets in $\Q_2$, and the length of the top intervals of these covering sets are comparable to $|I_{E_1}|$. Therefore $S_{2,lac}(h\circ F)(E_1) \lesssim 2^k S_{2,lac}(h)(F(E_1))$. Therefore by pull-back (essentially the same proof as \cite[Proposition 3.2]{dt2014}) we obtain
\begin{eqnarray*}
 \|h\circ F\|_{\mathcal L^{p_2}(\Q_1, S_{2}, \mu_1)} &\lesssim& 2^{2k} \|h\|_{\mathcal L^{p_2}(\Q_2, S_{2,lac}, \mu_2)} \ \ .
\end{eqnarray*}
On the other hand, notice that for any $P\in \Q_1$  we have
\begin{eqnarray*}
c(P,P') &\lesssim_N& 2^{-Nk} \frac 1{|I_P|} \int (\widetilde \chi_{I_P}  \widetilde \chi_{I_{P'}})^N \sup_{j: N_j\in \omega_{1,P,upper}} |d_j| |f_3|
\end{eqnarray*}
Therefore using (perhaps a version of) Lemma~\ref{l.ms-select},  it follows that
\begin{eqnarray*} 
\|c\|_{\mathcal L^{p_3}(\Q_1, S_\infty, \mu_1)} &\lesssim_N& 2^{-Nk} \sup_{R\in \overline{\Q_1}} (\frac{1}{|I_R|} \int_{F_3}(\widetilde \chi_{I_R}\widetilde \chi_{I_{R'}})^{p_3' N})^{1/p_3'}   \|f_3\|_{p_3} \\
&\lesssim_N& 2^{-Nk} M_N(F_3) \|f_3\|_{p_3} \ \ ,
\end{eqnarray*}
here we used $p_3>r/(r-2)$. This leads to the desired estimate for $N_{diag,k}$.

Consequently, for the purpose of proving the desired estimates for $N$ the roles of $\Q_1$ and $\Q_2$ are symmetric, and we may assume that 
$$M_N(\Q_1,F_3)\mu_1(\Q_1) \le M_N(\Q_2,F_3) \mu_2(\Q_2) \ \ .$$ 
It follows from Lemma~\ref{l.Ktilde-sparse-CC} that
\begin{eqnarray*}
&&\|Kh\|_{\mathcal L^{p'_1}(\Q_1, S_{2,lac},\mu_1)}\\ &\lesssim& \big(M_N(\Q_1,F_1)^{1/p_1+1/p_2}\mu_1(\Q_1)^{1/p_2}+M_N(F_3)\mu_2(\Q_2)^{1/p_2}\big) \|h\|_{\mathcal L^\infty(S_{2,lac})} \|f_3\|_{p_3}\\
&\lesssim& M_N(F_3) \mu_2(\Q_2)^{1/p_2}\|h\|_{\mathcal L^\infty(S_{2,lac})} \|f_3\|_{p_3} \ \ .
\end{eqnarray*}
Recall that on $\Q_1$ and $\Q_2$ we have $S_2\lesssim S_{2,lac}$, so using a combination of outer Radon-Nikodym and outer H\"older inequalities, we obtain
\begin{eqnarray*}
N_{offdiag} &=& \sum_{P\in \Q_1}|I_P| g(P) (K h)(P) \\
&\lesssim& \|g\|_{\mathcal L^{p_1}(\Q_1,\S_{2,lac},\mu_1)} \|K h\|_{\mathcal L^{p_1'}(\Q_1,S_{2,lac},\mu_1)}\\
&\lesssim& M_N(F_3) \mu_1(\Q_1)^{1/p_1} \mu_2(\Q_2)^{1/p_2}\|g\|_{L^{\infty}(S_{2,lac})}  \|h\|_{\mathcal L^\infty(S_{2,lac})} \|f_3\|_{p_3}\\
&\lesssim&  M_N(F_3) \  (\sum_{m} |I_{E_m}|)^{1/p_1}  (\sum_n |I_{G_n}|)^{1/p_2}   \|g\|_{L^{\infty}(S_{2,lac})}  \|h\|_{\mathcal L^\infty(S_{2,lac})} \|f_3\|_{p_3} \ \ ,
\end{eqnarray*}
as desired. \endproof

\begin{lemma}\label{l.Ktilde-sparse-CC} Let $f$ be supported on $\Q_2$. Then
\begin{eqnarray*}
&&\|K f\|_{\mathcal L^{p_1'}(\Q_1,S_{2,lac},\mu_1)} \\ 
&\lesssim&  (M_N(\Q_1,F_3)^{1/p_3'} \mu_1(\Q_1)^{1/p_2} + M_N(F_3) \mu_2(\Q_2)^{1/p_2})\|f\|_{\mathcal L^\infty(S_{2,lac})} \|f_3\|_{p_3} \ \ .
\end{eqnarray*}
\end{lemma}

To prove Lemma~\ref{l.Ktilde-sparse-CC}, we first show the following estimate for $\|K f\|_\infty$.

\begin{lemma}\label{l.Kf-lac-CC-infty}   Uniform over $\P\subset \P_1$, it holds for any $s<2$ and $N>0$ that
\begin{eqnarray*}
&&\| K f\|_{\mathcal L^\infty(\P, S_{2,lac},\mu_1)}\\  &\lesssim_{s,N}&   m_{\frac r {r-2},N}(\P,(f_3 d_j)) \|f\|_{\mathcal L^\infty(\P_2,S_{2,lac},\mu_2)}+ m_{s,N}(\P,(f_3 d_j T^\ast_j f)) 
\end{eqnarray*}
\end{lemma}

Below we deduce Lemma~\ref{l.Ktilde-sparse-CC} from Lemma~\ref{l.Kf-lac-CC-infty}.  By interpolation, it suffices to prove weak-type estimates. We may assume  $\|f\|_{\mathcal L^{\infty}(\Q_2,S_{2,lac},\mu_2)}=1$ by scaling. Using Lemma~\ref{l.ms-select} for $s=p_1$, for any $\lambda>0$ we may find $A\subset \Q_1$ such that
$$m_{p_1',N}(\Q_1\setminus A, (f_3d_jT^\ast_j f))  \quad \le \quad \lambda  \ \ ,  \ \ \text{and} $$
$$\lambda \mu_1(A)^{1/p_1'}  \quad \lesssim \quad  M_N(\Q_1,1_{F_3})^{1/p_1} \|(\sum_j |f_3d_j T^\ast_jf|^{p_1'})^{1/p_1'}\|_{p_1'} \ \ .$$
Thus, arguing as in Section~\ref{s.Kf-overlap-CC}  we obtain
\begin{eqnarray*} 
\lambda \mu_1(A)^{1/p_1'}  &\lesssim&  M_N(\Q_1,F_3)^{1/p_1} M_N(\P_2,F_3)^{1/p_2} \|f\|_{\mathcal L^{p_2}(\Q_2,S_{2,lac},\mu_2)}\|f_3\|_{p_3} \\
&\lesssim& M_N( F_3)  \mu_2(\Q_2)^{1/p_2}  \|f_3\|_{p_3} \ \ .
\end{eqnarray*}
Since $p_3>r/(r-2)$, it follows from Lemma~\ref{l.ms-select} that we could find $B\subset\Q_1$ with
$$m_{\frac r{r-2}} (\Q_1\setminus B, (f_3 d_j)) \le \lambda  \ \ \text{and}$$
$$\lambda \mu_1(B)^{1/p_3} \lesssim  M_N(\Q_1,1_{F_3})^{1/p_3'} \|f_3\|_{p_3}  \ \ .$$
Clearly we have $\mu_1(B) \le \mu_1(\Q_1)$, therefore
$$\lambda \mu_1(B)^{1/p_1'} = \lambda \mu_1(B)^{1/p_3} \mu_1(B)^{1/p_2} \lesssim M_N(\Q_1, 1_{F_3})^{1/p_3'}  \|f_3\|_{p_3} \mu_1(\Q_1)^{1/p_2} \ \ ,$$
and this completes our proof of Lemma~\ref{l.Ktilde-sparse-CC}.

\proof[Proof of Lemma~\ref{l.Kf-lac-CC-infty}] 
Let $E\subset \P$ be lacunary. Using sparseness of $\P_2$, it is clear that $\P_2(E)$ is a lacunary subset of $\P_2$. Let $u_j:=f_3d_j T^\ast_j f$ for convenience. Similar to the proof of Lemma~\ref{l.var-infty-lac}, it suffices to show that if $\|g\|_{\mathcal L^\infty(E, S_{2,lac},\mu_1)}=1$ then
\begin{eqnarray*}
&&\frac{1}{|I_E|}\sum_{P\in E} |I_P| Kf(P)g(P) \\
&\lesssim_{s,N}&     m_{\frac r {r-2},N}(\P,(f_3 d_j)) \|f\|_{\mathcal L^\infty(\P_2,S_{2,lac},\mu_2)}+ m_{s,N}(\P,(u_j))  \ \ .
\end{eqnarray*}

Let $\J$ and $b_P$ be defined as usual,  we also start with
\begin{eqnarray}\label{e.Jdecomp-CC-lac}
\sum_{P\in E} |I_P| Kf(P) g(P) &\lesssim& \sum_{J}  \|\sum_{P\in E} |I_P| g(P) \phi_{1,P}b_P\|_{L^1(J)} \quad \le \quad A+ B\ \ ,
\end{eqnarray}
where $A$ denote the contribution of $|I_P|\le 2|J|$ and the rest is in $B$. Using $|b_P(x)|\le \sup_{j: N_j\in \omega_{1,P,upper}}|u_j|$, it is not hard to see that
\begin{eqnarray*}
A &\lesssim& |I_E| m_{\infty,N}(E, (u_j)) \ \ .
\end{eqnarray*}

We now estimate $B$. For convenience, let $G_E(x) = \sum_{P\in E} |I_P| g(P) \phi_{1,P}(x)$ and
$$V^s(g_E)(x) := \sup_{K,n_0<\dots<n_K} (\sum_{j=1}^K |(\Pi_{n_j} -\Pi_{n_{j-1}} )g_E(x) |^s)^{1/s}$$
for $0<s<\infty$, where $\Pi_j$ denotes a suitable Fourier projection on to the relevant frequency scale of $E$ (see also proof of Lemma~\ref{l.var-infty-lac}).   

For each $J$, let $\omega_J = \bigcup_{P\in E: |I_P| \ge 4|J|} \omega_{1,P,upper}$. As in the proof of Lemma~\ref{l.var-infty-lac}  there exists $R \in \overline{\P}$ with $|I_R| \approx |J|$, $\dist(I_R,J)\lesssim |J|$, and $\omega_J \subset \widetilde \omega_R$. By decomposing $T_{j,|I_P|/16}=T_{j,|J|/4} + (T_{j,|I_P|/16} - T_{j,|J|/4})$ we obtain the decomposition  $b_P=\ell_{P,J}+s_{P,J}$, and we will estimate the contribution of each term.

\underline{Contribution of $\ell_{P,J}$}: We  decompose further $\ell_{P,J} = \ell_{P,J,core}+\ell_{P,J,tail}$   by decomposing $f=f_{I_E,core}  + f_{I_E,tail}$ where $f_{I_E,core}$ is the restriction of $f$ to the subset of $\P_2(E)$ containing all $P'$ with $I_{P'} \subset 5I_E$.
By the H\"older inequality, we have
\begin{eqnarray*}
\sum_J \|\sum_{|I_P|\ge 4|J|} |I_P| g(P) \phi_{1,P} \ell_{P,J,core}\|_{L^1(J)}&\lesssim& \sum_J \|B_{J}(g,f_{I_E,core}) (\sum_{j} |f_3d_j|^{\frac r{r-2}})^{\frac{r-2}r}\|_{L^1(J)} \ \ ,
\end{eqnarray*}
$$B_{J}(g,f) :=\big(\sum_{j} \big|\sum_{P, P' \ constraints} |I_P| |I_{P'}| g(P) f(P') \widetilde \phi_{1,P,j} \widetilde \phi_{2,P',j} \big|^{r/2}\big)^{2/r}\ \ ,$$
the constraints in the sum  are $|I_{P'}| \le |I_P|/16$, $|I_{P}|\ge 4|J|$,  $|I_{P'}|\ge |J|/4$. Let $F_{E,core}=\sum_{P'} |I_{P'}| f_{I_E,core}(P') \phi_{2,P'}$ and define $F_{E,tail}$ similarly. Clearly,
\begin{eqnarray*}
B_{J}(g, f_{I_E,core})(x) &\lesssim& M_J(V^{r/2}( G_E, F_{E,core})) \ \ , \  \ \text{($r>2$ is needed for convexity)} \ \ , 
\end{eqnarray*}
where as in \cite{dmt2012} we define the bilinear variation-norm $V^{r/2}(h_1,h_2)$ to be
$$\sup_{K: N_0<\dots<N_K} \Big(\sum_{k=1}^K |\sum_{N_{k-1}<i<j\le N_k} (\Delta_i h_1)(\widetilde \Delta_j h_2)|^{r/2}\Big)^{2/r} \ \ , $$
and $(\Delta_j)_{j\ge 0}$ and $(\widetilde \Delta_j)_{j\ge 0}$ are two suitable families of Littlewood-Paley projections relative to $\xi_E$.   By variation-norm estimates for paraproducts  \cite{dmt2012},  it follows that
\begin{eqnarray}
\nonumber &&\sum_{J} \|\sum_{|I_P|\ge 4|J|} |I_P| \phi_{1,P} \ell_{P,J,core}\|_{L^1(J)}\\
\nonumber &\lesssim&  m_{\frac r{r-2},N}(E,(f_3d_j)) \|M(V^{r/2}(G_E,F_{E,core}))\|_{L^{1,\infty}(3I_E)}\\
\label{e.varpar} &\lesssim&  m_{\frac r{r-2},N}(E,(f_3d_j)) \|G_E\|_{L^{p'}(\R)} \|F_{E,core}\|_{L^p(\R)} \ \ .
\end{eqnarray}
By Lemma~\ref{l.std-embed-tree} we have  $\|F_{E,core}\|_{L^p(\R)}\lesssim |I_E|^{1/p} \|f\|_{L^\infty(\P_2,S_{2,lac},\mu_2)}$ and a similar estimate for $\|G_E\|_{L^{p'}(\R)}$, giving the desired estimate for the contribution of $\ell_{P,J,core}$.

For contribution of $\ell_{P,J,tail}$, notice that for every $x\in J \subset 3I_E$ it holds that
\begin{eqnarray*}
B_J(g, f_{E,tail})(x)
&\lesssim& (\sum_j |\sum_P |I_P|g(P) \widetilde \phi_{1,P,j}|^r)^{1/r} (\sum_{P'} |I_{P'}| |f_{I_E,tail}(P')\phi_{2,P'}|) \\
&\lesssim_N& M_J(V^{r}(G_E)) \sup_{P': I_{P'}\cap 4I_E=\emptyset, |I_{P'}|\le |I_E|/16}|f(P')|  \widetilde \chi_{I_{P'}}(c(I_E))^{N+5}  \ \ .
\end{eqnarray*}
It follows that 
\begin{eqnarray}
\nonumber &&\sum_J \|\sum_{P\in E: |I_P| \ge 4|J|} |I_P| g(P) \phi_{1,P} \ell_{P,J,tail}\|_{L^1(J)} \lesssim\\
\label{e.ellPJtail} &\lesssim&  m_{\frac r{r-2},N}(E,(f_3d_j)) \|M(V^r(G_E))\|_{L^{1,\infty}(3I_E)}   \sup_{P'} |f(P')|\widetilde \chi_{I_{P'}}(c(I_E))^{N+5}\ \ ,
\end{eqnarray}
which implies the desired estimate for the contribution of $\ell_{P,J,tail}$ via   the continuous L\'epingle inequality (see e.g. \cite{bourgain1989, jsw2008}),. 

\underline{Contribution of $s_{P,J}$}: Since $|I_P|\ge 4|J|\ge 16|I_{P'}|$, we could remove the constraint $|I_{P'}|\le |I_P|/16$ in $s_{P,J}$. 
Using  the  H\"older inequality, it follows that
\begin{eqnarray*}
&& \sum_J \|\sum_{P\in E: |I_P|\ge 4|J|} |I_P| g(P) \phi_{1,P} s_{P,J}\|_{L^1(J)} \lesssim \\
&\lesssim& \sum_J \Big\|(\sum_j |\sum_{|I_P|\ge 4|J|} |I_P| g(P)\widetilde \phi_{1,P,j}|^{s'})^{1/s'}  (\sum_{j:N_j\in \omega_J} |u_j|^{s})^{1/s} \Big\|_{L^1(J)} \\
&\lesssim& \sum_J  |J| M_J(V^{s'}(G_E)) m_{s,N}(E, (u_j)) \ \ .
\end{eqnarray*}
The desired estimate then follows from the continuous   L\'epingle inequality. \endproof

\subsection{Proof of Theorem~\ref{t.CC-estimate}}

\subsubsection{The basic range}
 
We first show the range $2<q_1,q_2<2r/(4-r)$ and $1<q_2<r/2$ of Theorem~\ref{t.CC-estimate}. Note that now $2r/(3r-4)<q_1'<2$ and $q_3'>r/(r-2)$. Using outer Radon Nikodym/H\"older inequalities and  Lemma~\ref{l.Kf-CC}, it follows that
\begin{eqnarray}
\nonumber \<T_{CC}(f_1,f_2),f_3\> &=&  \Lambda_{\P_1,\P_2} := \sum_{P\in \P_1} |I_P| a_1(P) K(a_2)(P) \\
\label{e.trilinearCC}&\lesssim& M_N(F_3)  \|a_1\|_{\mathcal L^{q_1}}  \|a_2\|_{\mathcal L^{q_2}}  \|f_3\|_{L^{q_3'}(\R)}
\end{eqnarray}
where in the above display (and in subsequent displays) the outer norm for $a_j$ is over $(\P_j,S_{2,lac},\mu_j)$ and the outer norm for $Ka_2$ is over $(\P_1,S_{1,overlap} + S_{2,lac},\mu_1)$. Thus the desired claim follows from  generalized Carleson embeddings and duality, and the trivial bound $M_N(F_3)\lesssim 1$.

\subsubsection{Extending the range}
Let ${\bf M} \subset \{\alpha_1 + \alpha_2 +\alpha_3=1\}$ have the following vertices
$$M_1(\frac{ 3r-4}{2r} , \ \frac 12, \  -\frac {r-2}r) \quad , \quad M_2(\frac 12,\  \frac{ 3r-4}{2r} , \  -\frac {r-2}r)  $$
$$M_3(0, \frac{ 3r-4}{2r}, \frac{4-r}{2r}) \quad , \quad M_4(0,0,1) \quad ,\quad M_5(\frac{ 3r-4}{2r}, 0, \frac{4-r}{2r}) \ \ . $$
Similar to Section~\ref{s.extendrange-BC}, to show Theorem~\ref{t.CC-estimate} it suffices to prove restricted-weak type estimates for one $\alpha=(\alpha_1,\alpha_2,\alpha_3)$ in any neighborhood of any given vertex of $\M$. 
Our starting point will be \eqref{e.trilinearCC}, where $\P_1$ and $\P_2$ are symmetric, thus it suffices to consider $M_1, M_3, M_4$.  Below let $F_1, F_2, F_3$ have finite positive Lebesgue measures. Let $C_0>0$ be a finte large absolute constant such that $\|M(f)\|_{1,\infty} < \frac{C_0}{4}\|f\|_1$ in the maximal inequality.

In the following, let $2<p_1,p_2 <2r/(4-r)$ and $p_3>\frac {r} {r-2}$ with $\sum 1/p_j=1$. 

\underline{Near $M_1$.} Let $B= \bigcup_{j=1,2} \{M(1_{F_j})>C_0|F_j|/|F_3|\}$, clearly $|B|<|F_3|/2$. Let $G_3=F_3\setminus B$. We may assume that for some $k_1,k_2\ge 0$ it holds that  \begin{eqnarray}
\label{e.P1decay}
2^{k_1}\le 1 + \frac{\dist(I_P,B^c)}{|I_P|} < 2^{k_1+1} |I_P| \ \ , \ \ \forall \  P\in \P_1
\end{eqnarray}
and similarly $2^{k_2} \le 1+ \frac{\dist(I_{P'}, B^c)}{|I_{P'}|} < 2^{k_2+1}$ for every $P'\in \P_2$, provided that we have sufficient decay in the estimate. Let $|f_1|\le 1_{F_1}$, $|f_2| \le 1_{F_2}$, and $|f_3|\le 1_{F_3 - B}$. 

By convexity, it follows from \eqref{e.trilinearCC} and  Carleson embeddings that
\begin{eqnarray*}
\Lambda_{\P_1,\P_2} &\lesssim&  |F_3|^{\frac 1{p_3}} \prod_{j=1,2}  M_N(\P_j, 1_{F_3\setminus B})^{\frac 1{p_j}} \|a_j\|_{\infty}^{1-\frac 2{p_j}} \|a_j\|_{2,\infty}^{\frac 2 {p_j}}\\
&\lesssim& 2^{-N(k_1+k_2)} |F_1|^{1-\frac 1 {p_1}} |F_2|^{1-\frac 1{p_2}} |F_3|^{\frac  1{p_3} - (1-\frac 2{p_1})  - (1-\frac 2{p_2})} 
\end{eqnarray*}
 thus $\Lambda_{\P_1,\P_2}$ satisfies restricted weak-type estimate for $\alpha=(1/p_1', 1/p_2', -1/p_3)$, which could be made arbitrarily close to $M_1$.

\underline{Near $M_3$.}   Let $B=   \{M(1_{F_3})>C_0|F_3|/|F_1|\}$, clearly $|B|<|F_1|/2$.  We may assume that \eqref{e.P1decay} holds for some $k_1 \ge 0$, provided that we have sufficient decay in the estimate.  Let $|f_1|\le 1_{F_1-B}$ and $|f_2|\le 1_{F_2}$ and $|f_3|\le 1_{F_3}$, we will show that
\begin{eqnarray}\label{e.M3}
\Lambda_{\P_1,\P_2}(f_1,f_2,f_3) &\lesssim& |F_2|^{1-\frac 1{p_2}} |F_3|^{\frac 1{p_2}}
\end{eqnarray}
which implies desired estimates by letting $p_2$ close to $2r/(4-r)$.

To show \eqref{e.M3}, we first use convexity and \eqref{e.trilinearCC} and  Carleson embeddings to obtain
\begin{eqnarray}
\nonumber \Lambda_{\P_1,\P_2}(f_1,f_2,f_3) &\lesssim&  |F_3|^{\frac 1{p_3}}  M_N(\P_1, 1_{F_3})^{\frac 1{p_1}} \|a_1\|_{\infty}^{1-\frac 2{p_1}} \|a_1\|_{2,\infty}^{\frac 2 {p_1}} M_N(\P_2,F_3)^{1/p_2}\|a_2\|_{p_2}\\
\nonumber &\lesssim& 2^{-Nk_1} |F_3|^{\frac 1{p_3}} (\frac{|F_3|}{|F_1|})^{\frac 1 {p_1}} |F_1|^{\frac 1 {p_1}} M_N(\P_2,F_3)^{\frac 1 {p_2}} \|a_2\|_{p_2} \\
\label{e.firstcutoff} &\lesssim& 2^{-Nk_1} M_N(\P_2,F_3)^{1/p_2} \|a_2\|_{p_2}  |F_3|^{\frac  1{p_2'}} \ \ .
\end{eqnarray}
It follows in particular that $\Lambda(f_1,f_2,f_3)\le 2^{-Nk_1} |F_2|^{1/p_2} |F_3|^{1/p_2'}$ for any $2<p_2<2r/(4-r)$, which would imply the desired estimate  \eqref{e.M3} if $|F_3|\le |F_2|$. When $|F_3|>|F_2|$ we will carry out essentially another layer of restricted weak-type interpolation, which we details below. Let  $\widetilde B_1 = F_3 \cap \{M 1_{F_2}>C_0|F_2|/|F_3|\}$, clearly $|\widetilde B_1|<|F_3|/2$. We will show that
\begin{eqnarray}\label{e.f2f3interpolation-1}
\Lambda_{\P_1,\P_2}(f_1,f_2,f_31_{\widetilde B_1^c}) &\lesssim& 2^{-Nk_1} |F_2|^{1-1/p_2} |F_3|^{1/p_3} \ \ .
\end{eqnarray}
To show \eqref{e.f2f3interpolation-1},  we may assume  that for some $k_2\ge 0$ it holds for every $P'\in \P_2$ that
$2^{k_2} \le 1+\frac{\dist(I_{P'}, \widetilde B_1^c)}{|I_{P'}|} < 2^{k_2+1}$, provided that we have enough decay in the estimates. It follows from \eqref{e.firstcutoff} that
\begin{eqnarray*} 
\Lambda_{\P_1,\P_2}(f_1,f_2,f_31_{\widetilde B_1^c}) &\lesssim& 2^{-Nk_1} M_N(\P_2, F_3 - \widetilde B_1)^{ \frac1 {p_2}} \|a_2\|_{\infty}^{1-\frac 2{p_2}} \|a_2\|_{2,\infty}^{\frac 2 {p_2}} |F_3|^{\frac 1{p_2'}}\\
&\lesssim& 2^{-Nk_1} 2^{-Nk_2} (\frac{|F_2|}{|F_3|})^{1-\frac 2 {p_2}} |F_2|^{\frac 1 {p_2}} |F_3|^{\frac 1 {p_2'}}\\
&\lesssim& 2^{-N(k_1+k_2)} |F_2|^{1-\frac 1{p_2}} |F_3|^{\frac 1 {p_2}}
\end{eqnarray*}
proving \eqref{e.f2f3interpolation-1}. Now, if $|\widetilde B_1|>|F_2|$ then we continue to let $\widetilde B_2 = \widetilde B_1\cap \{M1_{F_2}> C_0|F_2|/|\widetilde B_1|\}$ and similarly $\widetilde B_3, \dots, \widetilde B_m$ such that $|F_3|>2|\widetilde B_1|>\dots > 2^{m}|\widetilde B_m|$ until $|\widetilde B_m|\le |F_2|$. (This process has to stop since $|F_3|/|F_2|$ is finite.) We obtain (here for convenience of notation let $\widetilde B_0\equiv F_3$):
\begin{eqnarray*}\Lambda_{\P_1,\P_2}(f_1,f_2,f_3) &\lesssim&    |\Lambda_{\P_1,\P_2}(f_1,f_2,f_3 1_{\widetilde B_m})|  + \sum_{j=0}^{m-1} |\Lambda_{\P_1,\P_2}(f_1,f_2, f_31_{\widetilde B_j - \widetilde B_{j+1}})| \\
&\lesssim&  2^{-Nk_1} |F_2|^{\frac 1 {p_2}} |\widetilde B_m|^{1- \frac 1 {p_2}} + \sum_{j=0}^{m-1} 2^{-Nk_1}|F_2|^{1-\frac 1{p_2}} (2^{-j}|F_3|)^{\frac 1 {p_2}} \\
&\lesssim& 2^{-Nk_1}|F_2|^{1-\frac 1{p_2}} |F_3|^{\frac 1 {p_2}}  \ \ , \ \ \text{as desired.}
\end{eqnarray*}

\underline{Near $M_4$.} Let $B=\{M1_{F_3}>C_0 |F_3|/|F_1|\}$ which satisfies $|B| < |F_1|/3$. Let $|f_1|\le 1_{F_1-B}$, $|f_2|\le 1_{F_2}$, $|f_3|\le 1_{F_3}$, it suffices to show that
\begin{eqnarray}\label{e.M4}
\Lambda_{\P_1,\P_2}(f_1,f_2,f_3) &\lesssim_\epsilon& |F_2|^{\epsilon} |F_3|^{1-\epsilon} \ \ ,
\end{eqnarray}
 for any $0<\epsilon \le 1/p_2$. As before, we will assume \eqref{e.P1decay}. For any $2<p_2<2r/(4-r)$ we also have \eqref{e.firstcutoff}, which would imply the desired estimate if $|F_3|\ge |F_2|$. When $|F_3| < |F_2|$ we will use interpolation. Similar to the analysis near $M_3$, it suffices to show that if $\widetilde B_1 = F_2 \cap \{M1_{F_3}>C_0 |F_3|/|F_2|\}$ (which satisfies $|\widetilde B_1|<|F_2|/2$) then
\begin{eqnarray*}
\Lambda_{\P_1,\P_2}(f_1, f_21_{\widetilde B_1^c}, f_3) &\lesssim& 2^{-Nk_1}|F_2|^{\epsilon} |F_3|^{1-\epsilon}
\end{eqnarray*}
and we will again assume that for some $k_2\ge 0$ it holds for every $P'\in \P_2$ that $2^{k_2}\le 1+ \frac{\dist(I_{P'},\widetilde B_1^c)}{|I_{P'}|} < 2^{k_2}$. Now it follows from \eqref{e.firstcutoff} that
\begin{eqnarray*}
\Lambda_{\P_1,\P_2}(f_1,f_21_{\widetilde B_1^c},f_3) &\lesssim& 2^{-Nk_1} M_N(\P_2, F_3)^{\frac 1{p_2}}  \|a_2\|_{\infty}^{1-\frac 2{p_2}} \|a_2\|_{2,\infty}^{\frac 2 {p_2}} |F_3|^{\frac 1{p_2'}}\\
&\lesssim&  2^{-Nk_1} 2^{-Nk_2} \min(1,\frac{|F_3|}{|F_2|})^{\frac 1 {p_2}}  |F_2|^{\frac 1 {p_2}} |F_3|^{\frac 1 {p_2'}}\\
&\lesssim& 2^{-N(k_1+k_2)} |F_2|^{\epsilon} |F_3|^{1-\epsilon} \ \ .
\end{eqnarray*}
This completes the proof of Theorem~\ref{t.CC-estimate}. \endproof

\section{Estimates for $LM$ model operators}\label{s.LM-estimates}
\begin{theorem}\label{t.LM-estimates} Let $T(f_1,f_2)$ be an LM model operator. Then $\|T\|_{L^{q_1}\times L^{q_2} \to L^{q_3}} < \infty$
for all  $1/q_3 = 1/{q_1} + 1/{q_2}$ with $q_1, q_2>2r/(3r-4)$ and $q_3>r'/2$. 
\end{theorem}
 
\proof We sketch the  proof of this theorem, which is a simple bilinear extension of the proof of Theorem~\ref{t.varCarl}.  For a tritile $Q\in \Q$ in the definition of $T$ we let $\omega_{Q}$ be the convex hull of $\omega_{Q_1}$,$\omega_{Q_2}$, $\omega_{Q_3}$. Let $D^{-1}x=x/2$ the dilation by $1/2$ with respect to the origin, and define $\omega_{Q,lower}=\omega_{1,Q,lower}\cap \omega_{2,Q,lower}\cap D^{-1}(\omega_{3,Q,lower})$ and define $\omega_{Q,upper}$ similarly.  Without loss of generality we may assume that $\omega_{Q,lower}$ is a half-line and $\omega_{Q,upper}$ is a finite interval for every $Q\in \Q$. We now could define $\widetilde \omega_Q$ to be the convex hull of $20\omega_{Q}$ and $20\omega_{Q,upper}$, and from here we may define generating subsets of $\Q$, and construct outer measure spaces on $\Q$ using the usual outer measure  and sizes as in Section~\ref{s.outerP}.

Let $A$ denote the collection of intervals $J$ such that for some $P_1,P_2\in \Q$ we have $I_{P_1}\subset J \subset 30 I_{P_2}$.   Let $p_1,p_2>2$ and $p_3>r/(r-2)$ such that $\sum 1/p_j=1$; note that this implies $p_1,p_2<2r/(4-r)$.  By routine applications of outer measure techniques and embedding theorems in Section~\ref{s.carl-embed}, we obtain 
\begin{eqnarray*}
\<T(f_1,f_2),f_3\>   &\lesssim_N&  |F_3|^{\frac 1{p_3}} \sup_{J\in A} M_J(1_{\supp(f_3)} \widetilde \chi_J^N)^{1-\frac 1{p_3}} \prod_{j=1,2} |F_j|^{\frac 1{p_j}} \sup_{J\in A}  M_{J} (f_j\widetilde \chi_{J}^N)^{1- \frac 2{p_j}} \ \ .
\end{eqnarray*}
Using similar arguments as in previous sections, it follows that restricted weak-type estimates holds for at least one $(\alpha_1,\alpha_2,\alpha_3)$ in any   neighborhood of any of the following points, which then implies the theorem (below $\beta_r:=(3r-4)/(2r)$):
$$H_1(\beta_r, -\beta_r, 1) \ \ ,  \ \ H_2(-\beta_r,\beta_r, 1) \ \ , \ \ H_3(1/2,  \beta_r,  1/2-\beta_r)   \ \  , \ \ H_4(\beta_r, 1/ 2,  1/2-\beta_r)    \ \ . \ \ \text{\endproof}$$

\appendix

\section{Reduction  to discrete operators}
We will show below that the (variation-norm) bilinear Fourier operators with symbols $m_{CC}$, $m_{C\times C}$, $m_{BC}$, $m_{CC}$ are controlled by the respective discrete operators.

Recall that $A$ is the set of all admissible triple $(side,m,n)$ in $\mathcal H$. It is clear that  
\begin{itemize}
\item If $\alpha \in A$ then $L_1\le \min(m_\alpha, n_\alpha)\le 2L_1$.
\item If $k>2L_1$ then $(left, 2L_1, k)$ and $(right, k, 2L_1)$ are not in $A$.
\end{itemize}
 
Let $\D_{left}$, $\D_{right}$ be the sets of  left and and right-sided dyadic intervals. By definition, if $\alpha=(side,m,n)$ then $\I_{M,N}(\alpha)$ consists of $I\in \D_{side}$ with $M\in  I_{lower,\alpha}:= I - (m+1)|I|$ and $N \in I_{upper,\alpha} := I+(n+1)|I|$.

\subsection{Discretization of $m_{CC}$}

We will discuss the discretization for $m_{CC,1}$, the discretization for other $m_{CC,j}$ are  similar. We first make several observations regarding the decomposition of $1_{M<\xi<N}$ in Section~\ref{s.interval-decomp}.

\begin{lemma}[Observation 1]\label{l.ob1} If $I\in \I(side, m,n)$ and  $\xi\in \frac 54 I$ then
\begin{equation}\label{e.relative1}
m/(4n)  <  |\xi-M|/|\xi -N|  <4m/n \ \ .
\end{equation}
\end{lemma}
\proof This follows from $(m- 1/8)|I| \le |\xi-M|, |\xi-N|  \le (m+2+ 1/8)|I|$. \endproof

\begin{lemma}[Observation 2]\label{l.ob2} Let $I\in  \I(side,m,n)$ with $n\ge 4m$. Assume that  $J\in \I(side',m',n')$ and $\sup J< \inf I$. Then $n'/m' > n/m$.

Similarly, if $m\ge 4n$ and $J$ is on the right of $I$ then $m'/n'>m/n$.
\end{lemma}
\proof Assume $n\ge 4m$, the other case is symmetric. Without loss of generality we may assume that $J$ is adjacent to $I$. Since $n>m$, it follows that $|J|=|I|$ or $|J|=|I|/2$. In the first case the desired estimate is trivial, and in the second case we have $m'\le 2m$ and $n' \ge 2n+2$, which also implies the desired estimate. \endproof

As a corollary of Lemma~\ref{l.ob2}, it follows that the intervals in $A_2$ are strictly in between the elements of $A_1$ and the elements of $A_3$.

The next observation concerns cancellation when summing bump functions $\phi_\alpha$ over the set of $\alpha\in A$ such that $I\in \I_{M,N}(\alpha)$, where $M,N,I$ are fixed.

\begin{lemma}[Observation 3] \label{l.ob3}  (i) Given any $C \ge 4$ we can find 
$O(L_1)$ many tuples $(a_j, b_j,  u_j, v_j)\in \Z^2\times \{1/2 , 1, 2\}^2$, with $L_1 +1 \le a_j, b_j \lesssim L_1$, such that the following holds for every left dyadic interval $I$ and every $M<N$:
$$\sum_{\alpha \in A: \ \ m_\alpha\ge Cn_\alpha}  \phi_{\alpha}(\xi) 1_{I\in \I_{M,N}(\alpha)} = \sum_j \phi_{u_j,v_j}(\xi) 1_{M < c(I)-(a_j-1/2)|I|} 1_{N\in I+ b_j|I|} \ \ .$$
Furthermore, an analogous statement also holds for right dyadic intervals.

(ii) A similar statement also holds for $\sum_{m_\alpha\le n_\alpha/C}$.

\end{lemma}
\underline{Remark:} the key idea here is that the left hand sides in the above equalities involve infinitely many terms, while the right hand sides  contain only $O(L_1)$  terms.

\proof   Let $I_-$ and $I_+$ be the left and right neighbors of $I$ in $\mathcal H$.  Below we only consider the  $\sum_{m_\alpha \ge C n_\alpha}$, the other sum could be handled similarly.

Since $C\ge 4$ we have $m_\alpha \ge    4L_1$, while $2L_1\ge n_\alpha\ge L_1$. Our key observation is the fact that: in the sum, the ratios $u(I)=|I_-|/|I|$ and $v(I)=|I_+|/|I|$ depends only on $side_\alpha$ and $n_\alpha$. In other words, knowing the $side$ of $I$ and knowing $n_\alpha$ (only a finite number of possible values) we could determine $\phi_\alpha (\xi)\equiv \phi_{u(I),v(I)}(\xi)$ completely and thus we have the freedom to sum the indicator constraints on $M$, $1_{M\in I-(m+1)|I}$, over $m\ge Cn_\alpha$. The following table details this observation

\begin{table}[h]
\begin{tabular}{|c|c||c|}
\hline
side of $I$  & values of $n_\alpha$  & Corresponding values of $(u(I), v(I))$\\ \hline
left & $L_1$ & $(2,\frac 1 2)$\\ \hline
left & $L_1<n_\alpha \le 2L_1$ & $(2,1)$\\ \hline
right & $L_1$ & $(1,\frac 1 2)$\\ \hline
right & $L_1<n< 2L_1$ & $(2,1)$\\ \hline
\end{tabular}
\end{table}
\noindent (note that if $I$ is right-sided then $n_\alpha<2L_1$ by definition of $A$).  \endproof

\subsubsection{Discretization for $m_{CC,1}$:} Using Lemma~\ref{l.ob3}, we may decompose $m_{CC,1}$ into
\begin{eqnarray*}
\sum_{|I| \le |J|/16} \phi_I(\xi_1) \varphi_J(\xi_2) 1_{M < c(I)-a_1|I|, N\in I+b_1|I| } 1_{M < c(J)-a_2|J|, N\in J+b_2|J| } 
\end{eqnarray*}
in the sum $I$ and $J$ belong to fixed collection of dyadic intervals, $\phi_I$, $\varphi_J$ are  given bump functions supported in $(5/4)I$ and $(5/4)J$,  and $L_1 \le a_1,b_1,a_2,b_2 \lesssim L_1$.

It follows from a standard Fourier sampling argument that we can decompose
$$\int \int e^{ix(\xi_1+\xi_2)} \phi_I(\xi_1) \varphi_J(\xi_2)\widehat f_1(\xi_1) \widehat f_2(\xi_2) d\xi_1d\xi_2$$
into finitely many wavelet sums
$$\sum_{P,P' \  tiles: \, \, \omega_P=I, \  \omega_{P'} = J}  |I_P| |I_{P'}|\<f_1,\phi_{1,P}\> \<f_2,\phi_{2,P'}\> \phi_{1,P}(x)\phi_{2,P'}(x) \ \  ,$$
where $\phi_{1,P}$, $\phi_{2,P'}$ are $L^1$-normalized wave functions adapted to the tiles $P$, $P'$. Thus the resulting bilinear operator for $m_{CC,1}$ is controlled by a finite sum over $CC$ discrete operators.

\subsection{Discretization of $m_{BC}$}

For convenience let $\xi_3=-\xi_1-\xi_2$ below. Thanks to fast decay of $a_k$, it suffices to consider the contribution of one fixed $k$.  
In other words, it suffices to  consider symbols
\begin{eqnarray*}
\widetilde m_{BC} 
&=& \sum_{\alpha\in A} \sum_{\ell(S)\le |I|} \phi_{1,S}(\xi_1) \phi_{2,S}(\xi_2) \phi_{3,S}(-\xi_3)\phi_{\alpha,I}(-\xi_3) 1_{\{2M\in I_{lower,\alpha}\}} \,  1_{\{2N \in I_{upper,\alpha}\}}
\end{eqnarray*}
 in the sum  $S$ is in a fixed collection of shifted dyadic cubes with the property \eqref{e.whitney}, and $I$ is in a fixed collection of dyadic intervals, and $\phi_{j,S}$ are uniformly  $C^n$ bump functions  supported in $\frac 56 S_i$, where $n$ could be chosen arbitrarily large.

For any $I$, $S$, $\alpha\in A$, and any Schwarz $f_1$, $f_2$, $f_3$, we have
\begin{eqnarray*}
&&\int\int \int e^{ix(\xi_1+\xi_2)} \phi_{3,S}(\xi_1+ \xi_2)  \phi_{\alpha,I}(\xi_1+\xi_2) \big(\prod_{j=1}^2 \phi_{j,S}(\xi_j) \widehat f_j(\xi_j)  d\xi_1 d\xi_2 \big)f_3(x)dx\\
&=& \int\int\int_{\xi_1+\xi_2+\xi_3=0} \phi_{\alpha,I}(-\xi_3) \phi_{3,S}(-\xi_3)\widehat f_3(\xi_3) \Big(\prod_{j=1}^2 \widehat {f_j}(\xi_j) \phi_{j,S}(\xi_j)\Big) d\xi_1d\xi_2d\xi_3  \\
&=& C\int  \Big(f_3*\widehat{\phi_{3,S}}*\widehat{\phi_{\alpha,I}}\Big)(2y) \prod_{j=1}^2 \Big(f_j*\widecheck{\phi_{j,S}}\Big)(y) dy\\
&=&C\sum_{m\in\Z} \int_0^1 \frac{1}{\ell(S)}  \Big(f_3*\widehat{\phi_{3,S}}*\widehat{\phi_{\alpha,I}}\Big)(\frac{2(t-m)}{\ell(S)}) \prod_{j=1}^2 \Big(f_j*\widecheck{\phi_{j,S}}\Big)(\frac{t-m}{\ell(S)}) dt\\
&=&C\int_0^1 \sum_{Q\in \Q_{t,S}}  |I_Q|  \<f_3*\widehat{\phi_{\alpha,I}}, \overline{\phi_{3,Q}}\>\prod_{j=1}^2 \<f_j,\phi_{j,Q}\> dt \ \ , \ \ \text{where}
\end{eqnarray*}
\begin{itemize}
\item for each $t\in [0,1)$, $\Q_{t,S}$ is the collection of all tritiles $Q=(Q_1,Q_2,Q_3)$ with 
$$Q_1=I_Q\times S_1 \ \ , \ \ Q_2 = I_Q \times S_2 \ \ , \ \ Q_3 = I_Q \times  S_3  \ \ ,$$ 
and $I_Q$ is a shifted dyadic interval of length $\ell(S)^{-1}$ (with $t$-dependent shift);
\item for any $Q \in \Q_{t,S}$ and for each $j=1,2,3$, define $\phi_{j,Q}(x) :=   \widecheck{\phi_{j,S}}(x-c(I_Q))$
which  are $L^1$-normalized wave packets adapted to $Q_j=I_Q\times S_j$ (with frequency support in $\frac 56 S_j$).
\end{itemize}

Recall that for every $S\in \S$ the distance between  $S_1$ and $S_2$ is comparable to $L_2 |\ell(S)| \sim L_1|\ell(S)|$ (due to the Whitney condition \eqref{e.whitney}). Clearly, the collection $\Q_t$ will be of rank $1$ (with uniform constants over $0\le t \le 1$) if $L_1$ is sufficiently large.  

Now, for any $b_3$-shifted dyadic interval $I$, by a Fourier sampling argument $(f_3*\widehat{\psi_{\alpha, I}})(x)$ equals a sum over finitely many terms of the form $\sum_{P \in \P_I}  \<f_3,\overline{\phi_P}\>\overline{\phi_P(x)}$, 
where $\P_I$ is the collection of all rectangles $P=J\times I$ formed using standard dyadic intervals $J$ of length $|I|^{-1}$, and $\{\phi_P, \, P \in \P_I\}$ is a collection of Fourier wave packets adapted to $P\in \P_I$ with $supp(\widehat \phi_{J\times I}) \subset \frac 54 I$).

It follows that, modulo a multiple by some absolute constant,
$$\int \big(\int \int e^{ix(\xi_1+\xi_2)} \widetilde m_{BC} \widehat f_1(\xi_1)\widehat f_2(\xi_2) d\xi_1d\xi_2\big) \ \  f_3(x)dx$$
can be decomposed into finitely many terms of the following form:
\begin{eqnarray*}
&&\int_0^1 \sum_{Q\in \Q_t} |I_Q| \<f_1,\phi_{1,Q}\> \<f_2,\phi_{2,Q}\> \<\mathcal{C}(f_3),\overline{\phi_{3,Q}}\>dt \ \ , \ \ \text{where}\\
&&\mathcal{C}(f_3)(x) :=\sum_{\alpha, I}\sum_{P\in \P_I} 1_{\{2M\in I_{lower,\alpha}\}} 1_{\{2N\in I_{upper,\alpha} \}} \<f_3,\overline{\phi_P}\>\overline{\phi_P(x)} \ \ ,
\end{eqnarray*}
and $\Q_t := \bigcup_{S\in \S}\Q_{t,S}$. Let $\P$ be the union of $\P_I$ over $I\in \D$.  Using  Lemma~\ref{l.ob3}, we can write $\mathcal{C}(f_3)(x)$ as a sum over finitely many  terms of the form
$$\pm \sum_{P\in \P}|I_P|\<f_3,\overline{\phi_P}\> \overline{\phi_P(x)} 1_{\{2M\in \omega_{P,lower}\}} 1_{\{2N \in \omega_{P,upper}\}}$$
where $\{(\omega_{P,lower}, \omega_P, \omega_{P,upper}), P \in \P\}$ is rigid. Thus to bound the resulting bilinear operator for $m_{BC}$,  it suffices to estimate  discrete model $BC$ operators.

\subsection{Discretization of $m_{LM}$}

In Lemma~\ref{l.CC-properties} and Lemma~\ref{l.BC-properties} we will show that the supports of $m_{CC}$ and $m_{BC}$ are inside $\{M\le \xi_1\le \xi_2\le N\}$ and they do not intersect except at possibly  $(\xi_1,\xi_2)=(M,M)$ and $(N,N)$. It will follow that
\begin{equation}\label{e.LM-factors}
m_{LM}=(\chi_{M<\xi_1<\xi_2<N}-m_{CC})(\chi_{M<\xi_1<\xi_2<N}-m_{BC}) \ \ ,
\end{equation}
which will be used in subsection~\ref{s.LM-properties} to reduce the bilinear multiplier operator with symbol $m_{LM}$ to LM model operators.

\subsubsection{Support of $m_{CC}$}\label{s.CC-properties}
It will follow from Lemma~\ref{l.CC-properties} below that $m_{CC} = 1$ on $R_1$ (defined in \eqref{e.R1-defn}), and   $m_{CC}$ is supported inside the following enlargement of $R_1$:
$$R'_1:= \Big\{(\xi_1,\xi_2) \in [M,N]^2: \min(|\xi_1-M|, |\xi_2-N|) \le \frac 1{3} |\xi_1-\xi_2| \Big\} \ \ .$$
Note that $\overline{R_1}\subset R'_1$, and $R'_1 \subset \{M \le \xi_1<\xi_2 \le N\}$.

\begin{lemma}\label{l.CC-properties}
(i) If a summand in \eqref{e.product}  is non-zero for some $(\xi_1,\xi_2) \in R_1$, 
then this summand must appear in one of $m_{CC,k}$, $1\le k \le 5$.

(ii) All summands of $m_{CC,k}$ are supported in $R'_1$.
\end{lemma}

\proof

(i) Suppose that for some $\alpha,\beta$ and $I\in \I(\alpha)$ and $J\in \I(\beta)$ and $(\xi_1,\xi_2)\in R_1$ we have $\phi_{\alpha,I}(\xi_1) \phi_{\beta,J}(\xi_2)   \ne 0$. 

Without loss of generality we may assume that $(\alpha,\beta) \not\in A_1\times A_3$, so $\alpha \not\in A_1$ or $\beta\not\in A_3$. By symmetry, we may assume  that $\beta \not\in A_3$. 

We first show that $\alpha \in A_1$. Since $\beta\not\in A_3$, it follows from Lemma~\ref{l.ob1} that $|\xi_2-M| \quad \le \quad   16 |\xi_2-N|$,  therefore
$$|\xi_2-N| \quad \ge \quad  |M-N|/15 \quad \ge \quad |\xi_1-\xi_2|/15\ \ .$$
 Since $(\xi_1,\xi_2)\in R_1$, it follows from the   definition of $R_1$ that 
\begin{equation}\label{e.xi1-M-small}
|\xi_1-M| \quad \le \quad  |\xi_1-\xi_2|/200  \quad \le \quad |M-N|/200 \ \ .
\end{equation}
It follows that $|\xi_1-M| \le |\xi_1-N|/199$. Using Lemma~\ref{l.ob1} again, it follows that 
$$m_\alpha/n_\alpha \quad \le \quad  4 |\xi_1-M|/|\xi_1-N| \quad < \quad 1/4$$
thus  $\alpha \in A_1$ as claimed above. 

It remains to show that $|I| < |J|/16$. As a corollary of the condition $\alpha \in A_1$ and $\beta\not\in A_3$,  we have $m_\beta \le 4 n_\beta \le 8L_1$, while clearly $m_\alpha\ge L_1$.  Using $\xi_1\in \frac 54 I$ and $\xi_j \in \frac 54 J$ and  \eqref{e.xi1-M-small} it follows that
\begin{eqnarray*}
\frac{|I|}{|J|} &\le& \frac{|\xi_1-M|/(m_\alpha-1/8)}{|\xi_2-M|/(m_\beta+9/8)} \quad  \le \quad \frac{11|\xi_1-M|}{|\xi_2-M|} \quad < \quad \frac 1 {16} \ \ .
\end{eqnarray*}
This completes the proof of claim (i).

(ii) Without loss of generality we may assume that $\alpha \in A_1$.

By the definition of $m_{CC,j}$'s we have $|I|\le |J|/16$.  It suffices to show that
$$(5/4) I \times (5/4) J \quad \subset  \quad R'_1 \ \ .$$
To see this, take any $(\xi_1,\xi_2) \in \frac 54 I \times \frac 54 J$.

Since $\alpha\in A_1$ implies that $m_\alpha = \min(m_\alpha,n_\alpha) \le  2L_1$, while clearly $m_\beta \ge L_1 \ge 1$.  It follows that
$$\frac{|\xi_1-M|}{|\xi_2-M|}   \quad \le \quad  \frac{|I|(m_\alpha+9/8)}{|J|(m_\beta-1/8)} \quad \le \quad \frac 1 4  \ \ .$$
Thus $|\xi_1-M| \le |\xi_1-\xi_2|/3$, as desired.
\endproof

\subsubsection{Support of $m_{BC}$}\label{s.BC-properties}
Consider the following enlargement of $R_2$ (defined in \eqref{e.R2-defn}):
$$R'_2:= \Big\{(\xi_1,\xi_2): |\xi_1-\xi_2|\le \frac 1 2 \min\big(|\frac 1 2(\xi_1+\xi_2)-M|, |\frac 1 2(\xi_1+\xi_2)-N|\big)\Big\}$$ 

\begin{lemma}\label{l.BC-properties} Let $m_{BC}$ be defined by Definition~\ref{d.BCsym}. Then:

(i) $m_{BC}$ is supported inside  $R'_2$.

(ii) Any summand of  \eqref{e.expansionBC} whose support intersects $R_2$ must appear in $m_{BC}$.

\end{lemma}
\proof
(i) Take any $(\xi_1,\xi_2)$ in the support of $m_{BC}$, then for some  $S\in \S$ and $I\in \I_{2M, 2N}(\alpha), \alpha\in A,$ such that $\ell(S)\le |I|$ we have
$$\phi_{1,S,k}(\xi_1) \phi_{2,S,k}(\xi_2) \phi_{3,S,k}(\xi_1 + \xi_2)
\phi_{\alpha,I}(\xi_1 + \xi_2)   \ne 0 \ \ .$$ 
It follows that $(\xi_1,\xi_2)\in \frac 56 S_1 \times \frac 56 S_2$, therefore  using \eqref{e.whitney} we obtain
\begin{eqnarray*}
|\xi_1-\xi_2| &=& \sqrt 2  \dist((\xi_1,\xi_2),\{\xi_1=\xi_1\})  \quad < \quad 5L_2  \ell(S) \ \ .
\end{eqnarray*}
On the other hand, since $\xi_1+\xi_2 \in \frac 54I$ and $m_\alpha,n_\alpha \ge L_1$, we have
\begin{eqnarray*}
\min\big(|\frac 1 2(\xi_1+\xi_2)-M|, |\frac 1 2(\xi_1+\xi_2)-N|\big)  
&\ge&    \frac 1 2 (L_1-1/8)|I| \ \ .
\end{eqnarray*}
Since $\ell(S)\le |I|$ and since $L_2=L_1/40$, it is not hard to check that $(\xi_1,\xi_2)$ satisfies the defining property of  $R'_2$.

(ii) Suppose that $(\xi_1,\xi_2)\in R_2$ such that
$$\phi_{1,S,k} (\xi_1)\phi_{2,S,k}(\xi_2) \phi_{3,Q,k} (\xi_1+\xi_2)
\phi_{\alpha,I}(\xi_1 + \xi_2)    \ne 0 \ \ .$$ 
We will show that $\ell(S) \le |I|$. 

First, using $(\xi_1,\xi_2)\in R_2$ and using $\xi_1+\xi_2 \in  (5/4) I$, it follows that
\begin{eqnarray*}
|\xi_1-\xi_2| &\le&
\frac 1 {100}\min\big(|\frac 1 2(\xi_1+\xi_2)-M|, |\frac 1 2(\xi_1+\xi_2)-N|\big) \\
&\le& \frac 1{100} (2L_1+1/8)|I| \ \ .
\end{eqnarray*}
Since $(\xi_1,\xi_2) \in (4/5)S_1 \times (4/5)S_2$, it follows that
\begin{eqnarray*}
|\xi_1-\xi_2| &=& \sqrt 2 \dist((\xi_1,\xi_2), \{\xi_1=\xi_2\}) \quad \ge \quad \sqrt{2} L_2 \ell(S) \ \ .
\end{eqnarray*}

Collectin estimates and using $L_2=L_1/40$, it is clear that $\ell(S)  \le |I|$, thus the corresponding summand is part of $m_{BC}$.
\endproof

\subsubsection{Reduction to discrete $T_{LM}$, part I: decomposition into simpler trilinear symbols}\label{s.LM-properties}

Recall that $D_t$ is the dilation $D_t A:= \{2^t x: x\in A\}$. 

\begin{lemma}\label{l.LM-properties}  $m_{LM}$ can be decomposed into finitely many symbols of the form
\begin{eqnarray*}
\sum_{\alpha,\beta,\gamma} &&\sum_{I\in \D_1} \sum_{J \in \D_2} \sum_{L \in \D_3} 1_{\{constraints \, on \, I, J, L\}}  \ \ \phi_{1,I}(\xi_1)\phi_{2,J}(\xi_2)\phi_{3,L}(-\xi_1-\xi_2) \times \\
&\times& 1_{\{M\in I_{lower,\alpha}, N \in I_{upper,\alpha}\}}  \times 1_{\{M\in J_{lower,\alpha}, N \in J_{upper,\alpha}\}} \times 1_{\{2M\in L_{lower,\alpha}, 2N \in L_{upper,\alpha}\}} \ \ , 
\end{eqnarray*}
\begin{itemize}
\item $\D_1$, $\D_2$, $\D_3$ are three fixed collections of dyadic intervals;
\item $\phi_{1,I},\phi_{2,J},\phi_{3,L}$ are $C^n$-bump functions  adapted to these intervals, with support in $\frac 54 I$, $\frac 54 J$, $\frac 54 L$; and $n$ could be chosen arbitrarily large
\item the constraints on $I,J,L$ read as follows: the length  $I,J,L$ are comparable, and the distance between $I,J,D_{-1}L$ are $O(|I|)$.
\end{itemize}
\end{lemma}

Note that the dependence on $M$ and $N$ are only in the last three factors of the summands.  To prove Lemma~\ref{l.LM-properties}, we analyze the  factors in the factorization \eqref{e.LM-factors}.

\noindent \underline{Part I: First factor.} Using Lemma~\ref{l.CC-properties}, it follows that
\begin{eqnarray}
\nonumber \chi_{M<\xi_1<\xi_2<N}-m_{CC} &=& \chi_{\xi_1<\xi_2}(\chi_{M<\xi_1,\xi_2<N}-m_{CC}) \\
\label{e.lmcc}
&=&\chi_{\xi_1<\xi_2}\sum_{\alpha,\beta \in A} \ \ \sum_{(I,J)\in \I(\alpha) \times \I(\beta) \text{, constraints}} \phi_{\alpha,I}(\xi_1)\phi_{\beta,J}(\xi_2) \ \ ,
\end{eqnarray}
the constraints on $(I,J) \in \I(\alpha) \times \I(\beta)$ read as follows: 
\begin{itemize}
\item If $\alpha\not\in A_1$ and $\beta\not\in A_3$ then there are no constraints.
\item If $\alpha \in A_1$ and $\beta\in A_3$ then no $I$,$J$ are allowed.
\item If $\alpha \in A_1$ and $\beta\not\in A_3$ then one requires $|I| > |J|/16$.
\item If $\alpha \not\in A_1$ and $\beta \in A_3$ then one requires $|J|> |I|/16$.
\end{itemize}

Recall that $\mathcal H$ is the union of $\I(\alpha)$ over $\alpha \in A$. We will show that

\begin{lemma}\label{l.lengthIJ} Assume that $(I,J)\in \I(\alpha) \times \I(\beta)$ contributes to the right hand side of \eqref{e.lmcc}. Then (i)   $|I| \sim |J|$ and (ii) if $\xi_1<\xi_2$ and $\phi_{\alpha,I}(\xi_1)\phi_{\beta,J}(\xi_2)\ne 0$   then
\begin{eqnarray*}
\dist((\xi_1,\xi_2), \{(M,M),(N,N)\}) &\sim& L_1 |I| \ \ .
\end{eqnarray*}
\end{lemma}

\proof (i) Without loss of generality we may assume that the ratio of the lengths of $I$ and $J$ is not in $\{1,\frac 12, 2\}$. As the lengths of adjacent intervals in $\I$ differ by a ratio of $1$ or $2$ or $1/2$, it follows that $I\ne J$ and they are not adjacent. 

Now, we show that $\sup I\le \inf J$. Assume towards a contradiction that $\sup J \le \inf I$. Since there is at least one interval in $\mathcal H$ between $I$ and $J$ and since the lengths of adjacents intervals in $\mathcal H$ differ by a factor in $\{1/2,1,2\}$, it follows that  $\inf \frac 54 I > \sup \frac 54 J$. Consequently, $(5/4)I \times (5/4) J\cap \{(\xi_1,\xi_2): \xi_1<\xi_2\} = \emptyset$,
contradicting the fact that $(I,J)$ contributes to \eqref{e.lmcc}. 

Now, if $\alpha,\beta\in A_2$ then clearly $|I|$ and $|J|$ are comparable to $|M-N|/L_1$, so they are comparable.  Therefore we may assume that either $\alpha\not\in A_2$ or $\beta \not\in A_2$.  

We will show that 
\begin{equation*} \text{either} \quad \alpha \in A_1 \quad  \text{or} \quad \beta\in A_3 \ \ . 
\end{equation*}Note that from the constraints we must have 
$(\alpha,\beta) \not\in A_1\times A_3$ (otherwise there won't be any $(I,J)$), thus the above two properties can not hold simultaneously. Assume towards a contradiction that  $\alpha \not\in A_1$ and $\beta \not\in A_3$. Since $\alpha \not\in A_2$ or $\beta\not\in A_2$, it follows that $\alpha \in A_3$ or $\beta \in A_1$. 
\begin{itemize}
\item If $\alpha \in A_3$ then using Lemma~\ref{l.ob2} and  $\sup I \le  \inf J$ it follows that 
$$m_\beta/n_\beta > m_\alpha/n_\alpha \ge 4 \ \ , $$
thus  $\beta \in A_3$, contradict to the above assumption. 
\item If $\beta\in A_1$ then $n_\alpha/m_\alpha \ge n_\beta/m_\beta \ge 4$ and thus $\alpha \in A_1$, contradiction again.
\end{itemize}

Below, without loss of generality assume that $\alpha\in A_1$. Thus $\beta \in A_1$ or $A_2$. 

 If $\beta \in A_1$ then since $J$ is on the right of $I$ we easily have $|I|\le |J|$, while $|J|<|I|16$ by the above constraints. Thus  $|I| \sim |J|$. 

If $\beta \in A_2$, then since $\bigcup_{\gamma\in A_2} \I(\gamma)$ has $O(1)$ elements of comparable lengths,  $|J|$ is comparable to $|I_0|$ the length of the  right most interval inside $\bigcup_{\gamma\in A_1}\I(\gamma)$. Note that either $I=I_0$ or $I$ is on the left of $I_0$, thus $|I| \le |J_0|$. Combining with the constraint $|J|<16|I|$, we obtain  $|J| \sim |I|$.

(ii) Let $\alpha,\beta\in A$ such that $I\in \I(\alpha)$ and $J \in \I(\beta)$. It is clear that
\begin{eqnarray*} &&\dist((\xi_1,\xi_2), \{(M,M),(N,N)\})\\
&\sim& \min(|\xi_1-M|+|\xi_2-M|, |\xi_1-N| + |\xi_2-N|) \ \ .
\end{eqnarray*}
If $(\alpha,\beta) \in (A_1\cup A_2)^2$ then using $|I|\sim |J$ we obtain
\begin{eqnarray*}&&\min(|\xi_1-M|+|\xi_2-M|, |\xi_1-N| + |\xi_2-N|)\\
 &\sim& |\xi_1-M|+|\xi_2-M|  \quad \sim \quad L_1 |I| \ \ .
\end{eqnarray*}
Similarly, if $(\alpha,\beta) \in (A_2\cup A_3)^2$  then the desired claim follows. Since $(\alpha,\beta) \not \in A_1\times A_3$ (by the constraints), the remaining case is  $(\alpha,\beta) \in A_3 \times A_1$, however this can't happen either since $I$ is on the left of $J$. \endproof

\noindent \underline{Part II: The second factor.}
For convenience of notation, let $\xi_3:=-\xi_1-\xi_2$. Thanks to Lemma~\ref{l.BC-properties}, we may write
\begin{eqnarray}
\nonumber &&\chi_{M<\xi_1<\xi_2<N} - m_{BC}  \quad = \quad \chi_{M<\xi_1,\xi_2<N}(\chi_{\xi_1<\xi_2}\chi_{2M< \xi_1+\xi_2 < 2N} - m_{BC} ) \\
\label{e.lmbc}
&=& \chi_{M<\xi_1,\xi_2<N}\sum_{k,\alpha} a_k\sum_{S \in \S} \sum_{\stackrel{\ell(S)\ge 2|L|}{L \in \I_{-2N,-2M}(\alpha)}} \phi_{\alpha,L}(-\xi_3)  \phi_{3,S,k}(-\xi_3)  \prod_{j=1}^3 \phi_{j,S,k}(\xi_j)  \ \ .
\end{eqnarray}

\begin{lemma}\label{l.lengthSL} Suppose that $(S,L) \in (\S, I_{2M,2N}(\alpha))$ and  contributes to the right hand side of \eqref{e.lmbc}. Then (i) $\ell(S) \sim |L|$ and (ii) if $(\xi_1,\xi_2) \in (M,N)^2$ is in the support of the corresponding summand then
$$\dist((\xi_1,\xi_2),\{(M,M),(N,N)\}) \sim L_1 \ell(S) \ \ .$$
\end{lemma}

\proof (i) Note that   $|L|\le \ell(S)/2$ by definition, thus it remains to show  $\ell(S) =O(|L|)$.  Without loss of generality, assume that $|L-2N|<|L-2M|$. Then
$$\dist(L, 2N) \sim L_1|L| \  \ .$$ 
Since $\phi_{j,S,k}$ is supported inside $(5/6)S_j$, it folows that
$(\frac 56 S_1 \times \frac 56 S_2)\cap (M,N)^2  \ne \emptyset$.
Take any $(\xi_1,\xi_2)$ in this intersection. Then $\xi_1<\xi_2<N$, therefore using $(\xi_1+\xi_2) \in supp(\phi_{\alpha,L}) \subset \frac 54 L$ it follows that
\begin{eqnarray*}
|\xi_1-\xi_2|  &\le& |\xi_1-N|+|\xi_2-N| \quad = \quad |\xi_1+\xi_2-2N|\\
&\le&    3|L|  +2\dist(L,2N) \quad \lesssim \quad L_1|L| \ \ .
\end{eqnarray*}
Now, let $\ell$ be the line $\{\xi_1=\xi_2\}$. Since $(\xi_1,\xi_2)\in S_1\times S_2$,  it follows from \eqref{e.whitney} that
\begin{eqnarray*}
\ell(S) \quad \le \quad L_2^{-1}\dist((\xi_1,\xi_2),\ell)  
&\lesssim& L_2^{-1} |\xi_1-\xi_2| \quad \lesssim \quad (L_1/L_2)|L|  \ \ .
\end{eqnarray*}
(Recall that $L_1=O(L_2)$.) This completes the proof of (i).

(ii) Let $A$ be the distance from $(\xi_1,\xi_2)$ to the set of two corners $(M,M)$ and $(N,N)$. Then the desired lower bound for $A$ follows from
\begin{eqnarray*}
A &\gtrsim& \dist((\xi_1,\xi_2),\ell) \quad \gtrsim \quad L_1\ell(S) \ \ .
\end{eqnarray*}
Using the triangle inequality and the Whitney property \eqref{e.whitney}, we also have
\begin{eqnarray*}
A &\le&   \dist((\xi_1,\xi_2),\ell) + \dist((\xi_1+\xi_2)/2, \{M,N\})\\
&\lesssim& L_2 \ell(S) + L_1 |L|  \ \ .
\end{eqnarray*}
Thanks to (i) we obtain  the upper bound $A\lesssim L_1 \ell(S)$, as desired. 
\endproof \\

\noindent \underline{Part III (final step).}
We now examine $m_{LM}$. Using \eqref{e.lmcc} and Lemma~\ref{l.lengthIJ} and Lemma~\ref{l.lengthSL}, it follows that $m_{LM}(M,N,\xi_1,\xi_2)$ could be written as
\begin{eqnarray*}
&=& \chi_{\xi_1<\xi_2} \chi_{M<\xi_1,\xi_2<N} \Big(\sum_{\alpha,\beta}\sum_{I,J} \phi_{\alpha,I}(\xi_1)\phi_{\beta,J}(\xi_2)1_{constraints \,\, on\,\, I,J}\Big) \times \\
& \times&\Big(\sum_{k,\gamma}a_k\sum_{S,L}  \phi_{1,S,k}(\xi_1) \phi_{2,S,k}(\xi_2) \phi_{3,S,k}(\xi_1+\xi_2)\phi_{\gamma,L}(\xi_1+\xi_2) 1_{constraints \,\, on \,\, S, L} \Big)
\end{eqnarray*}
here there are constraints on $I,J,S,L$ relative to $\alpha,\beta,\gamma$. Observe that the  summation over $(k,\gamma,S,L)$ is zero outside $\{\xi_1 < \xi_2\}$, while the  summation over $(\alpha,\beta, I,J)$ is zero outside $\{M<\xi_1,\xi_2<N\}$. Therefore we could drop the factor $\chi_{\xi_1<\xi_2}\chi_{M<\xi_1,\xi_2<N}$ in the right hand side and obtain
\begin{eqnarray*}
m_{LM} &=& \sum_{\alpha,\beta,\gamma,k} a_k \sum_{I \in \I_{M,N}(\alpha)} \sum_{J \in \I_{M,N}(\beta)} \sum_{L\in \I_{2M,2N}(\gamma)} \sum_{S \in \S}1_{\text{constraints on I,J,S,L}} \times\\
&\times& \Big[\phi_{\alpha,I}(\xi_1)\phi_{1,S,k}(\xi_1)\Big] \times \Big[\phi_{\beta,J}(\xi_2)  \phi_{2,S,k}(\xi_2) \Big]\times [\phi_{\gamma,L}(\xi_1+\xi_2)\phi_{3,S,k}(\xi_1+\xi_2)]  \ \ .
\end{eqnarray*}

 Since $a_k$ decays rapidly, for the purpose of proving Lemma~\ref{l.LM-properties} we may drop the summation over $k$ and consider only the contribution of one $k$.

Recall that $D_y$ denotes the dilation by $2^y$ with respect to $0$, i.e. $D_y \xi = 2^y \xi$. In particular $D_{-1}L = \{\frac 1 2 x: x\in L\}$. We claim that in any non-zero summand in $m_{LM}$, it holds that

(i) $|I|\sim |J|\sim \ell(S) \sim |L|$.

(ii) The spatial distances between any two of $I, J, S_1, S_2, D_{-1}(S_3), D_{-1}(L)$ are bounded above by $O(L_1|I|)$.

Indeed, (i) follows from  Lemma~\ref{l.lengthIJ} and Lemma~\ref{l.lengthSL}:
\begin{eqnarray*}\ell(S) &\sim |L| &\sim \frac 1 {L_1}\dist((\xi_1,\xi_2), \{(M,M),(N,N)\})\\
|I| &\sim |J| &\sim  \frac 1{L_1}\dist((\xi_1,\xi_2), \{(M,M),(N,N)\}) \ \ .
\end{eqnarray*}

For (ii),  by the Whitney property \eqref{e.whitney}   $\dist(S_1, S_2)=O(L_2\ell(S))$, thus the distances between $D_{-1}(S_3)$ and $S_1,S_2$ are  $O(L_1\ell(S))$. The desired claim now follows from examining the factors in the summation.

It follows that by decomposing $m_{LM}$ into $O(1)$ sums we may assume that $J,L, S_1, S_2, S_3$ are completely determined from $I$: comparable length and nearby location. 

Since $a_k$ decays rapidly we can ignore the summation over $k$, and below we will even drop the dependence on $k$ of the inner sums for brevity of notations. We end up with a symbol of the form
$$\sum_{I\in \D_1}\sum_{J\in \D_2}\sum_{L\in \D_3} \psi_{1,I}(\xi_1) \psi_{2,J}(\xi_2) \psi_{3,L}(\xi_1+\xi_2) 1_{constraints}$$
here $\D_1,\D_2,\D_3$ could be $\D_{left}$ or $\D_{right}$.are three fixed collections of (standard) dyadic intervals, $\psi_{1,I}$ is a $C^n$ bump function adapted to $I$ and is supported on $\frac 54 I$, and $\psi_{2,J}$ and $\psi_{3,L}$ satisfy similar properties.  The contraints read as follows: for fixed bounded integers $m_1,n_1,m_2,n_2$ it holds that
\begin{itemize}
\item $|J|=2^{m_1}|I|$ and $c(I)+n_1|I| \in J \ne\emptyset$.
\item $|L|=2^{m_2}|I|$ and $c(I)+n_2|I| \in D_{-1}(L)$.
\item $M\in I_{lower,\alpha} \cap J_{lower,\beta}$ and $2M \in  L_{lower,\gamma}$;
\item $N\in I_{upper,\alpha} \cap J_{upper,\beta}$ and $2N \in L_{upper,\gamma}$.
\end{itemize}
 This completes the proof of Lemma~\ref{l.LM-properties}. \endproof

\subsubsection{Reduction to $T_{LM}$,  part II: completion of the proof.} Using Lemma~\ref{l.LM-properties}, we may decompose $m_{LM}$ into boundedly many $m_{m_1,n_1,m_2,n_2}$, defined by
\begin{eqnarray*}
&&\sum_{\alpha,\beta,\gamma}\sum_{I, J, L}  \psi_{1,I}(\xi_1) \psi_{2,J}(\xi_2) \psi_{3,L}(\xi_1 + \xi_2) 1_{constraints \,\, on \,\, I, J, L} \times\\
&\times& 1_{\{M\in I_{lower,\alpha}, N\in I_{upper,\alpha}\}} 1_{\{M\in J_{lower,\beta}, N\in J_{upper,\beta}\}} 1_{\{2M\in L_{lower,\gamma}, 2N\in L_{upper,\gamma}\}} \ \ ,
\end{eqnarray*}
and the `constraints' on $I,J,L$ specify the location and the length of $J$ and $L$ relative to $I$ using $m_1,n_1,m_2,n_2$ as discussed in the last section. Note that the decomposition of $m_{LM}$ is independent of $M$ and $N$.

For each fixed $I,J,L$, we will   sum the summands
over $\alpha$, $\beta$, $\gamma$. Using Lemma~\ref{l.ob3}, we will  divide $m_{m_1,n_1,m_2,n_2}$ into $8$ symbols, such that in the summation  each of $I,J,L$ are required to be in one of $\D_{left}, \D_{right}$, and each of these $8$ symbols could be further decomposed into  $O(L_1)$ symbols having the following structure:
\begin{eqnarray*}
\sum_{I \in \D_1}\sum_{ J\in \D_2} \sum_{L\in \D_3} && \psi_{1,I}(\xi_1) \psi_{2,J}(\xi_2) \psi_{3,L}(-\xi_1-\xi_2) \times\\
&\times& 1_{\{M\in I_{1,lower}, N\in I_{1,upper}\}} 1_{\{M\in J_{2,lower}, N\in J_{2,upper}\}} 1_{\{2M\in L_{3,lower},  2N\in L_{3,upper}\}}
\end{eqnarray*}
where $\psi_j$ are bump functions supported in $[-c, c]$,  ($1/2<c<5/8$) and the collection of interval triples $\{
(I_{1,lower}, I, I_{1,upper}), I\in \D_1\}$ is rigid, and the same holds for the other two collections of interval triples. Note that the rigidity type of these collections are not opposite, i.e. we won't have both $(infinite, finite)$ and $(finite, infinite)$. Due to the flexibility of $(finite,finite)$ rigidity (which could be converted into two terms of any of the other types), we could assume that the rigidity type of the three collections are the same.

We now split the intervals $I,J,L$ if necessary to ensure $|I|=|J|=|L|$: to do this we will split  $\psi_{1,I}$, $\psi_{2,J}$, $\psi_{3,L}$  and we will correspondingly split all bounded intervals in $I_{1,lower}$, $I_{1,upper}$, $J_{2,lower}$, $J_{2,upper}$, $L_{3,lower}$, $L_{3,upper}$ (but we won't split the halflines). Note that  if we split $\psi_{1,I}$ into $2^k$ bump function adapted to the corresponding subintervals of $I$ ($k\ge 0$),  then each new bump function in theory could be supported in an interval as large as the $(1+2^k (c-1/2))$ enlargement of  the adapting subinterval. Since $k$ is bounded, by choosing $c$ sufficiently close to $1/2$ when   partitioning $1_{M<\xi<N}$  we could ensure that $1+2^k(c-1/2) < 5/4$.

After the splitting, the rest of the discretization is similar to  the discretization of $m_{BC}$. We omit the details.

\end{document}